\documentclass[review,17pt,reqno]{elsarticle}
\usepackage{amssymb,amsmath}
\usepackage{cases}
\usepackage{exscale}
\usepackage{relsize}
\usepackage{subfigure}
\usepackage{natbib}
\usepackage[colorlinks=true,linkcolor=blue,citecolor=blue,urlcolor=blue,]{hyperref}
\usepackage{tabularray}
\usepackage{caption}
\usepackage{booktabs}
\usepackage{makecell}
\usepackage{graphicx}
\usepackage{pythonhighlight}
\usepackage{array, threeparttable}
\usepackage{lineno}

\usepackage{algorithmicx}
\usepackage[ruled]{algorithm2e}    
\usepackage{algpseudocode}
\usepackage{mathrsfs} 
\usepackage{arydshln}
\usepackage{multirow}
\usepackage{graphicx} 

\usepackage{ulem}

\usepackage{epsfig}                          
\usepackage{epstopdf}      
\usepackage{placeins}                  
\usepackage{appendix}

\begin{document}
\newtheorem{mytheorem}{Theorem}[section]
\newtheorem{myproposition}{Proposition}[section]
\newtheorem{mylemma}{Lemma}[section]
\newtheorem{mycorollary}[mytheorem]{Corollary}
\newtheorem{myprop}[myproposition]{Proposition}
\newtheorem{myremark}{Remark}[section]
\newtheorem{example}{Example}[section]
\renewcommand{\baselinestretch}{1.1}
\newcommand{\figref}[1]{Fig.~\ref{#1}}
\makeatletter
\newcommand\figcaption{\def\@captype{figure}\caption}
\newcommand\tabcaption{\def\@captype{table}\caption}
\makeatother
\begin{frontmatter}
\title{Variable-preconditioned transformed primal--dual method for generalized Wasserstein gradient flows}

\author[label1]{Jin Zeng\corref{cor1}}
\ead{202511110511@std.uestc.edu.cn}

\author[label1]{Dawei Zhan\corref{cor1}}
\ead{2021110903008@std.uestc.edu.cn}

\author[label2]{Ruchi Guo}
\ead{ruchiguo@scu.edu.cn}

\author[label1]{Chaozhen Wei\corref{cor2}}
\ead{cwei4@uestc.edu.cn}

\address[label1]{School of Mathematical Sciences, University of Electronic Science and Technology of China, Chengdu, Sichuan 611731, China}
\address[label2]{School of Mathematics, Sichuan University, Chengdu, Sichuan 610065, China}
\cortext[cor1]{J. Z. and D. Z. contributed equally to this work.}
\cortext[cor2]{Corresponding author.}


\begin{abstract}
We propose a Variable-Preconditioned Transformed primal--dual (VPTPD) method for solving generalized Wasserstein gradient flows based on the structure-preserving JKO scheme. This is a nontrivial extension of the TPD method [Chen et al. (2025) SIAM J. Sci. Comput.] incorporating proximal splitting techniques to address the challenges arising from the nonsmoothness of the objective function.  Our key contributions include: (i) a semi-implicit-explicit iterative scheme that combines proximal gradient steps with explicit gradient steps to treat the nonsmooth and smooth terms respectively; (ii) variable-dependent preconditioners constructed from the Hessian of a regularized objective to balance iteration count and per-iteration cost; (iii) a proof of existence and uniqueness of bounded solutions for the generalized proximal operator with the chosen preconditioner, along with a convergent and bound-preserving Newton solver; and (iv) an adaptive step-size strategy to improve robustness and accelerate convergence under poor Lipschitz conditions of the energy derivative. Comprehensive numerical experiments spanning from 1D to 3D settings demonstrate that our method achieves superior computational efficiency---achieving up to a 20$\times$ speedup over existing methods---thereby highlighting its broad applicability through several challenging simulations. 
\end{abstract}

\end{frontmatter}


\section{Introduction}
Consider gradient flows with concentration-dependent mobility over a subdomain $\Omega \subset \mathbb{R}^{d}$ of the form
\begin{equation}\label{Wasserstein gradient flow}
	\partial_{t}\rho = -\nabla\cdot\left(M\left(\rho\right)\boldsymbol{v}(\rho)\right),\quad 
	v(\rho) = -\nabla\displaystyle{\frac{\delta}{\delta\rho}\mathcal{J}\left(\rho\right)}, \quad  \left(x,t\right)\in\Omega\times\left(0,T\right),
\end{equation}
where $\rho\left(x, t\right)$ is the density function, $M\left(\rho\right)$ is the concentration-dependent mobility, $\boldsymbol{v}(\rho)$ is the velocity field, and $\delta \mathcal{J}/\delta\rho$ is the chemical potential defined as the first variation of $\mathcal{J}(\rho)$. The energy functional $\mathcal{J}$, in general, may include but not limited to the following terms
\begin{equation}
	\label{energy_function}
	\mathcal{J} = \int_{\Omega}\Big(\mathcal{H}\left(\rho\right) + \mathcal{V}(x)\rho\Big)\,\mathrm{d}x + \frac{\epsilon^{2}}{2}\int_{\Omega} \left\vert\nabla\rho\right\vert^{2}\, \mathrm{d}x + \displaystyle\frac{1}{2}\iint_{\Omega^{2}} \mathcal{W}\left(x - y\right) \rho\left(x\right) \rho\left(y\right)\,\mathrm{d}x\,\mathrm{d}y,
\end{equation}
where the first term represents the internal and potential energy, the second term represents the Dirichlet energy, and the third part represents the nonlocal interaction energy. 

The mobility $M(\rho)$ is often degenerate for specific $\rho$ values. For example, $M(\rho)=\rho$ has a degeneracy at $\rho=0$ in aggregation-drift-diffusion equations \cite{carrillo2015finite}; $M(\rho)=1-\rho^2$ is degenerate at $\rho=\pm 1$ in surface diffusion flows \cite{JiangWei2019SolidDewetting}. Due to this degeneracy, the solution to \eqref{Wasserstein gradient flow} is shown to have global bounds \cite{Ambrosio2005Gradient,Elliott1996CahnHilliard}. Together with no-flux boundary conditions, the admissible set $\mathcal{M}(\Omega)$ of the solution satisfies:
\begin{equation}\label{measure space}
	\mathcal{M}\left(\Omega\right) = \left\{\rho: \rho\; \text{measurable},\; M\left(\rho\right) \geq 0 \;\;\text{and} \int_{\Omega}\rho\left(x\right)\,\mathrm{d}x = \text{const.}\right\}.
\end{equation}  

The coupling between the nonlinear, nonlocal energy functional and the degenerate concentration-dependent mobility poses significant challenges to the computation of such gradient flow in \eqref{Wasserstein gradient flow}. To overcome the stability issue stemming from the strong nonlinearity and high-order derivatives in the chemical potential, a popular line of research concerns the development of energy stable methods, such as the convex splitting method \cite{eyre1998unconditionally, baskaran2013convergence}, the stabilized semi-implicit method \cite{xu2006stability}, the Invariant Energy Quadratization (IEQ) method \cite{Yang2020IEQAllenCahn}, the Scalar Auxiliary Variable (SAV) method \cite{shen2018scalar}, and the exponential time differencing (ETD) method \cite{ju2015fast,hou2024energy}. 
Bound-preserving methods that maintain the global bounds of discrete solutions also attract intense interest and have been developed, including but are not limited to the flux-limiting discontinuous Galerkin method \cite{frank2020bound} and the Lagrangian multiplier method \cite{cheng2022new1,cheng2022new2}. More recently, some improved variants of SAV and ETD methods are proposed to achieve bound preservation while maintaining energy stability \cite{huang2021bound, huang2023structure,huang2025maximum}. However, these methods require intricate design for specific energy functionals and mobilities to simultaneously address energy dissipation, bound preservation and mass conservation for \eqref{Wasserstein gradient flow} at the discrete level, especially for nonlocal energy and nonlinear mobilities.

Recently, a novel structure-preserving approach based on minimizing movements \cite{DeGiorgi} for \eqref{Wasserstein gradient flow} with respect to the generalized Wasserstein distances has been developed for linear \cite{carrillo2022primal}, concave quadratic \cite{carrillo2024structure} and general nonlinear mobilities \cite{deng2025PDFB}. For the generalized Wasserstein gradient flows \eqref{Wasserstein gradient flow}, the celebrated JKO scheme \cite{jordan1998variational} approximates the solution $\rho^{k}$ at $t = k\tau$ by solving the following minimization
\begin{equation}\label{JKO}
	\rho^{k+1}\in \arg\inf\limits_{\rho} \left\{{\dfrac{1}{2 \tau}}d_{\mathcal{W}_m}^2\left(\rho, \rho^k \right)+\mathcal{J}\left(\rho\right)\right\}
\end{equation}
where $d_{\mathcal{W}_m}$ is the generalized Wasserstein distance induced by the mobility $M(\rho)$. It has been proved that the minimizer of \eqref{JKO} always exists for suitable step size $\tau$ \cite{Ambrosio2005Gradient}, and the temporal discrete sequence $\{\rho^k\}$ weakly converges to the solution of \eqref{Wasserstein gradient flow} with linear \cite{jordan1998variational} and concave mobilities \cite{Carrillo2010NonlinearMobility,lisini2012cahn} as $\tau\rightarrow 0$. By utilizing the Benamou-Brenier dynamic formulation of transport distances \cite{benamou2000computational,Dolbeault2009NewTransport}, the scheme \eqref{JKO} can be further reframed as a convex optimization problem subject to linear PDE constraints. 
Specifically, by introducing the transport dynamics with an auxiliary time variable $s\in\left[0, 1\right]$, the squared distance between two measures $\rho\left(\cdot, 0\right)$ and $\rho\left(\cdot, 1\right)$ in the Wasserstein space is defined as the minimal total kinetic energy $\Phi$ along the optimal transport geodesic connecting them
\begin{align}\label{W_distance}
	\Phi\left(\rho, \boldsymbol{m}\right) = \int_{0}^{1}\int_{\Omega}\varphi\left(\rho,\boldsymbol{m}\right)\,\mathrm{d}x\,\text{d}s,\quad
	\varphi\left(\rho,\boldsymbol{m}\right) = \begin{cases}
		\displaystyle\frac{{\left\Vert \boldsymbol{m}\right\Vert}^{2}}{M\left(\rho\right)} ,& \mathrm{if}\; M > 0,\\
		0,& \mathrm{if}\; (M,\boldsymbol{m}) = (0,\boldsymbol{0}),\\
		+\infty,&\mathrm{otherwise},
	\end{cases}
\end{align}
where $\boldsymbol{m}(x,s)$ is the momentum variable that controls the transport geodesic, i.e., $\rho$ and $\boldsymbol{m}$ satisfy the continuity equation. 
Then the dynamic JKO scheme \cite{carrillo2022primal,carrillo2024structure} approximates $\rho^{k+1}(x) = \rho(x, 1)$ by the constrained minimization 
\begin{align}\label{dynamic-JKO}
	\begin{cases}
		$$\left(\rho\left(x,s\right),\boldsymbol{m}\left(x,s\right)\right)\in \arg\inf\limits_{\left(\rho,\boldsymbol{m}\right)} \mathcal{F}(\rho,\boldsymbol{m}):=\dfrac{1}{2 \tau}\Phi\left(\rho, \boldsymbol{m} \right)+\mathcal{J}\left(\rho\left(\cdot,1\right)\right)$$,\\
		$$\mathrm{s.t.}\quad \partial_{s}\rho+\nabla\cdot \boldsymbol{m}=0,\;\rho\left(x,0\right)=\rho^{k}\left(x\right),\; \boldsymbol{m} \cdot \boldsymbol{n} = 0$$,
	\end{cases}
\end{align}
where $\boldsymbol{n}$ represents the unit outer normal vector. This scheme \eqref{dynamic-JKO} naturally inherits several key properties: original energy dissipation from the underlying variational structure, mass conservation from the no-flux boundary condition, and bound preservation guaranteed by the definition of the action function $\varphi$ in \eqref{W_distance} \cite{carrillo2024structure}. 

Convex optimization techniques have been widely used for the resultant linearly constrained convex minimization of JKO scheme \cite{benamou2016augmented,papadakis2014optimal}. 
However, the main computational challenge arises from the non-smoothness of $\varphi(\rho,m)$. A feasible approach is to introduce the Fisher information regularization \cite{LiWuchen2020Fisher} in the objective function to overcome the lack of smoothness and strict convexity, which gives access to the sequential quadratic programming method and enforces the boundedness of solution, yet introduces artificial numerical diffusion at the PDE level. Instead, Carrillo et al. \cite{carrillo2022primal} introduces a penalty term in the objective function to handle constraints and employs a primal--dual three-operators (PD3O) splitting method \cite{Yan2018NewPDThreeFunctions} to solve one-step JKO scheme (PDJKO), where the proximal operator of $\varphi$ is exactly solved and hence the boundedness of $\rho$ is truly preserved. While the PDJKO method performs well for $M(\rho)=\rho$, it may suffer from slow convergence for nonlinear mobility. Therefore, they \cite{carrillo2024structure} further proposed a preconditioned primal dual method (PrePDJKO) by introducing a matrix preconditioner for dual variable \cite{liu2021acceleration} and combining with FFT-based efficient matrix inversion, which greatly improves the computational efficiency compared to the PDJKO method. However, large-scale simulation (e.g., 3D simulation or 2D simulation with finer grid resolution) is still challenging given that the number of required iterations for convergence can drastically increase with the size of problem. 
Recently, the G-prox PDHG method \cite{Jacobs2019Gprox} achieves grid-resolution-independent convergence for nonsmooth optimization by modifying the primal update with a specific preconditioning. However, a direct application of this approach may be impractical for our problem \eqref{dynamic-JKO} with general energy functional involving differential operators or nonlocal convolutions, as the resulting generalized proximal operator can be computationally inefficient.
Notably, a recent PDFB splitting method \cite{deng2025PDFB} for computing the JKO scheme \eqref{dynamic-JKO} can also obtain a grid-resolution-independent convergence rate, yet it involves more complex sub-problems that may degrade the applicability to large-scale simulations. Finally, we point out that there is also some literature on primal--dual splitting methods for generalized JKO scheme with higher-order numerical discretizations and their extension to reaction-diffusion systems and mean-field control dynamic problems \cite{Fu2023JCP,Fu2024JFM}.


Roughly speaking, all the approaches above still belong to the primal--dual (PD) family \cite{chambolle2011first}. 
Beyond this line of work, the newly proposed Transformed primal--dual (TPD) method \cite{chen2023transformed1,chen2023transformed2,chen2025novel} accelerates convergence by recasting the saddle-point system and its underlying dynamics.
The key idea is to construct a suitable preconditioner that transforms the saddle-point system into an equivalent form that is nearly upper triangular and diagonally dominant, thereby almost decoupling the primal and dual variables.
However, the current TPD method is merely designed for smooth and strongly-convex objective functional, not applicable to the present case.

In this paper, we propose a Variable-Preconditioned Transformed primal--dual method (VPTPD) for solving generalized Wasserstein gradient flows \eqref{Wasserstein gradient flow}. The VPTPD is a nontrivial extension of the TPD method to handle the objective functions in \eqref{dynamic-JKO} comprising a convex nonsmooth term $\Phi$ and a smooth differentiable term $\mathcal{J}$. The novelty is multifold.  
First, we integrate the TPD method with proximal splitting method and obtain a semi-implicit-explicit iterative scheme, which treats the nonsmooth term with a proximal gradient step and the smooth term with an explicit gradient step. Second, we construct variable-dependent preconditioners through the Hessian of a regularized smooth objective function, tailored to obtain a delicate balance between the number of required iterations for convergence and the computational complexity of each iteration. 
Third, we prove the existence and uniqueness of bounded solution for the generalized proximal operator with the chosen preconditioner, and propose a Newton solver which is both convergent and bound-preserving. 
At last, we develop an adaptive step-size strategy for the TPD method to mitigate performance degradation caused by poor Lipschitz continuity of the energy derivative, and accelerate its global convergence rate. 
Through comprehensive numerical experiments spanning from 1D to 3D settings, we demonstrate superior computational efficiency of our VPTPD method---achieving up to a 20$\times$ speedup over existing methods. 

The rest of our paper is organized as follows. In \autoref{Fully_discrete_dynamic_JKO_scheme}, we introduce the fully discrete dynamic JKO scheme. \ref{TPD} introduces the TPD method and compares it with the PDJKO and PrePDJKO methods. \ref{num method} elaborates on the VPTPD method adapting to the JKO scheme, in terms of the iterative scheme with tailored preconditioners, computation of proximal operators and adaptive strategy. In \autoref{experiment}, a suite of numerical examples are reported to demonstrate the high efficiency and applicability of the VPTPD method to a number of challenging simulations. 

\section{Fully discrete dynamic JKO scheme}\label{Fully_discrete_dynamic_JKO_scheme}

We now consider the discretization of \eqref{dynamic-JKO} over a computational domain of a $d$-dimensional box $\Omega = \mathop{\times}_{j=1}^{d}\left[l_j,r_j\right]\allowbreak\in\mathbb{R}^{d}$. The domain $\Omega$ is partitioned into $N = \prod_{j=1}^{d} N_j$ uniform grid cells, with a grid spacing $\Delta x_j = \frac{r_j-l_j}{N_j}$ in the $j$-th dimension and cell volume $\Delta V = \Pi_{j=1}^{d} \Delta x_{j}$. Let $(\rho_{\boldsymbol{i}},\boldsymbol{m}_{\boldsymbol{i}})$ represent the discrete values $(\rho(x_{\boldsymbol{i}}),\boldsymbol{m}(x_{\boldsymbol{i}}))$ at the center of the grid cell $\Omega_{\boldsymbol{i}}$, where $\boldsymbol{i} = [i_{1},\cdots,i_{d}]$ is the index vector ($1\leq i_{j}\leq N_{j},\,j = 1,\cdots,d$) and $x_{\boldsymbol{i}}=\{x_j:x_j=l_j+(i_j-1/2)\Delta x_j\}_{j=1}^{d}$.  

For the constraints of contuity equation with initial and boundary conditions in \eqref{dynamic-JKO}, we apply the central difference scheme for the divergence term $\nabla\cdot\boldsymbol{m}$ and a one-step finite difference in the auxiliary time $s$. This discretization preserves the first-order accuracy of the overall JKO scheme \cite{LiWuchen2020Fisher} and results in the following system:
\begin{align} \label{discrete constraint}
	\begin{cases}
		\rho + \sum_{j=1}^{d} \dfrac{\boldsymbol{m}_{\boldsymbol{i}+ e_{j}} - \boldsymbol{m}_{\boldsymbol{i}-e_{j}}}{2\Delta x_{j}} = \rho^{k},\\
		\boldsymbol{m}_{\boldsymbol{i}-\boldsymbol{e}_{j}} = -\boldsymbol{m}_{\boldsymbol{i}} \quad \mathrm{if} \; i_{j} = 1,\quad \mathrm{and} \quad 
		\boldsymbol{m}_{\boldsymbol{i}+\boldsymbol{e}_{j}} = -\boldsymbol{m}_{\boldsymbol{i}} \quad \mathrm{if} \; i_{j} = N_{j},
	\end{cases}
\end{align}
where $\boldsymbol{e}_{j}$ denotes the basis vector in the $j$-th dimension. The Wasserstein distance functional \eqref{W_distance} is approximated by using the midpoint rule in space and a one-step approximation in the auxiliary time ($\Delta s = 1$), leading to its discrete form:
\begin{equation}\label{Phi_discrete}
	\Phi_{h}\left(u\right) = \sum_{\boldsymbol{i}}\varphi\left(\rho_{\boldsymbol{i}},\boldsymbol{m}_{\boldsymbol{i}}\right)\Delta V
\end{equation}  
where we define the variable $u := [\{\rho_{\boldsymbol{i}}\}, \{\boldsymbol{m}_{\boldsymbol{i}}\}]$. Denoting $\mathcal{J}_h(u)$ the discrete energy functional that approximates $\mathcal{J}(\rho)$ by an appropriate second-order quadrature rule, we obtain the fully discrete dynamic JKO scheme:	
\begin{align}\label{generalized discrete JKO}
	\begin{cases}
		u^{k+1}\in\arg\min\limits_{u}\left\{\dfrac{1}{2\tau}
		\Phi_{h}\left(u\right) + \mathcal{J}_{h}\left(u\right)\right\},\\
		\mathrm{s.t.} \quad Bu = b,
	\end{cases}
\end{align} 
where the discretization of constraints \eqref{discrete constraint} is written in a matrix-vector form with $B$ being a linear surjective operator.

Given that the discretization of the scheme and its convergence has been elaborated elsewhere \cite{carrillo2022primal,carrillo2024structure}, which is not the focus of this paper, we refer the readers to the aforementioned previous works for details. The proposed fully discrete dynamic JKO scheme has first-order accuracy in time and second-order accuracy in space, which can be validated by numerical experiments in \ref{JKO_accuracy}. 

\section{Transformed Primal Dual (TPD) method}\label{TPD}
In this section, we first briefly introduce the TPD method and then conduct a comparative discussion on the TPD method with the PD-type methods.

\subsection{Review of the TPD method}\label{TPD-method}
The TPD method \cite{chen2023transformed1,chen2023transformed2} tackles the linearly constrained optimization with a smooth, convex objective function $f(u)$ in the general form
\begin{equation}\label{l-c problem}
\min\limits_{u} f\left(u\right), \quad
\mathrm{s.t.}\ Bu=b,
\end{equation}
by introducing a Lagrangian multiplier $p$ and considering the saddle point problem
\begin{equation}\label{minimax problem}
	\min\limits_{u}\max\limits_{p} \ \mathcal{L}(z) =f\left(u\right) + \left\langle Bu-b,p \right\rangle, \quad z = (u,p).
\end{equation}
The saddle point $z^*=\left(u^{*}, p^{*}\right)$ of the above Lagrangian satisfies the first-order optimal condition:
\begin{equation}\label{PD_saddle}
	\mathcal{K}(z^*) :=
	\begin{bmatrix}
		\nabla f & B^{\top}\\
		B & 0
	\end{bmatrix}
	\begin{bmatrix}
		u^{*}\\
		p^{*}
	\end{bmatrix}=
	\begin{bmatrix}
		\mathbf{0}\\
		b
	\end{bmatrix},
\end{equation}
where the gradient operator is understood as $\langle \nabla f, u\rangle := \nabla f(u)$. 

The classical PD-type methods searches the saddle point via the following primal--dual flow with certain discretization approaches
\begin{equation}\label{PD_flow}
	\begin{bmatrix}
		u'\\
		p'
	\end{bmatrix} = 
	\begin{bmatrix}
		-\nabla f & -B^{\top}\\
		B & \mathbf{0}
	\end{bmatrix}
	\begin{bmatrix}
		u\\
		p
	\end{bmatrix}-
	\begin{bmatrix}
		\mathbf{0}\\
		b
	\end{bmatrix}.
\end{equation}
However, the coefficient matrix, i.e., the operator $\mathcal{K}$, has a vanishing lower-right block, and consequently lacks diagonal dominance and strong monotonicity in general; this typically leads to slow trajectories and poor conditioning in the discrete iterations.

The TPD method alleviates this difficulty by employing the following Schur-complement-based decomposition for $\mathcal{K}$:
\begin{equation}\label{eq_decomp0}
\begin{bmatrix}
\nabla f & B^{\top}\\
B & \mathbf{0}
\end{bmatrix}
=
\underbrace{\begin{bmatrix}
I & \mathbf{0}\\
B I_u^{-1} & I
\end{bmatrix}}_{=:L^{-1}}
\underbrace{\begin{bmatrix}
\nabla f & B^{\top}\\
B\!\left(I - I_u^{-1}\nabla f\right) & -B I_u^{-1} B^{\top}
\end{bmatrix}}_{=:\mathcal{K}_{\mathrm{TPD}}}.
\end{equation}
where $I_{u}$ is a linear and symmetric positive definite (SPD) operator as the approximation of $\nabla f$ such that its inverse $I_{u}^{-1}$ and the Schur complement $BI_{u}^{-1}B^{\top}$ are easily computed. 
By applying the \textit{transformation} matrix $L$ to the original system \eqref{PD_saddle},
the TPD method basically solves the following equivalent transformed system:
\begin{equation}\label{TPD_saddle}
	\mathcal{K}_{\mathrm{TPD}}(z^*)=
	\begin{bmatrix}
		\nabla f & B^{\top}\\
		B\left(I - I_{u}^{-1}\nabla f\right) & -BI_{u}^{-1}B^{\top}
	\end{bmatrix}
	\begin{bmatrix}
		u^{*}\\
		p^{*}
	\end{bmatrix}=
	\begin{bmatrix}
		\mathbf{0}\\
		b
	\end{bmatrix}.
\end{equation}
Then, optimization algorithms can be derived by discretizing the following so called TPD flow:
\begin{equation}\label{TPD flow}
	\begin{bmatrix}
		u'\\
		p'
	\end{bmatrix} = 
	\begin{bmatrix}
		-\nabla f & -B^{\top}\\
		B\left(I - I_{u}^{-1}\nabla f\right) & -BI_{u}^{-1}B^{\top}
	\end{bmatrix}
	\begin{bmatrix}
		u\\
		p
	\end{bmatrix}-
	\begin{bmatrix}
		\mathbf{0}\\
		b
	\end{bmatrix}.
\end{equation}
Notice that the lower-left block of the operator $\mathcal{K}_{\mathrm{TPD}}$ will become small provided that $I_u$ is close enough to $\nabla f$, rendering $\mathcal{K}_{\mathrm{TPD}}$ nearly upper triangular. 
Therefore, this transformation (partially) decouples the primal and dual variables and adds strong convexity in the dual variable, which can significantly accelerate the convergence rate \cite{chen2023transformed1}. 
In particular, when $\nabla f$ is an invertible linear operator and $I_{u}^{-1}=(\nabla f)^{-1}$, the primal and dual variables are fully decoupled, i.e., $B\left(I - I_{u}^{-1}\nabla f\right) = 0$. 
To illustrate the superiority of TPD over PD method, we present the following toy example, where we will demonstrate the essential role of the transformation for fast convergence.
\begin{example}
Consider a dynamical system in $\mathbb{R}^2$:
\begin{equation}
\label{eq_toy}
\begin{bmatrix}
	u'\\
	p'
\end{bmatrix} = 
	\begin{bmatrix}
	-\epsilon & 1\\
	-1 & 0
\end{bmatrix}
\begin{bmatrix}
	u\\
	p
\end{bmatrix},
\end{equation}
where $\epsilon$ might be small for making the matrix non-diagonal-dominant.
In the following figures, we present the numerical results for the convergence of the PD flow \eqref{PD_flow} and TPD flow \eqref{TPD flow} both through explicit discretization.
\begin{figure}[htbp]
	\centering  
	\includegraphics[width=0.32\textwidth]{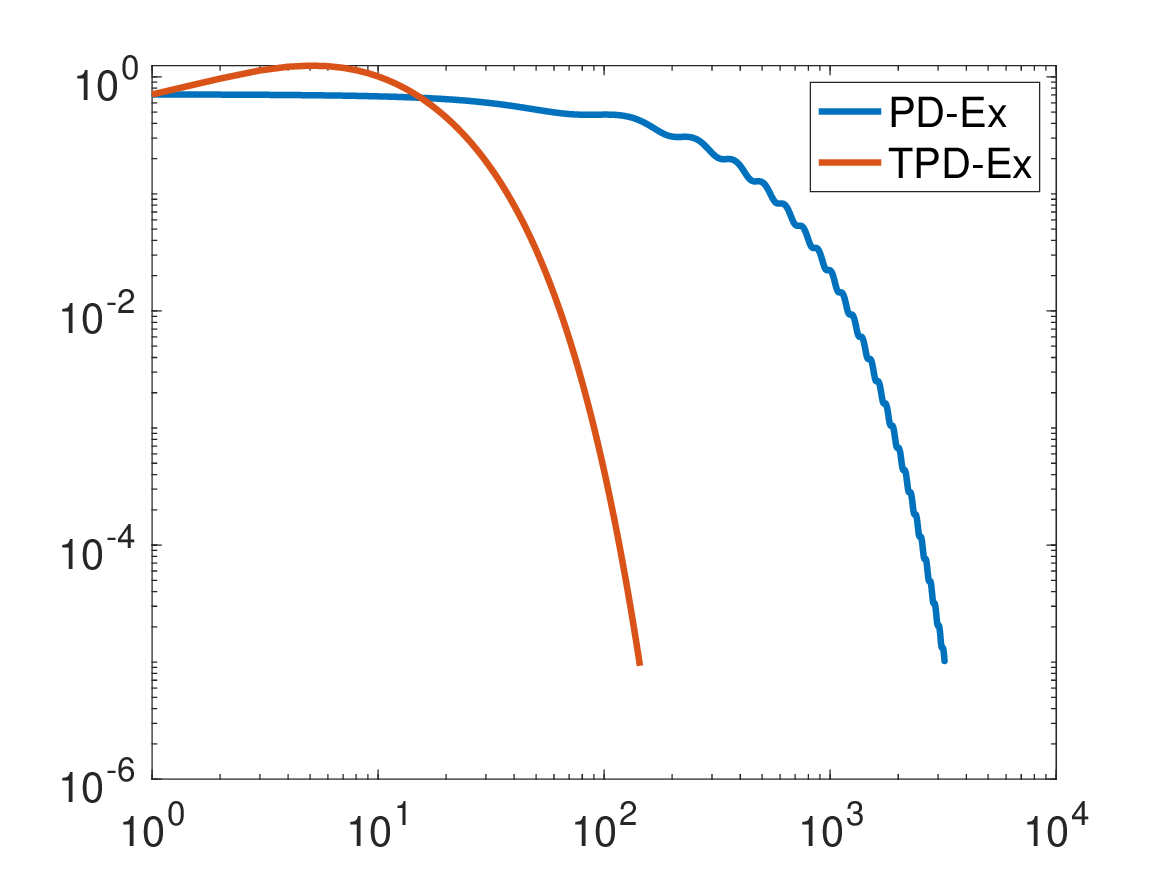}
	\includegraphics[width=0.32\textwidth]{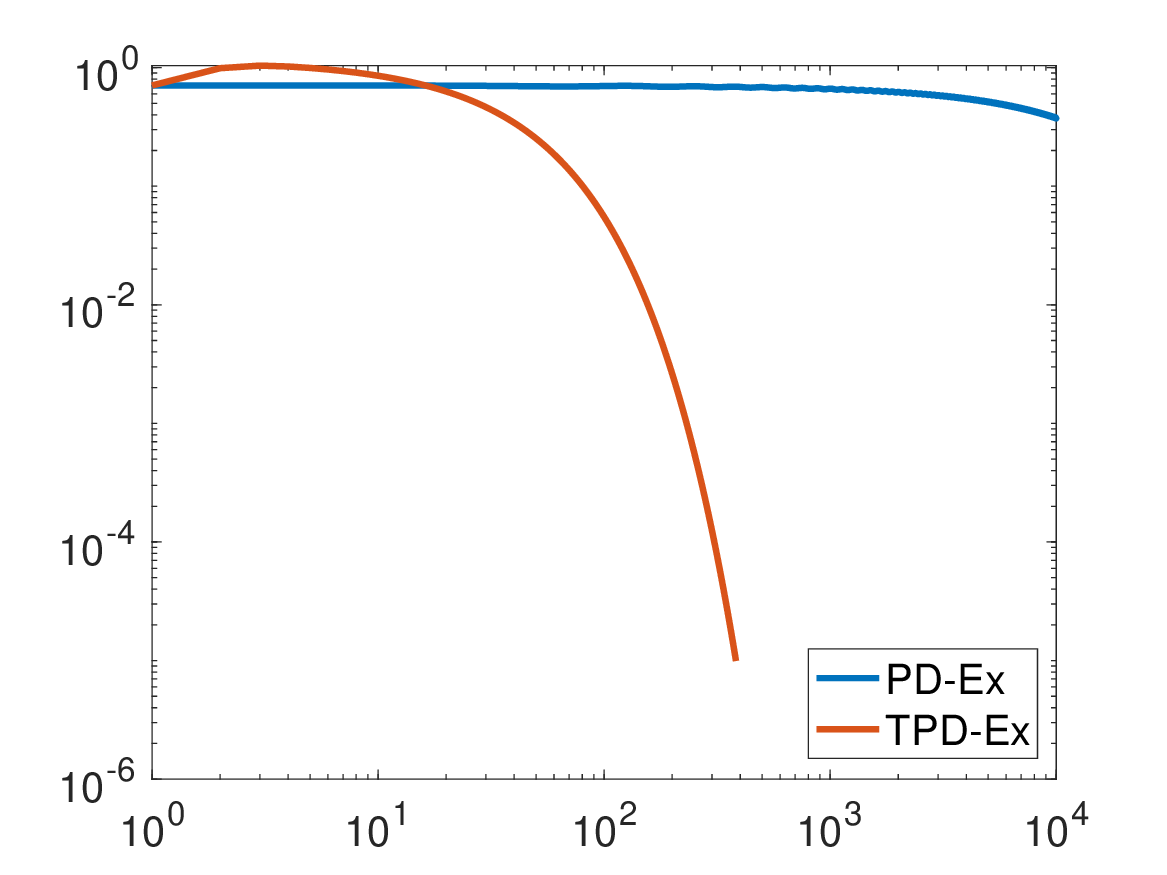}
	\caption{Evolution of solution errors: $\epsilon=0.3$ (left) and $\epsilon=0.03$ (right).}
	\label{fig:toy}
\end{figure}
\end{example}

When $f$ is highly nonlinear and $B$ involves with discrete differential operators, TPD further embeds preconditioners into flow, leading to the preconditioned TPD dynamics \cite{chen2023transformed2}:
\begin{equation}\label{TPD2_flow}
	\begin{bmatrix}
		u'\\
		p'
	\end{bmatrix} = 
	\begin{bmatrix}
		I_{u}^{-1} & 0\\
		0 & I_{p}^{-1}
	\end{bmatrix}
	\left(\begin{bmatrix}
		-\nabla f & -B^{\top}\\
		B\left(I - I_{u}^{-1}\nabla f\right) & -BI_{u}^{-1}B^{\top}
	\end{bmatrix}
	\begin{bmatrix}
		u\\
		p
	\end{bmatrix}-
	\begin{bmatrix}
		\mathbf{0}\\
		b
	\end{bmatrix}
	\right),
\end{equation}
where $I_{p}^{-1}$ is an appropriate SPD of $BI_{u}^{-1}B^{\top}$. In general, $I_u^{-1}$ and $I_{p}^{-1}$ should be variable-dependent (or time-dependent) and easily computable to achieve better efficiency. 
Applying explicit Euler discretization to the dynamics \eqref{TPD2_flow} with step sizes $(\lambda,\sigma)$ yields the basic iteration:
\begin{align}\label{TPD_discrete}
	\begin{cases}
		\dfrac{u_{n+1} - u_{n}}{\lambda}= -(I_{u})^{-1}\left(\nabla f\left(u_{n}\right)+B^{\top}p_{n}\right),\\
		\dfrac{p_{n+1} - p_{n}}{\sigma}= (I_{p})^{-1}\left(Bu_{n+1} - b -B(I_{u})^{-1}\left(\nabla f \left(u_{n+1}\right)+B^{\top}p_{n}\right)\right).
	\end{cases}
\end{align}


\subsection{Comparison with the PDHG-type methods}
When $f$ is nonsmooth and its explicit gradient is not available, the saddle point system is typically solved by primal--dual hybrid gradient (PDHG) method \cite{chambolle2011first} and its preconditioned variants \cite{liu2021acceleration,Jacobs2019Gprox}, which are essentially implicit Euler discretization of the PD flow \eqref{PD_flow}. The saddle point system \eqref{minimax problem} can be rewritten into the standard PDHG framework:
\begin{equation}\label{minimax problem2}
	\min\limits_{u}\max\limits_{p} \ \mathcal{L}(z) =f\left(u\right) + \langle Bu,p \rangle - g^{*}\left(p\right),
\end{equation}
where $g^{*}\left(p\right)=\langle b,p \rangle$ is the convex conjugate of the singleton indicator $g(x)=\delta_{\{b\}}$.
For solving \eqref{minimax problem2}, the proximal-point method with respect to a specific metric induced by an SPD matrix $P$ is given by:
\begin{equation}\label{prox_L}
z_{n+1} = \arg \min \limits_{u} \max \limits_{p} \ \mathcal{L}(z) + \dfrac{1}{2}\left\Vert z-z_{n}\right\Vert_{P}^{2}, \quad \text{where } \left\Vert z\right\Vert_{P}^{2}=z^{\top}Pz.
\end{equation}
PDHG-type methods with preconditioning \cite{liu2021acceleration,Jacobs2019Gprox} can be interpreted as operator-splitting realizations of the above proximal update with
block-diagonal $P_{\mathrm{PDHG}}=\mathrm{diag}\{\frac{1}{\lambda}I_u, \frac{1}{\sigma}I_p\}$ as preconditioners, yielding the iteration:
\begin{align}\label{PDHG}
	\begin{cases}
		& u_{n+1} = \arg \min \limits_{u} \  f(u) + \langle u,B^{\top} \bar{p}_{n} \rangle + \dfrac{1}{2\lambda}\left\Vert u-u_{n}\right\Vert_{I_u}^{2}, \\
		& p_{n+1} = \arg \min \limits_{p} \ g^{*}\left(p\right) - \langle Bu_{n+1},p \rangle + \dfrac{1}{2\sigma}\left\Vert p-p_{n}\right\Vert_{I_p}^{2}, 
	\end{cases}
\end{align}
where $\bar{p}_{n}$ is certain approximation of $p_{n+1}$ by extrapolation. For classical PDHG methods ($I_u=I$ and $I_p=I$), step sizes must satisfy the joint restriction $\lambda\sigma\|BB^{\top}\|<1$ for convergence, which causes severe difficulty when $B$ is a specific differential operator. In particular, G-prox PDHG \cite{Jacobs2019Gprox} performs preconditioning with $P_{\mathrm{Gprox}}=\mathrm{diag}\{\frac{1}{\lambda}B^{\top}B,\frac{1}{\sigma}I\}$ and achieves grid-resolution-independent convergence under certain assumptions on $B$ and $g^*$ (see \ref{Gprox_prox}).

Moreover, the implicit discretization of the preconditioned TPD flow \eqref{TPD2_flow} is equivalent to the proximal update \eqref{prox_L} with a block-triangular metric matrix:
\begin{equation}\label{P_VPTPD}
	 P_{\mathrm{TPD}} = L^{-1}\tilde{P} := \begin{bmatrix}
		I & 0\\
		-BI_u^{-1} & I
	\end{bmatrix}^{-1}
	\begin{bmatrix}
		\dfrac{1}{\lambda}I_u & \mathbf{0} \\
		\mathbf{0} & \dfrac{1}{\sigma}I_p
	\end{bmatrix} = 
	\begin{bmatrix}
		\dfrac{1}{\lambda}I_u & \mathbf{0}\\
		\dfrac{1}{\lambda}B & \dfrac{1}{\sigma}I_p
	\end{bmatrix},
\end{equation}
where $L$ is the transformation matrix defined in \eqref{eq_decomp0} that transforms the original saddle-point system to the nearly upper-triangular system while preserving the solution, and $\tilde{P}$ is a variable-dependent block-diagonal preconditioner that improves the conditioning of the transformed system. Notice that the equivalence here may be only understood in form, as the metric $P$ in \eqref{prox_L} is usually SPD. 
With implicit discretization, an iterative scheme is given by:  
\begin{align}\label{VPTPD}
	\begin{cases}
		 u_{n+1} = \arg \min \limits_{u}  \  f(u) + \langle u,B^{\top} p_{n} \rangle + \dfrac{1}{2\lambda}\left\Vert u-u_{n}\right\Vert_{I_u}^{2}, \\
		 p_{n+1} = \arg \min \limits_{p} \ g^{*}\left(p\right) - \langle Bu_{n+1},p \rangle + \dfrac{1}{2\sigma}\left\Vert p-p_{n}\right\Vert_{I_p}^{2} 
		  + \langle BI_u^{-1}(\nabla f(u_{n+1})+B^{\top}p_n),p \rangle, 
	\end{cases}
\end{align}
which can be viewed as a generalization of the TPD method \cite{chen2023transformed2,chen2023transformed1} to nonsmooth problems.
An important observation is that the triangular metric leads to an additional dual correction term (i.e., the last term on the right-hand side) inherited from the transformation that alleviates the primal--dual coupling and accelerates the convergence of the iteration, which is absent in block-diagonal preconditioned PDHG \eqref{PDHG}.
The transformation endows the TPD method with substantial potential for efficiently solving the nonlinear saddle point problem, yet poses new challenges when adapting to solving the JKO scheme of generalized Wasserstein gradient flows, such as construction of suitable preconditioners, computation of proximal operator, and strategy of adaptively selecting step size.





\section{Variable-Preconditioned Transformed primal--dual method for JKO scheme}\label{num method}

In this section, we will propose the Variable-Preconditioned Transformed primal--dual (VPTPD) method for solving the JKO scheme of generalized Wasserstein gradient flows, which is a nontrivial extension of the TPD method incorporating proximal splitting techniques to address the challenges arising from the nonsmoothness of the objective function. To clarify the different notations and their relationships involved in the method, we provide a summary of notations in Table \ref{tab:notation}.

\begin{table}[htbp]
\centering
\caption{Summary of notations.}
\label{tab:notation}
\small
\setlength{\tabcolsep}{4pt}      
\setlength{\extrarowheight}{2pt} 
\begin{tabular}{>{\raggedright\arraybackslash}p{0.48\textwidth} | >{\raggedright\arraybackslash}p{0.48\textwidth}}
\hline
\multicolumn{2}{c}{\textbf{Symbol: Description}} \\[4pt]
\hline
\(\mathcal{F}\): continuous objective of JKO scheme &
\(\mathcal{F}_h\): discrete objective of JKO scheme \\[4pt]
\(\hat{\mathcal{F}}_h\): regularized discrete objective for preconditioning &
\(\tilde{\nabla}\mathcal{F}_h\): an explicit subgradient of $\mathcal{F}_h$ \\[4pt]
\(\mathcal{J}\): continuous energy functional &
\(\mathcal{J}_h\): discrete energy functional \\[4pt]
\(\Phi\): Wasserstein distance functional &
\(\Phi_h\): discrete Wasserstein distance functional \\[4pt]
\(B\): discrete temporal-spatial divergence operator &
\(I_u, I_p\): primal/dual variable-preconditioner \\
\hline
\end{tabular}
\end{table}


\subsection{The VPTPD iterative scheme} 
Denoting the discretized objective function in \eqref{generalized discrete JKO} $\mathcal{F}_{h}:=\Phi_{h}/2\tau + \mathcal{J}_{h}$, the nonsmoothness of $\Phi_{h}$ makes the iterative scheme \eqref{TPD_discrete} of the TPD method with explicit gradient no longer applicable. By the idea of operator splitting, we replace the explicit gradient by the two-step proximal-splitting descent, i.e., $\partial \mathcal{F}_h = \nabla \mathcal{J}_{h}(u_n) + \frac{1}{2\tau} \partial\Phi_{h}(u_{n+1})$, where we apply the explicit gradient for the smooth part $\mathcal{J}_{h}$ and apply the proximal operator for the nonsmooth part $\Phi_{h}$, and obtain the following semi-implicit-explicit iterative scheme:
\begin{align}\label{semi-VPTPD}
	\begin{cases}
		u_{n+1} = \mathrm{Prox}_{\frac{\lambda}{2\tau}\Phi_{h}}^{I_{u}}\left(\bar{u}_{n} - \lambda (I_{u})^{-1}\left(\nabla\mathcal{J}_{h}\left(\bar{u}_{n}\right)+B^{\top}\bar{p}_{n}\right)\right),\\
		p_{n+1} = \bar{p}_{n} + \sigma (I_p)^{-1}\left(Bu_{n+1} - b - B(I_{u})^{-1}\left(\tilde{\nabla}\mathcal{F}_{h}\left(u_{n+1}\right)+B^{\top}\bar{p}_{n}\right)\right),\\
		\bar{u}_{n+1} = \left(1 + \kappa_{1}\right) u_{n+1} - \kappa_{1} u_{n}, \quad
		\bar{p}_{n+1} = \left(1 + \kappa_{2}\right) p_{n+1} - \kappa_{2} p_{n},
	\end{cases}
\end{align}
where $(\bar{u}_n,\bar{p}_n)$ are the estimates of $(u_{n+1},p_{n+1})$ by an extrapolation with (typically positive) coefficients $\kappa_{1}$ and $\kappa_{2}$ for faster convergence, and the generalized proximal operator for a convex functional $f$ with a SPD preconditioner $M$ is defined as
\begin{equation}\label{illustration-VPTPD}
	\mathrm{Prox}_{f}^{M}\left(v\right) := \arg \min \limits_{u}\left\{\dfrac{1}{2}\left\Vert u-v\right\Vert_{M}^{2} + f\left(u\right)\right\}, \quad \text{where } \left\Vert u\right\Vert_{M}^{2}=u^{\top}Mu.
\end{equation}

In particular, we define $\tilde{\nabla}\mathcal{F}_{h}\left(u_{n+1}\right)$ in the update of the dual variable \eqref{semi-VPTPD} as an explicit computable subgradient of $\mathcal{F}_h$: 
\begin{equation}\label{nabla_F}
\tilde{\nabla} \mathcal{F}_h(u_{n+1}) \in \partial \mathcal{F}_h (u_{n+1}) = \nabla\mathcal{J}_h(u_{n + 1}) + \frac{1}{2\tau} \partial \Phi_h (u_{n+1}).
\end{equation}
Based on the exact inclusion relation derived directly from the first line in \eqref{semi-VPTPD}:
\begin{equation}
	u_{n+1} \in \bar{u}_{n} - \lambda I_{u}^{-1} \left(\nabla\mathcal{J}_{h}\left(\bar{u}_{n}\right)+ B^{\top}\bar{p}_{n}\right) - \dfrac{\lambda}{2\tau} I_{u}^{-1} \partial\Phi_{h}\left(u_{n+1}\right),
\end{equation}
which gives a computable element of $\partial\Phi_h(u_{n+1})$, yielding the explicit expression of the subgradient $\tilde{\nabla} \mathcal{F}_h$:
\begin{equation}
\tilde{\nabla} \mathcal{F}_h(u_{n+1}) = \frac{1}{\lambda} I_{u}\left(\bar{u}_{n}-u_{n+1}\right)-B^{T}\bar{p}_{n}-\nabla\mathcal{J}_{h}\left(\bar{u}_{n}\right)+\nabla\mathcal{J}_{h}\left(u_{n+1}\right) \in \partial\mathcal{F}_h(u_{n+1}).
\end{equation}

As described in \autoref{TPD-method}, the key of the TPD method is the selection of appropriate linear operators $I_u$ and $I_p$ based which the decoupling transformation and preconditioning are performed. As suggested by \cite{chen2023transformed2}, a typical choice of $I_u^{-1}$ is the approximation of the inverse of the Hessian $\nabla^2 \mathcal{F}_h(u)$ when $\mathcal{F}_h$ is twice-differentiable, which unfortunately does not align with our case. Therefore, we construct a regularized objective function $\widehat{\mathcal{F}}_{h}$:
\begin{equation}\label{modify}
	\widehat{\mathcal{F}}_{h}\left(u\right) = \dfrac{1}{2\tau} \sum_{\boldsymbol{i}}\widehat{\varphi}\left(\rho_{\boldsymbol{i}},\boldsymbol{m}_{\boldsymbol{i}}\right)\Delta V + \mathcal{J}_{h}\left(u\right), \quad \text{where } \widehat{\varphi}\left(\rho,m\right) =
	\displaystyle\frac{{\left\Vert \boldsymbol{m}\right\Vert}^{2}}{M\left(\rho\right)+r},
\end{equation}
where $r$ is a small regularization parameter. Note that $\widehat{\mathcal{F}}_{h}\left(u\right)$ is twice differentiable and strictly convex with suitable $\tau$ when $M(\rho)\geq 0$, and we could take $\nabla^2\widehat{\mathcal{F}}_{h}\left(u_{n}\right)$ as the preconditioner for achieving optimal preconditioning for nonlinear $\nabla \mathcal{F}_h$. However, this iteration-evolving and non-diagonal preconditioner poses significant difficulties in computing the proximal operator and matrix inversion (see Remark~\ref{inverse method}). To achieve the delicate balance between the convergence rate of iteration and the computational complexity of each iteration, we adopt the matrix comprising solely its diagonal entries as the preconditioner that is PDE-time-dependent but iteration-invariant. Specifically, for the $(k+1)$-th JKO step \eqref{generalized discrete JKO}, we take $I_u=\mathrm{diag}(\nabla^2\widehat{\mathcal{F}}_{h}(u^{k}))$ and $I_{p} = B(I_{u})^{-1}B^{T}$, where $u^k$ is the minimizer of the previous JKO step and satisfies $M(\rho^k)\geq 0$. 

We emphasize that the regularization in \eqref{modify} is merely for constructing a stable preconditioner $I_u$, and the transformation and preconditioning induced thereby does not change the equivalence between the original saddle point system \eqref{PD_saddle} and the transformed system \eqref{TPD_saddle}. We actually still solve the primal proximal update for the original objective (without regularization) in \eqref{semi-VPTPD}. In practice, the regularization parameter $r$ with a wide range of values (typically around $10^{-3}\sim 10^{-8}$) robustly provides satisfactory preconditioning performance (see \autoref{2D-CH}). 
With the tailored preconditioners, our VPTPD demonstrates remarkably improved convergence and high computational efficiency as shown in \autoref{compare example}. 


\begin{myremark}[Selection of $I_u$ and $I_p$]\label{inverse method}
	Although the full Hessian would be an ideal preconditioner for nonlinear $\nabla \mathcal{F}_h$, it is impractical: it destroys the separability of the primal proximal step (leading to solving a large coupled nonlinear system) and yields a far more coupled (potentially dense) Schur-complement system in the dual update, thereby greatly increasing per-iteration cost and memory usage. An observation is that $\mathcal{J}_h$ mainly contributes diagonal or Laplacian-like components and stencil couplings come from $\hat{\varphi}$ which are typically moderate compared to the diagonal contribution. Hence, the diagonal preconditioner offers a favorable compromise to strike a better balance between iteration count and per-iteration complexity.
	
	With diagonal $I_u$, its inversion (twice per each iteration) is easily obtained, and the resulting primal proximal operator in \eqref{semi-VPTPD} decouples into component-wise subproblems that can be solved efficiently by a convergent, bound-preserving Newton solver and are amenable to parallelization (see \autoref{subsec:prox}).

	
	For the dual update, since $I_{p} = B(I_{u})^{-1}B^{T}$ is a sparse matrix with only $O(N)$ nonzeros, its inversion can be efficiently calculated via a sparse Cholesky factorization with approximate minimum degree ordering for moderate sizes \cite{chen2008algorithm}, or via incomplete Cholesky-preconditioned conjugate gradient method for large-size systems (typically within 5 iterations for 3D simulations). Importantly, since $I_p$ is fixed within each JKO step and only updating across successive PDE time steps, only one matrix decomposition is required per JKO step, while PDE-time-evolving preconditioners are sufficient for effective preconditioning given that consecutive solutions typically vary mildly. 

\end{myremark}

\begin{myremark}[PDHG-type methods for JKO scheme]\label{compare_VPTPD}
In previous works \cite{carrillo2022primal,carrillo2024structure} for computing the dynamic JKO scheme \eqref{generalized discrete JKO}, the PDJKO and PrePDJKO methods utilizes PD3O splitting method \cite{Yan2018NewPDThreeFunctions} with a relaxed penalty term $g(Bu)=i_{\delta}\left(Bu\right)$ ($\delta>0$) for the constraint of continuity equation, which introduces $g^*(p)= \langle b, p\rangle + \delta \lVert p\rVert$ in the saddle point formulation \eqref{minimax problem2}. This extra term in fact adds regularization on the diagonal to make the operator $\mathcal{K}$ in \eqref{PD_saddle} diagonal-dominant such that resulting system exhibits better properties for PD-type methods to be convergent. PDJKO \cite{carrillo2022primal} employs identical preconditioners $P_{\mathrm{PDJKO}}=\mathrm{diag}\{\frac{1}{\lambda}I,\frac{1}{\sigma}I\}$ in \eqref{PDHG} and suffers from slow convergence due to the joint restriction on step sizes $\lambda\sigma\|BB^{\top}\|<1$. Inspired by \cite{liu2021acceleration}, PrePDJKO \cite{carrillo2024structure} performs preconditioning in the dual update with a constant matrix-valued preconditioner $P_{\mathrm{PrePDJKO}}=\mathrm{diag}\{\frac{1}{\lambda}I,\lambda BB^{\top}\}$ and achieve better performance for nonlinear mobilities. Interestingly, we can explicitly show that PrePDJKO coincides with G-Prox PDHG \cite{Jacobs2019Gprox} applied to a reformulation of the original one-step JKO problem. In practice, however, the PrePDJKO implemented in this work adopts a proximal-splitting treatment for the underlying three-operator structure, and can therefore be viewed as a proximal-splitting variant of G-prox PDHG (see details in \ref{Gprox-PrePD}).

However, constant-valued preconditioners still limits the computational efficiency of nonlinear PDE systems, especially for large-size systems. In contrast, our VPTPD method employs a decoupling-transformation together with a PDE-time-dependent preconditioner, which leads to a significant acceleration of performance. In \autoref{compare example}, we compare the performance of VPTPD with PDHG-type methods, and demonstrate the superiority of VPTPD in terms of convergence performance and computational efficiency for JKO scheme.
\end{myremark}

\begin{myremark}[Convergence of VPTPD method]\label{VPTPD_convergence}
In \cite{chen2023transformed2,chen2023transformed1}, the authors established the exponential convergence of the continuous dynamical system \eqref{TPD2_flow} and the corresponding linear convergence of its explicit Euler discretization \eqref{TPD_discrete}, under the condition of $\mathcal{F}$ being smooth and strongly convex.
In particular, they have explicitly established the convergence rate which relies on the factor $1- \mathcal O(\kappa_f^{-3}(I_u))$, provided that $I_p$ presents a reasonably good approximation to $BI^{-1}_uB^{\top}$, where $\kappa_f(I_u) = L_f(I_u)/\mu_f(I_u)$ with $L_f(I_u)$ and $\mu_f(I_u)$ being Lipschitz and convexity constants of the functional $\mathcal{F}$ measured under the metric $I_u$.
If those preconditioners are constructed good enough, such that the resulting condition numbers are independent of mesh size, 
then the resulting TPD algorithm may admit a mesh-independent convergence rate. 

Moreover, the convergence of the VPTPD method \eqref{semi-VPTPD} becomes even more challenging due to the nonsmoothness of the objective function $\mathcal{F}_h$. While a complete convergence theory for this setting is not yet available, existing analyses for the smooth TPD framework \cite{chen2023transformed2,2025GuoZou} suggest that admissible step sizes should be governed by the Lipschitz continuity of the smooth part $\nabla \mathcal{J}_h$ measured in the $I_u$-metric and the condition number of the Schur-complement $BI_u^{-1}B^{\top}$ measured in the $I_p$-metric. 
In practice, we would expect that the synergy between decoupling-transformation and variable preconditioning can substantially relax restrictions on step sizes and enhances convergence acceleration, as supported by numerical results in \autoref{sec:FP1D}. We also observe robust convergence of VPTPD for a broad range of step sizes $(\lambda,\sigma)$ (typically of order $\mathcal{O}(1)$) in various numerical experiments, which empowers an effective adaptive step-size strategy to achieve faster convergence.

Recently, G-prox PDHG \cite{Jacobs2019Gprox} has been shown to achieve mesh-independent convergence for certain nonsmooth PDE-constrained optimization problems. While this is a highly attractive feature, its application to more compelx problems, such as the one-step JKO minimization involving three operators, may lead to expensive proximal subproblems; see the discussion in \ref{Gprox_prox}. 
One may consider more advanced iterative schemes, such as the semi-implicit-explicit scheme in \eqref{semi-VPTPD},
but the convergence analysis together with the mesh-independent property becomes unknown yet.
Indeed, we actually observe the mesh-dependent convergence of PrePDJKO and VPTPD for certain problems shown in Figure \ref{merge}. In our comprehensive numerical tests in \autoref{compare example}, VPTPD admits efficiently computable proximal operators and accelerated convergence, but may still exhibit mesh dependence. 
In short, a complete understanding of VPTPD convergence behavior is highly nontrivial and certainly deserves further investigation in future study.


\end{myremark}


\subsection{Computation of proximal operator}\label{subsec:prox}
The ease of computing the proximal operator in the algorithm \eqref{semi-VPTPD} is crucial for its effectiveness and efficiency:
\begin{equation}\label{d-Proximal}
	\mathrm{Prox}_{\frac{\lambda}{2\tau}\Phi_{h}}^{I_{u}}\left(\hat{u}\right) = \arg\min\limits_{u}\ 
	\dfrac{1}{2}\left\Vert u - \hat{u}\right\Vert_{I_u}^{2} + \frac{\lambda}{2\tau}\Phi_{h}\left(u\right).
\end{equation}
Thanks to the separable structure of $\Phi_h(u) =\sum_{\boldsymbol{i}}\varphi(\rho_{\boldsymbol{i}},\boldsymbol{m}_{\boldsymbol{i}})\Delta V$, we are able to compute the proximal operator with scalar preconditioning $I_u=\frac{1}{\lambda}I$ in our previous work \cite{carrillo2022primal,carrillo2024structure}. However, non-diagonal preconditioners will lead to a large coupled nonlinear system, which can not be efficiently computed. Therefore, we choose a diagonal preconditioner $I_u$ that preserves the separable structure of the proximal operator:
\begin{equation}\label{prox_sep}
	\mathrm{Prox}_{\frac{\lambda}{2\tau}\Phi_{h}}^{I_{u}}(\hat{u}) = \left(\mathrm{Prox}_{\frac{\zeta}{2}\varphi}^{\mathcal{D}_{\boldsymbol{i}}}(\hat{\rho}_{\boldsymbol{i}},\hat{\boldsymbol{m}}_{\boldsymbol{i}})\right)_{\boldsymbol{i}},
\end{equation}
where $\zeta = \lambda \Delta V/\tau$ and $\mathcal{D}_{\boldsymbol{i}}=\mathrm{diag}(\mathcal{D}_\rho,\{\mathcal{D}_{m_j}\}_{j=1}^d)$ consists of the elments of $I_u$ corresponding to $(\rho_{\boldsymbol{i}},\boldsymbol{m}_{\boldsymbol{i}})$ for certain $\boldsymbol{i}$. This separable structure \eqref{prox_sep} not only reduces the original problem \eqref{d-Proximal} to many independent scalar subproblems, but also enables the applicability of parallelization. By dropping the index $\boldsymbol{i}$, we now focus on the following component-wise subproblem:
\begin{equation}\label{subproblem}
	\mathrm{Prox}_{\frac{\zeta}{2}\varphi}^{\mathcal{D}}(\hat{\rho},\hat{\boldsymbol{m}})=\arg\min\limits_{\rho,\boldsymbol{m}} \dfrac{1}{2}\mathcal{D}_{\rho}(\rho-\hat{\rho})^{2} + \dfrac{1}{2}\sum_{j=1}^{d}\mathcal{D}_{m_{j}}(m_{j}-\hat{m}_{j})^{2}+ \dfrac{\zeta}{2}\varphi\left(\rho, \boldsymbol{m}\right).
\end{equation}

One of our main contribution is concluded in the following theorem, that proves the existence and uniqueness of bounded solution to \eqref{subproblem}, and provides a convergent and bound-preserving Newton iteration method to solve it. This result is an extension of our previous work \cite{carrillo2024structure} to account for two general types of concave mobilities for generalized Wasserstein gradient flows \cite{lisini2012cahn} and diagonal preconditioners.

\begin{mytheorem}
	For concave mobilities in the following form:
	\begin{align*}
	\mathrm{(I)}\quad  &  M(\rho) = (\rho - \alpha)^{\xi}(\beta - \rho)^{\xi},  && \text{with $\alpha<\beta$, $\xi\in(0,1]$}, && \text{for $\alpha\leq\rho\leq\beta$}, \\
	\mathrm{(II)}\quad &  M(\rho) = \rho^{\xi},  && \text{with $\xi\in(0,1]$}, && \text{for $\rho\geq 0$},
	\end{align*}
	the problem \eqref{subproblem} admits a unique solution $(\rho^{*},\boldsymbol{m}^{*})$ within the bounds such that $M(\rho^*)\geq 0$. The distribution of the solution $\rho^*$ depending on the input $(\hat{\rho},\hat{\boldsymbol{m}})$ is provided in \autoref{Newton initial table}. The Newton's method with the tailored initial values in \autoref{Newton initial table} is guaranteed to converge to the desired bounded solution $(\rho^{*},\boldsymbol{m}^{*})$.
\end{mytheorem}
\textbf{Proof:}
	Given that $\mathcal{D}_{\rho}, \mathcal{D}_{m_j}, \eta$ are positive and $\varphi(\rho,\boldsymbol{m})$ in \eqref{W_distance} is convex, the ojective function of \eqref{subproblem} is strongly convex and hence yields a unique solution. Moreover, the definition of $\varphi(\rho, \boldsymbol{m})$ guarantees the bounds of the solution $\rho^*$ such that $M(\rho^{*}) \geq 0$. 
	
	To discuss the convegence of Newton's method, we first consider the minimization \eqref{subproblem} restricted to the interval where $M(\rho)> 0$. First-order optimal condition yields:
	\begin{align}\label{derivative of N(u)}
		\begin{cases}
			\mathcal{D}_{\rho}\left(\rho - \hat{\rho}\right) -\dfrac{\zeta M'\left(\rho\right)}{2}\left(\dfrac{\sum_{j=1}^{d} m_{j}^{2}}{M^{2}\left(\rho\right)}\right) = 0,\\
			\mathcal{D}_{m_{j}}\left(m_{j} - \hat{m}_{j}\right)+\zeta\dfrac{m_{j}}{M\left(\rho\right)}= 0,
			\qquad j = 1,2,\cdots,d.
		\end{cases}
	\end{align}
	By the second equation, we obtain $m_{j}=\frac{M(\rho)\mathcal{D}_{m_j}\hat{m}_{j}}{\mathcal{D}_{m_j}M(\rho)+\zeta}$. By inserting the relation into the first equation, we obtain the equation for $\rho$ only:
	\begin{equation}\label{L_rho}
		L\left(\rho\right) := \mathcal{D}_{\rho}\left(\rho-\hat{\rho}\right) - \dfrac{\zeta M'\left(\rho\right)}{2}\sum_{j=1}^{d} \left(\dfrac{\mathcal{D}_{m_j}\hat{m}_{j}}{\mathcal{D}_{m_{j}}M\left(\rho\right)+\zeta}\right)^{2} = 0.
	\end{equation}
	It is easy to check that $L'>0$ provided that $M(\rho)>0$ and $M''(\rho)<0$, and hence $L(\rho)=0$ has at most a root within the interval where $M(\rho)>0$. In particular, if there exists a root within the bounds, it is the desired solution of the proximal operator $\rho^*$; otherwise, the solution of proximal operator must be obtained at the endpoints where $M(\rho)=0$. With $\rho^*$, the momentum $\boldsymbol{m}^{*}$ can be obtained accordingly. 
	
	By evaluating $L(\rho)$ at specific points, the interval containing the root of equation \eqref{L_rho} can be identified. Together with the information of $L''(\rho)$ over the interval, we obtain a tailored plan of initial values for the Newton's method to guarantee the convergence to the solution $\rho^*$. Depending on the input $(\hat{\rho},\hat{\boldsymbol{m}})$, the interval of the root, the sign of $L''(\rho)$, and the correponding initial value $\rho_0$ for Newton's method or the optimal solution $\rho^*$ to \eqref{subproblem} are summarized in \autoref{Newton initial table}, where the parameters $\rho_m$, $C_1$ and $C_2$ are given by:
	\begin{equation}
		\rho_m = \dfrac{\alpha + \beta}{2}, \quad
		C_{1} =  \dfrac{-1}{2\zeta\mathcal{D}_{\rho}}\sum_{j=1}^{d}\left(\mathcal{D}_{m_{j}}\hat{m}_{j}\right)^{2}, \quad
		C_{2} = \dfrac{\beta-\alpha}{2\zeta\mathcal{D}_{\rho}}\sum_{j=1}^{d}\left(\mathcal{D}_{m_{j}}\hat{m}_{j}\right)^{2}.
	\end{equation}
	\qed
	
	\begin{table}[htbp]
		\centering
		\footnotesize
		\caption{Initial value plan for convergent Newton iteration method}
		\label{Newton initial table}
		\begin{tblr}{
				colspec = {X[c] X[c] X[c] X[c]},
				hlines = {dashed},
				hline{3, 5-6, 8-9, 12-13} = {solid},
				hline{1, 2, 18} = {1pt, solid},
				cell{1}{1-4} = {font=\bfseries},
			}
			$\mathrm{Range}(\hat{\rho})$ & $\mathrm{Root\,interval\,of\,L}$ & $\mathrm{sign}(L'')$& $\rho_0,\rho^*$ \\
			\SetCell[c=4]{c} $M(\rho) = \rho^{\xi}, \xi\in(0, 1), \rho\geq 0$  & & &\\
			$\hat{\rho} \leq \delta $ & $(0,+\infty)$ & $L''<0$ & $\rho_0 = \delta$ \\
			$\hat{\rho} > \delta $ & $(\hat{\rho},+\infty)$ & $L''<0$ & $\rho_0 = \hat{\rho}$\\
			\SetCell[c=4]{c} $M(\rho) = \rho, \rho \geq 0$  & & &\\
			$\hat{\rho} \leq C_1 $ & $(-\infty, 0]$ & * & $\rho^{*} = 0$\\
			$\hat{\rho} > C_1 $ & $(0,+\infty)$ & $L''<0$ & $\rho_{0} = 0$\\
			\SetCell[c=4]{c} $M(\rho) = (\rho-\alpha)^{\xi}(\beta-\rho)^{\xi}, \xi\in(0, 1), \alpha \leq \rho \leq \beta$  & & &\\
			$\hat{\rho} < \rho_{m}$ & $(\alpha, \rho_{m})$ & $L''<0$ & $\rho_0 = \alpha + \delta$ \\
			$\hat{\rho} = \rho_{m}$ & $\left\{\rho_{m}\right\}$ & * & $\rho^* = \rho_{m}$ \\
			$\hat{\rho} > \rho_{m}$ & $(\rho_{m},\beta)$ & $L''>0$ & $\rho_0 = \beta - \delta$ \\
			\SetCell[c=4]{c} $M(\rho) = (\rho-\alpha)(\beta-\rho), \alpha \leq \rho \leq \beta$  & & &\\
			$\hat{\rho} \leq \alpha - C_2$ & $(-\infty,\alpha]$ & * & $\rho^* = \alpha$\\
			$\alpha - C_2 < \hat{\rho} < \rho_{m}$ & $(\alpha, \rho_{m})$ & $L''<0$ & $\rho_{0} = \alpha$\\
			$\hat{\rho} = \rho_{m}$ & $\left\{ \rho_{m} \right\}$ & * & $\rho^{*} = \rho_{m}$\\
			$\rho_{m} < \hat{\rho} < \beta + C_2$ & $(\rho_{m}, \beta)$ & $L''>0$ & $\rho_{0} = \beta$\\
			$\hat{\rho} \geq \beta + C_2$ & $[\beta, +\infty)$ & * & $\rho^{*} = \beta$\\
		\end{tblr}
	\end{table}

\subsection{Adaptive step-size strategy}\label{ada strategies}
For the step size of the TPD method \eqref{TPD_discrete}, there exists a theoretical upper bound (when $\lambda=\sigma$) to guarantee its convergence \cite{chen2023transformed1,chen2023transformed2}.  
For the semi-implicit-explicit iterative scheme \eqref{semi-VPTPD} of our VPTPD method, the convergence condition on the step sizes should depend on the properties of objective functional $\mathcal{F}_h$ under the metric induces by preconditioners (see Remark~\ref{VPTPD_convergence}). In practice, the decoupling transformation and variable preconditioning improve the conditioning of the underlying saddle system, thereby relaxing the step-size constraints and allowing for a more flexible adaptive update of $\lambda$ and $\sigma$ to accelerate convergence.
Specifically, we propose an adaptive-step strategy for $\lambda$ based on the search directions $\Delta u_{n} = u_{n} - u_{n-1}$ from previous iterates and the convergence monitors defined as follows:
\begin{align}\label{criteria}
	\begin{cases}
		e^{1}_{n+1}:=\frac{\left\Vert Bu_{n+1} - b \right\Vert_{2}}{\sqrt{N}}, \quad
		e^{2}_{n+1}:=\mathrm{rel}(\mathcal{J}_{h}(u_{n+1})),\quad
		e^{3}_{n+1}:=\mathrm{rel}(\Phi_{h}(u_{n+1})),\\
		e^{4}_{n+1}:=\mathrm{rel}(u_{n+1}),\quad e^{5}_{n+1}:=\mathrm{rel}(p_{n+1}),
	\end{cases}
\end{align}
where $\mathrm{rel}(v_{n+1}) = \frac{\lvert v_{n+1} - v_{n} \rvert}{\lvert v_{n+1} \rvert}$ represents the relative error of variable $v$. 
The idea of the adaptive strategy is to increase the step size $\lambda$ when the search directions of previous two iterations $\Delta u_{n}$ and $\Delta u_{n - 1}$ are sufficiently aligned, and to decrease $\lambda$ when the convergence monitors $(e^{1}, e^{2}, e^{3}, e^{4}, e^{5})$ exhibit an increasing trend. Based on these mechanisms, the adaptive step-size strategy is given by:
\begin{align}\label{step-size-ada}
	\begin{split}
		\lambda_{n}^{*} &= \begin{cases}
			\min(\lambda_{n-1}\gamma_{+}, \lambda_{\mathrm{max}}), & \mathrm{if}\ \langle \Delta u_{n}, \Delta u_{n-1}\rangle > \theta_{\mathrm{max}} \lvert \Delta u_{n} \rvert  \lvert \Delta u_{n-1} \rvert\\
			\lambda_{n-1}, &\mathrm{otherwise},
		\end{cases}\\
		\lambda_{n} &= \begin{cases}
			\max(\lambda_{n}^{*}\gamma_{-}, \lambda_{\mathrm{min}}), & \mathrm{if} \, \exists \ i, \,e^{i}_{n} > e^{i}_{n-1},\\
			\lambda_{n}^{*}, & \mathrm{otherwise},
		\end{cases}
	\end{split}
\end{align}
where $\theta_{\mathrm{max}}$ is a threshold parameter, and $\gamma_{+} > 1$ and $\gamma_{-} \in (0, 1)$ are the step-size increment and decrement factors, respectively, satisfying $\gamma_{+} \cdot \gamma_{-} < 1$.


\subsection{Complete algorithm of VPTPD method} In the end, we summarize the VPTPD method (with options of adaptive strategy) as follows:
\begin{algorithm}
	\caption{VPTPD method for $(k+1)$-th step of dynamic JKO}
	\label{VPTPD-alg}
	\SetAlgoLined
	\KwIn{$(u_0, p_0)=(u^{k},p^{k})$ \quad $\triangleright$ initial values by solution of $k$-th JKO step}
	\KwIn{$\tau$, $\lambda_0$, $\sigma$, $\kappa_{1}$, $\kappa_{2}$, $\mathrm{tol}$, $\mathrm{TOL}$, $\mathrm{Iter}_{\max}$}
	\BlankLine
	Set: $\bar{u}_{0} = u_{0}$, $\bar{p}_{0} = p_{0}$, $\lambda = \lambda_{0}$, $\sigma = \sigma_{0}$, $n = 0$\;
	\While{$n < \mathrm{Iter}_{\max}$}{
		\Repeat{stopping criteria are achieved: $e^{1}_{n+1} < \mathrm{tol}$ and $\max_{2\leq i\leq 5}(e^{i}_{n+1}) < \mathrm{TOL}$}{
			\eIf{Adaptive Step Size}{
				Set $\lambda = \lambda_{n}$ by \eqref{step-size-ada}\;
			}{
				Set $\lambda = \lambda_{0}$\;
			}
			$u_{n+1} = \mathrm{Prox}_{\frac{\lambda_{n}}{2\tau}\Phi_{h}}^{I_{u}}\left(\bar{u}_{n} - \lambda_{n} (I_{u})^{-1}\left(\nabla\mathcal{J}_{h}\left(\bar{u}_{n}\right)+B^{\top}\bar{p}_{n}\right)\right)$\;
			$p_{n+1} = \bar{p}_{n} + \sigma (I_p)^{-1}\left(Bu_{n+1} - b - B(I_{u})^{-1}\left(\tilde{\nabla}\mathcal{F}_{h}\left(u_{n+1}\right)+B^{\top}\bar{p}_{n}\right)\right)$\;
			$\bar{u}_{n+1} = \left(1 + \kappa_{1}\right) u_{n+1} - \kappa_{1} u_{n}$\;
			$\bar{p}_{n+1} = \left(1 + \kappa_{2}\right) p_{n+1} - \kappa_{2} p_{n}$\;
		}
	}
	\KwOut{$(u^{k+1}, p^{k+1})=(u_{n+1},p_{n+1})$ \quad $\triangleright$ solution of $(k+1)$-th JKO step}
\end{algorithm}

\section{Numerical experiments}\label{experiment}
In this section, we perform a variety of numerical experiments to showcase the high computational efficiency of our VPTPD method compared to existing methods, and the flexibility and applicability of our method to several challenging problems.

Unless otherwise specified, in the following numerical experiments, we set the error tolerance for stopping criteria tol$=10^{-5}$ and TOL$=10^{-5}$, the regularization parameter for constructing preconditioners $r=10^{-5}$, the extrapolation coefficient $\kappa_{1}=0.3$ and $\kappa_{2}=0.8$, the initial step size $\lambda_{0} = 1$ and $\sigma_{0} = 0.3$, the adaptive step size parameters $\lambda_{\max} = 3\lambda_{0}, \lambda_{\min} = \lambda_{0}$, $\gamma_{+} = 1.2, \gamma_{-} = 0.8$, and $\theta_{\max} = 0.99$.

\subsection{Comparison experiments}\label{compare example}
We conduct a comparative study on the computational efficiency of the proposed VPTPD method \autoref{VPTPD-alg} and the aforementioned PDHG-type methods (such as PDJKO and PrePDJKO) (see Remark~\ref{compare_VPTPD}), in terms of the number of iterations and CPU time, via several benchmark simulations spanning from 1D to 3D problems.

\subsubsection{1D generalized Fokker-Planck equation}\label{sec:FP1D}
Consider the 1D Fokker-Planck equation of the following general form:
\begin{equation}\label{1Dsaturation}
	\rho_{t} = \nabla\cdot\left(M(\rho)\nabla\left(\mathcal{H}^{'}(\rho)+V(x)\right)\right),
\end{equation}
where the corresponding energy is $\mathcal{J}\left(\rho\right) = \int_{\Omega}\left(\mathcal{H}(\rho) + V(x)\rho\right)\mathrm{d}x$ with a confining drift potential $V(x)$.

\begin{figure}[htbp]
    \centering
    \begin{minipage}[c][0.2\textheight][t]{0.38\textwidth}
        \includegraphics[width=0.9\textwidth]{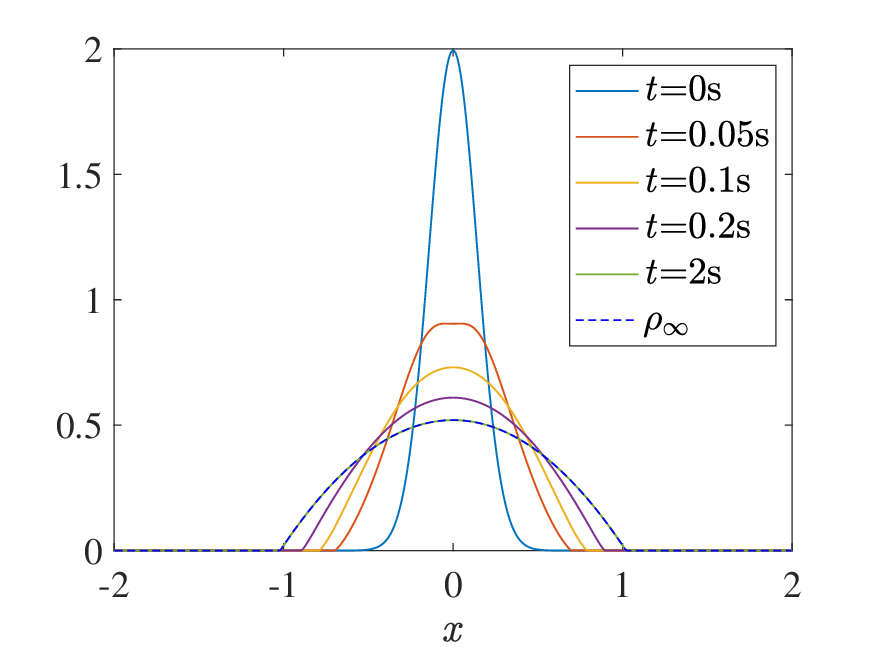}
        \figcaption{Evolution of solutions for the 1D Fokker-Planck equation.}
        \label{1Dfokkerplanck}
    \end{minipage}
    \begin{minipage}[c][0.2\textheight][b]{0.6\textwidth}
        \tabcaption{Number of iterations and CPU time for solving 1D Fokker-Planck equation \eqref{1Dsaturation} for $t \in [0, 2]$. $N_{x} = 400,\;\tau = 0.05, \; \lambda_0 = 2.5$, \; $r=10^{-3}$.}
        \label{tab:1D fokker-planck}
        \centering
        \renewcommand{\arraystretch}{1.8}
        \resizebox{1.0\columnwidth}{!}{
        \begin{tabular}{ccccc} 
            \hline
            \multicolumn{5}{c}{\textbf{1D Fokker-Planck equation}}\\
            \hline
            & PDJKO & PrePDJKO & VPPDJKO & VPPDJKO(s)\\ 
            \hline
            Number of iterations  & 68040  & 7398(-89.1\%) & 6238(-90.8\%)   & * \\
            \hline
            & TPD & PreTPD & VPTPD & VPTPD(s)\\ 
            \hline
            Number of iterations  & *  & 25924(-61.9\%) & 4713(-93.1\%) & 3559(-94.8\%) \\ 
            \hline
        \end{tabular}}
    \end{minipage}
\end{figure}

In the first example, we take $M(\rho)=\rho$, $\mathcal{H}(\rho)=\rho^2$, $V(x)=\left\vert x\right\vert^{2}$, and solve the equation for $t \in [0,2]$ with the initial condition $\rho_0(x)=\frac{5}{\sqrt{2\pi}}e^{-25x^2}$ in the computational domain $\Omega = [-2,2]$. The evolution of solution is shown in \autoref{1Dfokkerplanck} and the numerical result is well-matched with the analytical equilibrium solution \cite{carrillo2022primal}. 

We mainly compare the convergence behavior of the PDJKO and TPD methods with two types of preconditioning: constant-valued preconditioning $I_u=I$ and $I_p=BB^{\top}$, and variable preconditioning $I_u=\mathrm{diag}(\nabla^2\widehat{\mathcal{F}}_{h}(u^{k}))$ and $I_{p} = B(I_{u})^{-1}B^{T}$. For clarity, we adopt a prefix-suffix notation to distinguish these variants: the prefixes ``Pre" and ``VP" indicate constant and variable preconditioning, respectively, whereas the suffixes ``PDJKO'' and ``TPD'' specify the underlying algorithmic framework; an appended ``(s)" indicates the use of an adaptive step-size strategy.
We mention that, even though both PrePD and PreTPD employ the same scaling preconditioner, PreTPD only applies an implicit discretization to the primal variable but not to the dual variable, which may slow down its convergence.

The convergence results are summarized in Table~\ref{tab:1D fokker-planck}. It can be seen that the transformation alone (TPD) is unstable and fails to converge for the highly nonlinear problem. However, when combined with preconditioning, TPD becomes stable and overall superior than PDJKO. In particular, compared to constant-valued preconditioning (PreTPD), variable preconditioning significantly enhances the acceleration (VPTPD). The adaptive step-size strategy further improves the efficiency (VPTPD(s)), but fails to work with PDJKO methods possibly due to the severe restrictions on step sizes for PDHG-type methods. 

In the second example, we consider the porous medium equation by taking $M(\rho)=\rho$, $\mathcal{H}(\rho)=\rho^2$, $V(x)=0$, and solve the equation for $t \in [0,2]$ with the initial condition $\rho_0(x)=\frac{1}{\sqrt{2\pi\theta^2}} e^{\frac{x^2}{\theta^2}}, \theta = 0.2$ in the computational domain $\Omega = [-4,4]$. With this choice of energy functional, we are able to efficiently compute the proximal operator of $\mathcal{F}=\Phi/2\tau+\mathcal{J}$ for implementing the G-prox PDHG \cite{Jacobs2019Gprox}; see details in \ref{Gprox_prox}.

\begin{figure}[htbp]
\centering
\includegraphics[width=0.48\textwidth]{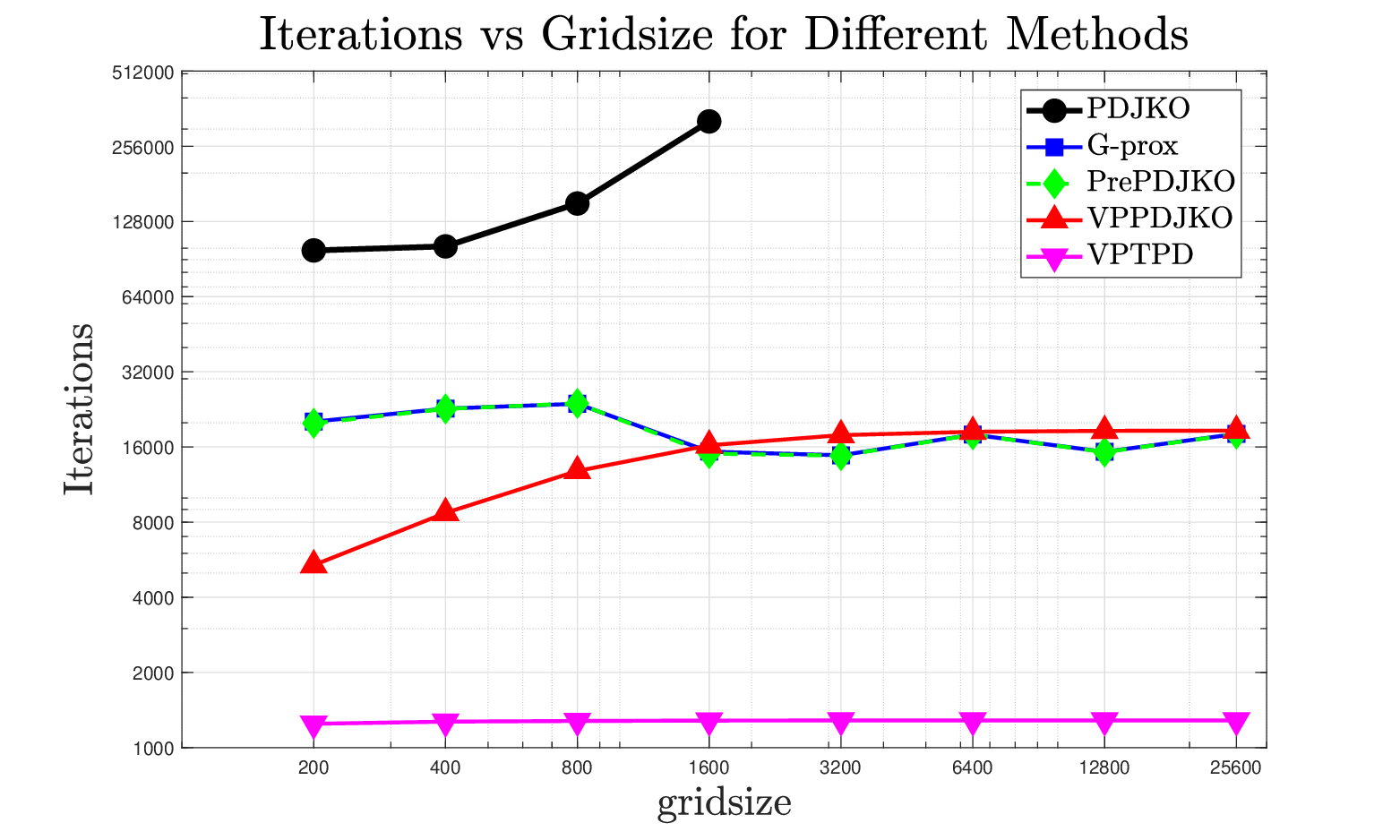}
\includegraphics[width=0.48\textwidth]{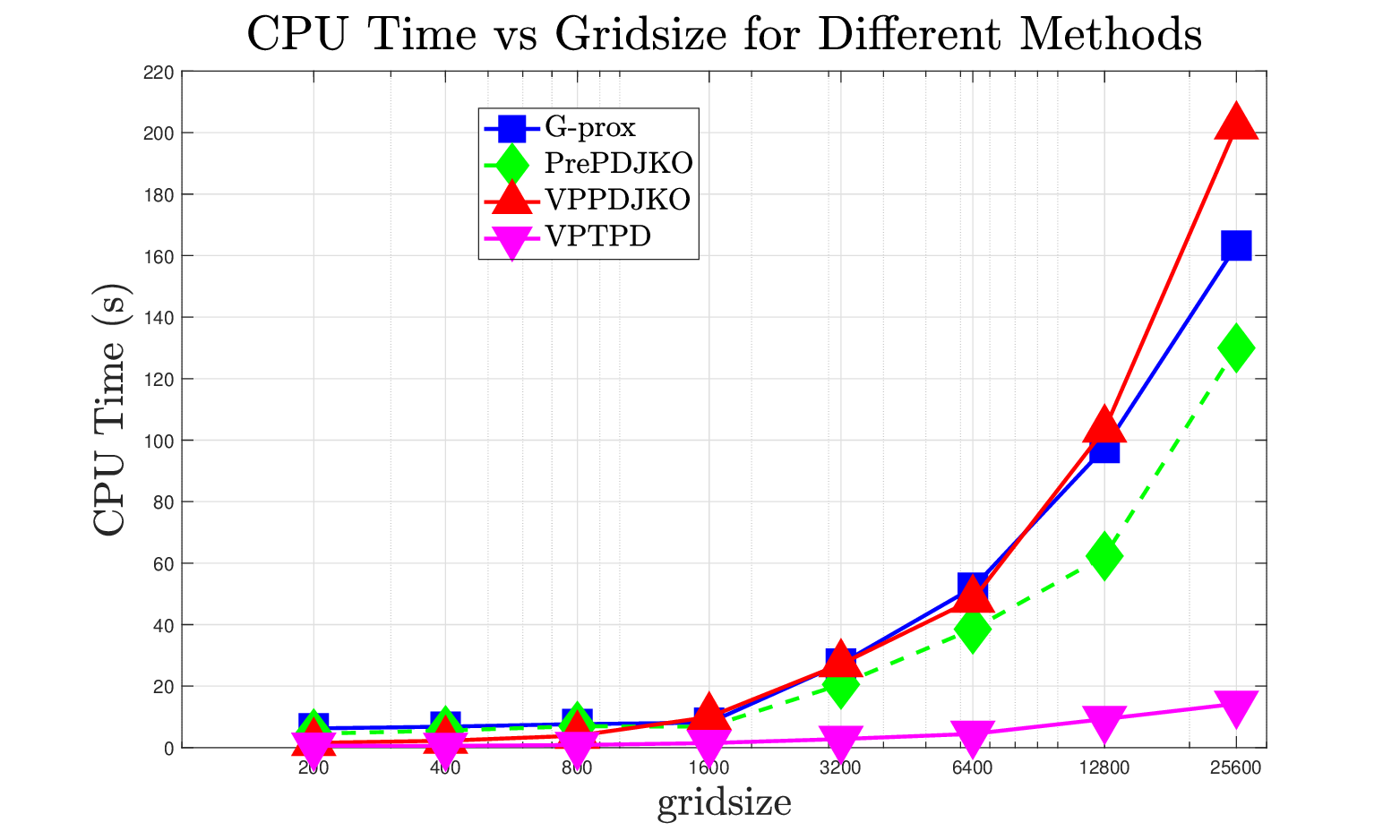}
\caption{Performance of VPTPD and PDHG-type methods on the porous medium equation under different grids with $\Delta t=0.05， \Omega=[-4, 4], T = 2$. For G-prox and PrePD: $\lambda = 0.8, 0.709, 0.69, 0.5, 0.65, 1.01, 2.5, 4$ for gridsizes $100\times 2^k$ ($k=1,\ldots,8$) respectively, and $\sigma = 1 / \lambda$. For VPPD: $\sigma = 0.5$ and $\lambda = 0.4$. For VPTPD: $\sigma = 0.3$ and $\lambda = 1$.}
\label{fig:performance_comparison}
\end{figure}

We compared the peformance of VPTPD with G-Prox PDHG, PrePDJKO, VPPDJKO on a sequence of refined grids. The results of iteration counts and CPU time are summarized in \autoref{fig:performance_comparison}. It can be observed that for porous medium equation, the G-prox PDHG and PrePDJKO methods exhibit similar trends in iteration counts. Under appropriate step sizes, their iteration counts remain essentially unchanged with grid refinement, indicating mesh-independent convergence. The CPU time of PrePDJKO is slightly lower, mainly because G-prox PDHG is a fully implicit iterative scheme and its proximal operator of $\mathcal{F}$ requires additional matrix-vector multiplication (details in \ref{Gprox_prox}). The iteration count for VPPDJKO gradually increases with mesh refinement and appears to eventually saturate at around gridsize 6400. VPTPD exhibits a similar trend, but requires substantially fewer iterations than VPPDJKO, evidently demonstrating the significant synergistic effect of the decoupling transformation and variable preconditioning in accelerating convergence.

In the third example, we implement a more challenging 1D Fokker-Planck equation with saturation effect by considering degenerate nonlinear mobility $M(\rho)=\rho(1-\rho)$, which indicates the global bounds of solutions $\rho\in[0,1]$, and logarithmic energy $\mathcal{H}=\rho\left(\ln\left(\rho\right)-1\right)$ with drift potential $V(x)=\frac{1}{2}\left\vert x\right\vert^{2}$. 
We solve the equation for $t \in [0, 15]$ with the initial value $\rho_{0}\left(x\right) = M/\mathrm{vol}\left(\Omega\right)$ ($\mathrm{vol}\left(\cdot\right)$ denotes the volume of $\cdot$) in the computational domain $\Omega = [-4, 4]$.

\begin{figure}[htbp]
	\centering  
	\includegraphics[width=0.32\textwidth]{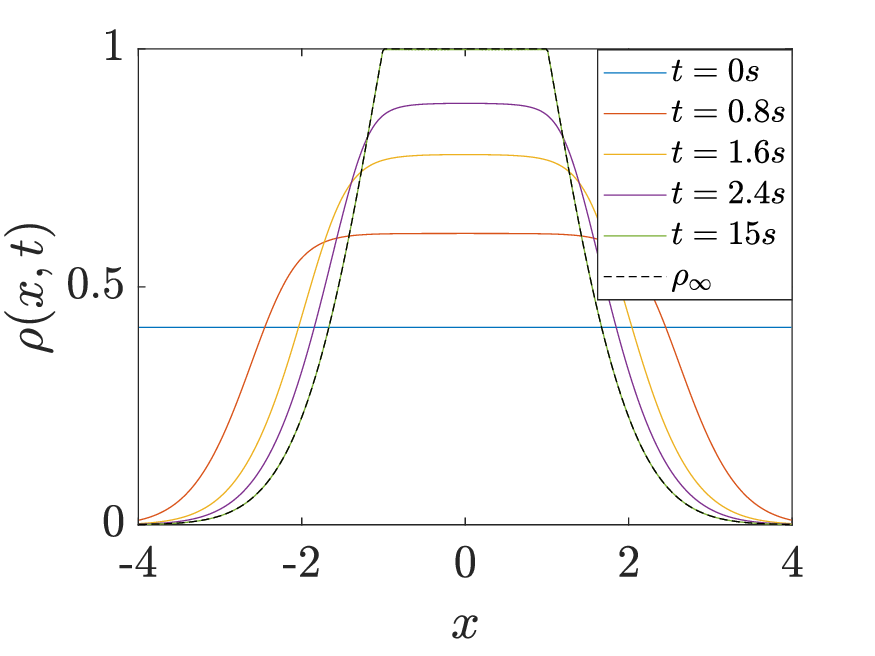}
	\includegraphics[width=0.32\textwidth]{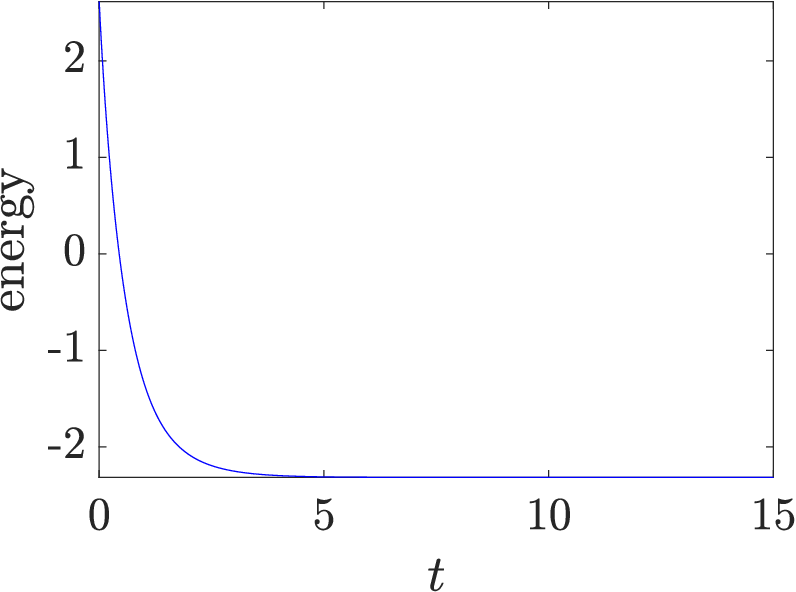}
	\includegraphics[width=0.32\linewidth]{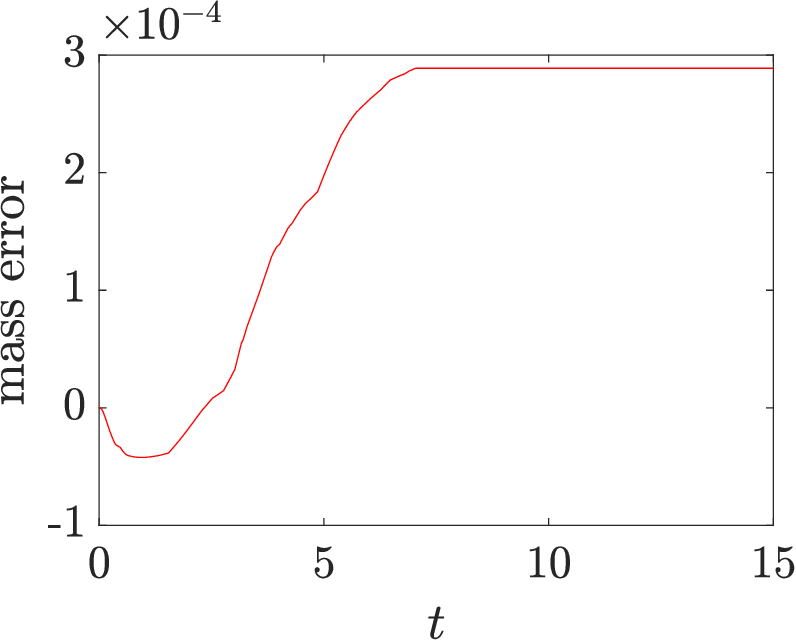}
	\caption{Evolution of solutions, energy and relative mass error for the 1D Saturation equation \eqref{1Dsaturation} for $t \in [0, 15]$. $M = 3.32,\; \Delta x = 0.02,\; \tau = 0.01, \; \lambda_0=0.05$. }
	\label{fig:main}
\end{figure}

The numerical results obtained by the VPTPD method are shown in \autoref{fig:main} (with the parameters given in the caption). The figure demonstrates that the numerical solution matches well with the analytical equilibrium solution $\rho_{\infty}$ \cite{bailo2023bound}, and more importantly the VPTPD method preserves the strict energy dissipation, global boundedness of solutions and mass conservation properties. We also compare the convergence and efficiency of the VPTPD method with the PDJKO and PrePDJKO methods in \autoref{tab:1D saturation}. Remarkably, the VPTPD method with constant step-size reduces the number of iterations by $99\%$ compared to PDJKO method. Compared to PrePDJKO method, the computational efficiency (in terms of CPU time) is improved by over 10 times for the VPTPD method with constant step-size and more than 10 times with the adaptive step-size strategy.

\begin{table}[htb]
	\centering
	\small
	\begin{tabular}{ccccc} 
		\hline
		\multicolumn{5}{c}{\textbf{1D Saturation equation}}\\
		\hline
		& PDJKO & PrePDJKO & VPTPD & VPTPD(s)\\ 
		\hline
		Number of iterations  & 8018525  & 165411(-97.9\%)       & 81198(-99.0\%)    &19500(-99.8\%)   \\
		CPU time  & 1000   & 51(-95.0\%)       & 13(-98.7\%)    &4(-99.6\%)  \\ 
		\hline
	\end{tabular}
	\caption{Number of iterations and CPU time for solving 1D saturation equation \eqref{1Dsaturation} for $t \in [0, 15]$. $N_{x} = 400,\;\tau = 0.01, \ \lambda_0=0.05$.  For adaptive step-size: $\lambda_0=1,\;\lambda_{max}=100,\;\lambda_{min}=0.0001,\;\gamma_{+}=10,\;\gamma_{-}=16$}
	\label{tab:1D saturation}
\end{table}

\subsubsection{2D Cahn-Hilliard equation} \label{2D-CH}
Now, we consider the 2D Cahn-Hilliard (CH) equation with nonlinear mobility:
\begin{equation}\label{2DCH_random}
	\rho_{t} = \nabla\cdot\left(M\left(\rho\right)\left(\mathcal{H}'\left(\rho\right) - \epsilon^{2}\Delta\rho\right)\right),
\end{equation}
where $M\left(\rho\right) = 1 - \rho^{2}$ and $\mathcal{H}\left(\rho\right) = \dfrac{1}{4}\left(1 - \rho^{2}\right)^{2}$. The corresponding energy $\mathcal{J}\left(\rho\right)$ is $\displaystyle\int_{\Omega}\left(\mathcal{H}\left(\rho\right) + \dfrac{\epsilon^{2}}{2}\left\vert\nabla\rho\right\vert^{2}\right)\mathrm{d}x$ with $\epsilon = 0.018$.

\begin{figure}[htb]
	\centering  
	\includegraphics[width=0.19\textwidth]{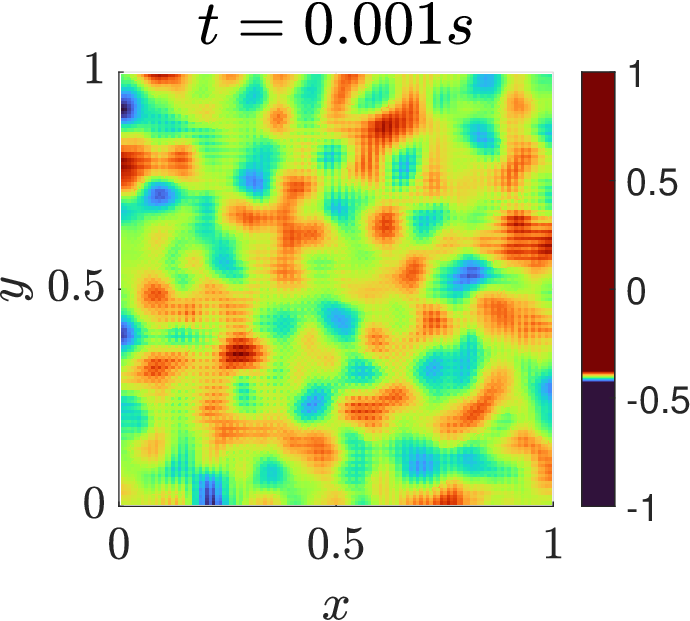}
	\includegraphics[width=0.19\textwidth]{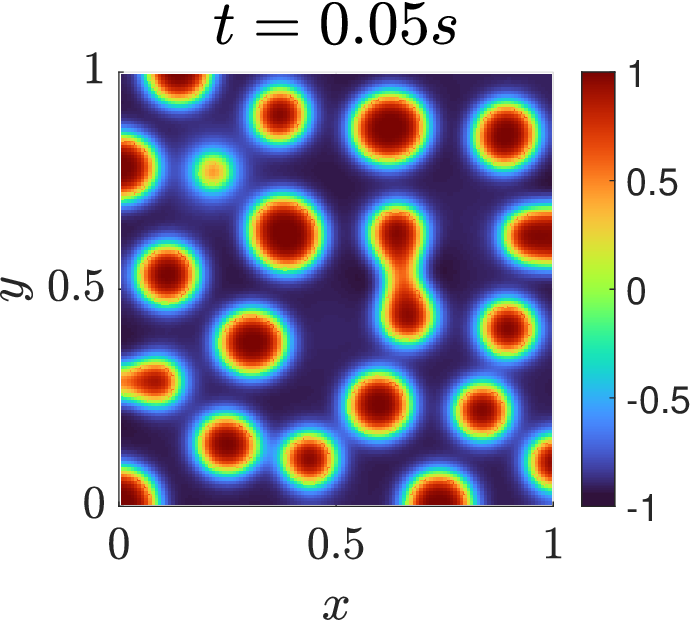}
	\includegraphics[width=0.19\linewidth]{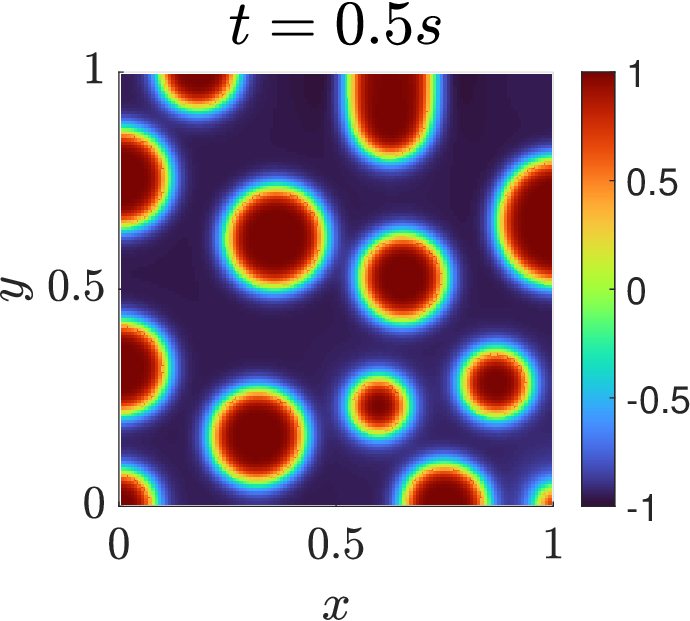}
	\includegraphics[width=0.19\linewidth]{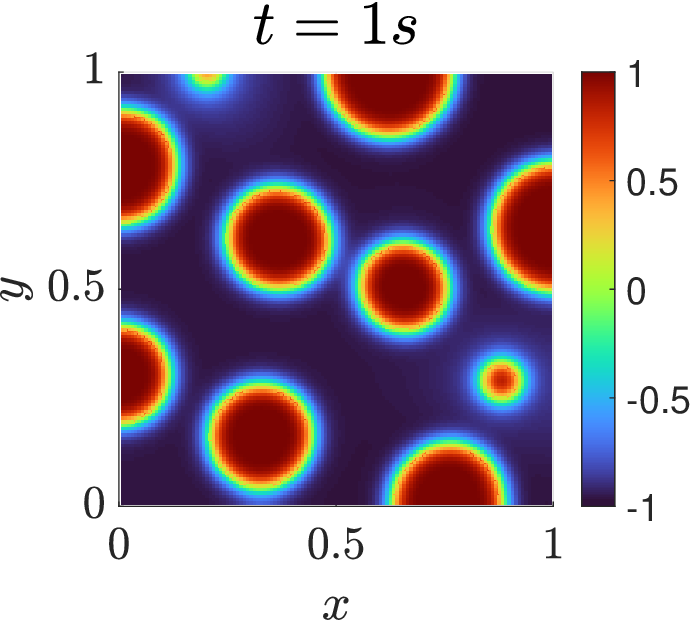}
	\includegraphics[width=0.19\linewidth]{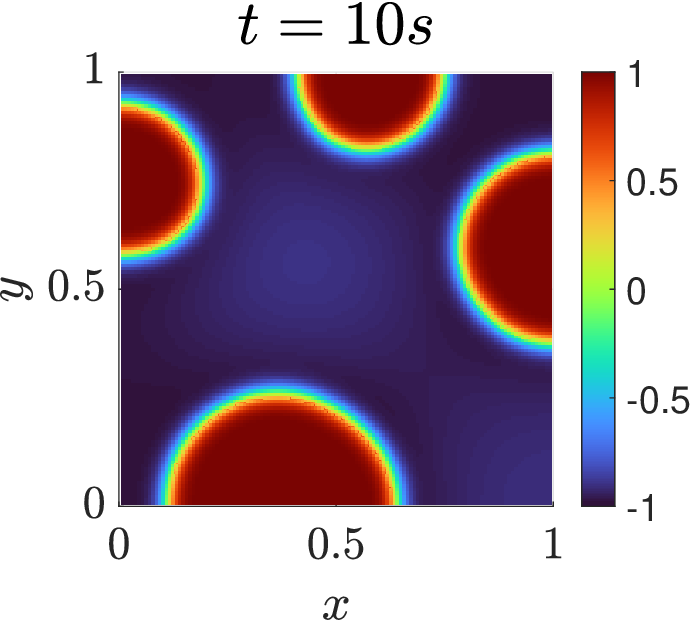}\\
	\includegraphics[width=0.35\linewidth]{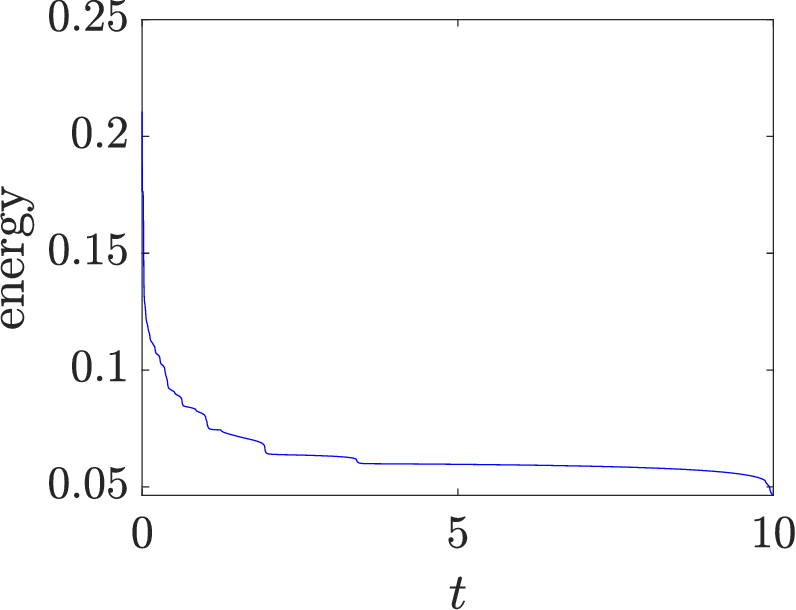}
	\includegraphics[width=0.35\linewidth]{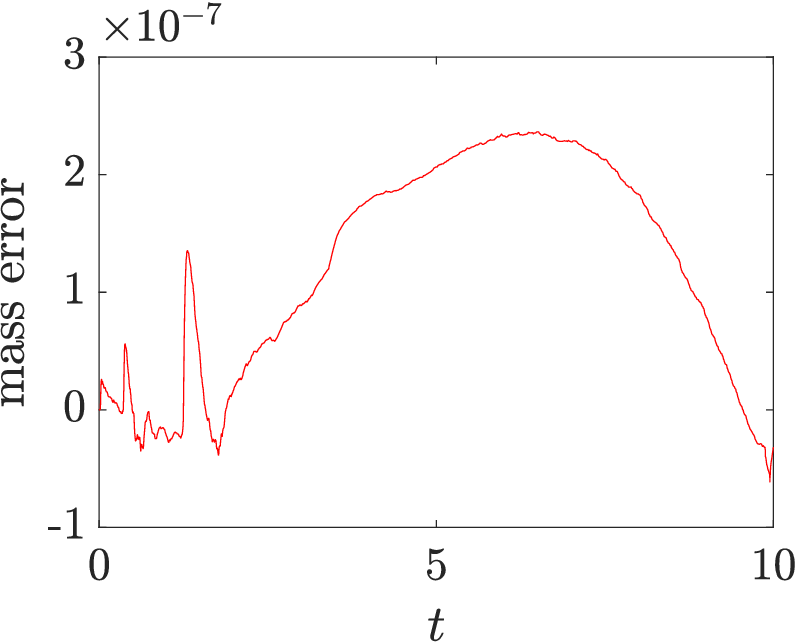}
	\caption{Evolution of solutions, energy and relative mass error for the 2D CH equation \eqref{2DCH_random} with perturbed initial value for $t \in [0, 10]$. $N_{x} = 128, \; N_{y} = 128, \; \tau = 0.001$.}
	\label{fig:CH_random}
\end{figure}

In the first experiment, we consider 2D phase-separation with a perturbed initial value to test the computational efficiency of VPTPD. The initial value is given by $\rho_{0}\left(x,y\right) = -0.4 + c$ over the domain $\Omega = [0, 1] \times [0, 1]$, where $c$ follows a uniform distribution in $[-0.1, 0.1]$. The numerical results obtained by the VPTPD method are shown in \autoref{fig:CH_random}. The computational efficiencies of different methods are presented in \autoref{tab:2D CH1}. The number of iterations and CPU time is significantly reduced by around $75\%$ for the VPTPD method compared to the PrePDJKO method. The adaptive strategy for VPTPD(s) method can further improve the computational efficiency of VPTPD method by $10\%$.  

\begin{table}[htbp]
	\centering
	\footnotesize
	\begin{tabular}{ccccc} 
		\hline
		\multicolumn{5}{c}{\textbf{2D Cahn-Hilliard equation}}\\
		\hline
		JKO steps & & PrePDJKO & VPTPD & VPTPD(s)\\ 
		\cline{2-5} 
		\hline
		\multirow{2}{*}{10} & Number of iterations  &9041   & 1081(-88.0\%) & 957(-89.4\%)\\
		& CPU time  & 88   & 10(-88.6\%) & 8.7(-90.1\%)\\
		\hline
		\multirow{2}{*}{500} & Number of iterations  & 153688  & 41602(-72.9\%) & 31253(-79.7\%)    \\
		& CPU time  & 1894   & 595(-69.6\%)       & 472(-75.1\%)\\
		\hline
	\end{tabular}
	\caption{Number of iterations and CPU time of different methods for solving 2D CH equation \eqref{2DCH_random} with perturbed initial value. $\;N_{x} = 128,\;N_{y} = 128,\;\tau = 0.001,\;\lambda_{0}=1$}
	\label{tab:2D CH1}
\end{table}

\begin{figure}[htbp]
\centering 
\begin{minipage}[c]{0.49\textwidth}
\centering
\includegraphics[width=0.49\linewidth]{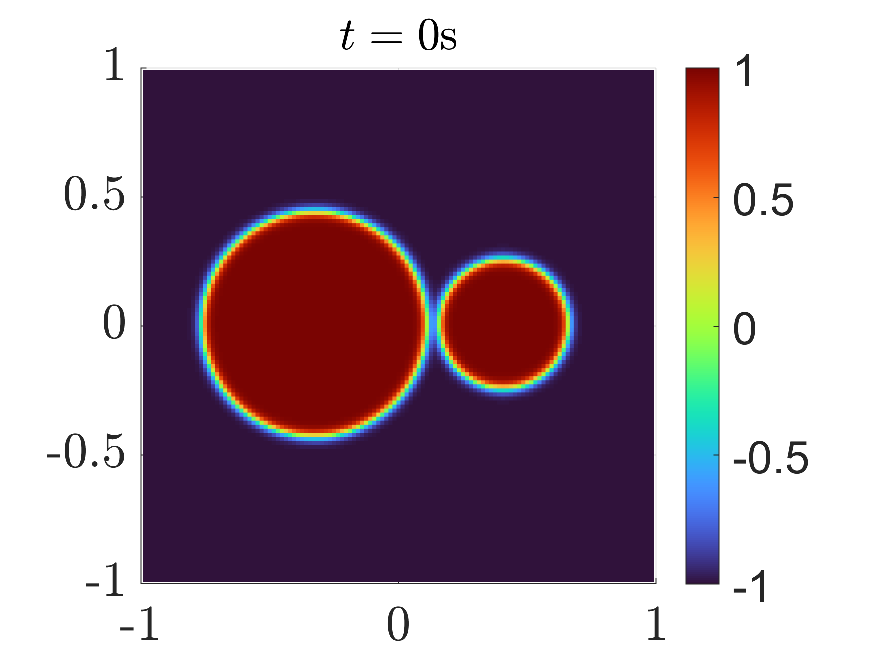}
\includegraphics[width=0.49\linewidth]{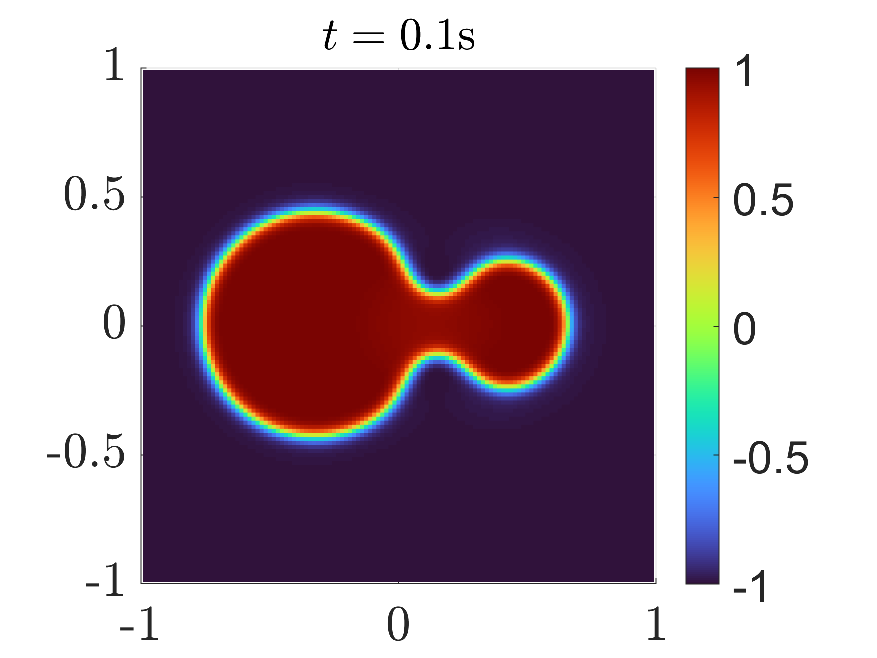}\\
\includegraphics[width=0.49\linewidth]{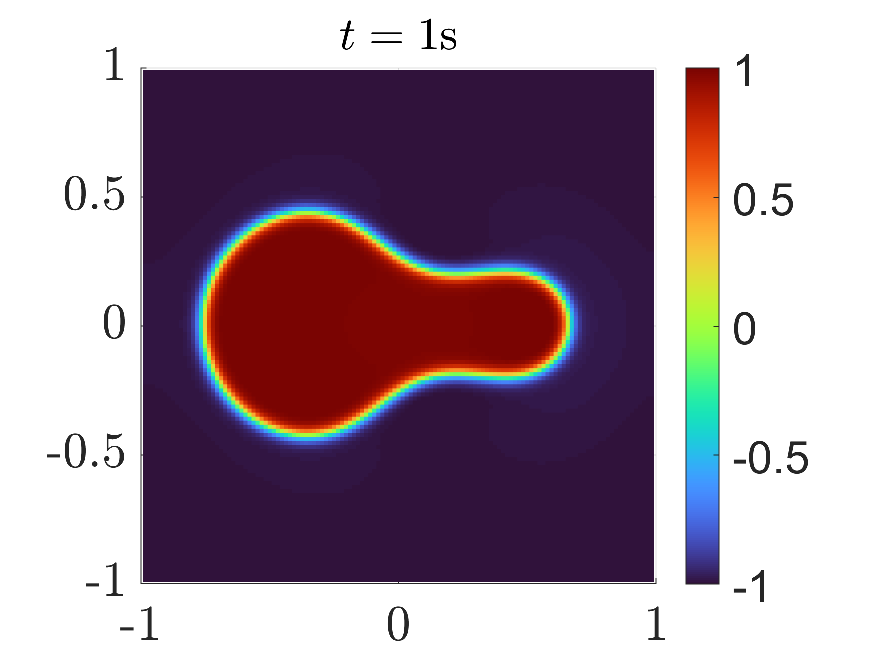}
\includegraphics[width=0.49\linewidth]{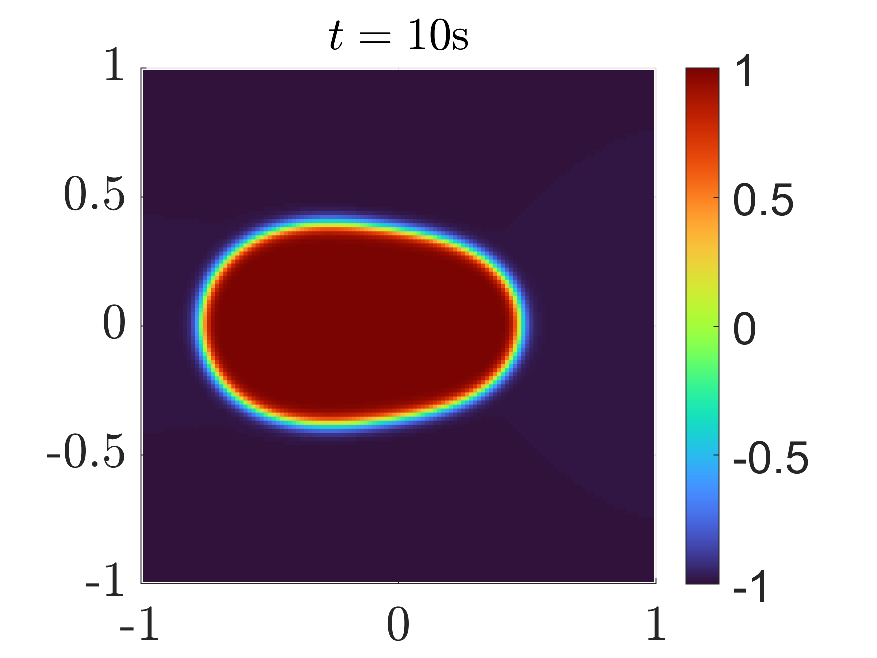}
\end{minipage}
\hfill
\begin{minipage}[c]{0.48\textwidth} 
\centering
\includegraphics[width=\linewidth]{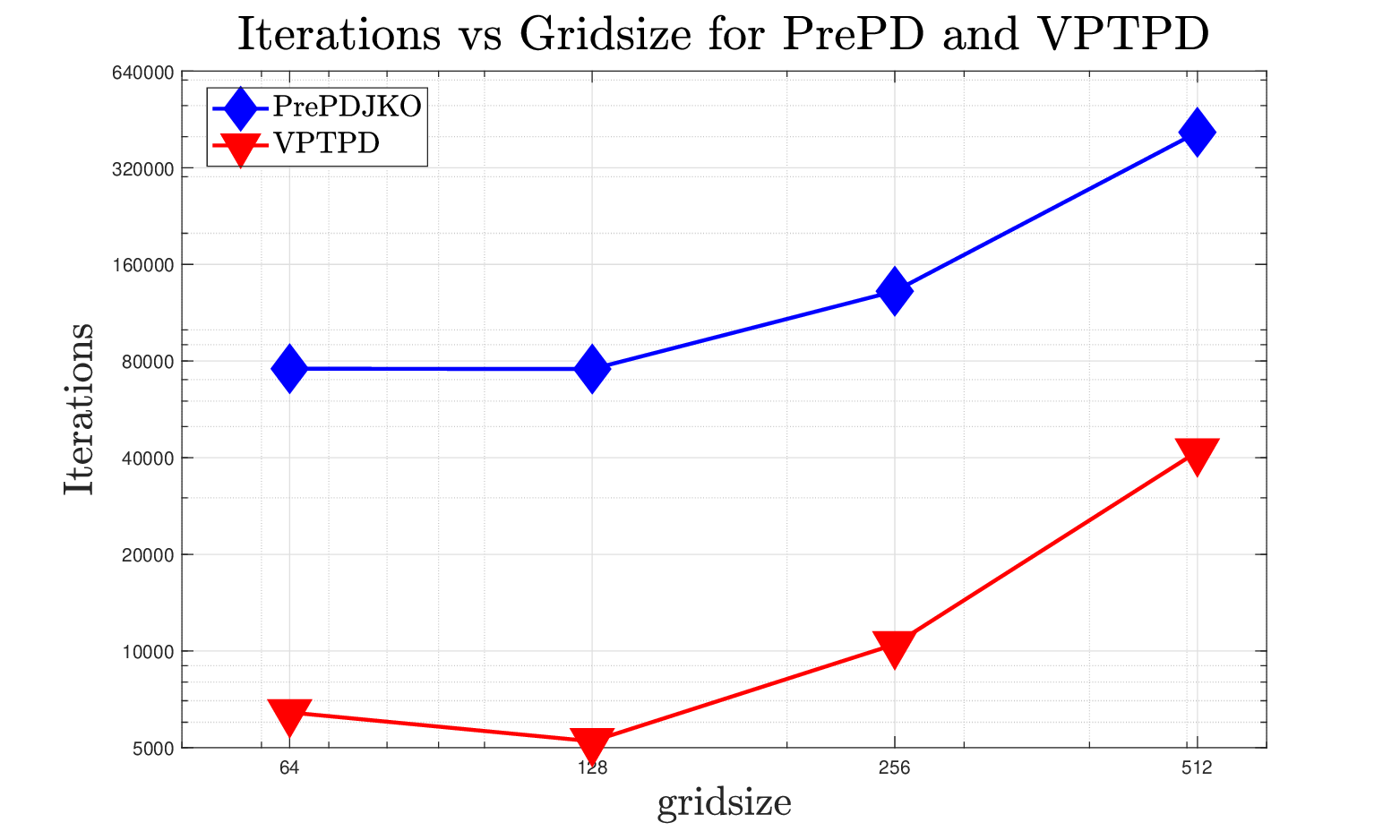}
\end{minipage}
\caption{Evolution of solutions (left) for the 2D CH equation \eqref{2DCH_random} with initial value \eqref{initial_merge} for $N_{x} = 128, \; N_{y} = 128, \; \tau = 0.01$. Iterations of VPTPD and PrePDJKO for 100-step JKO with various grid sizes (right). }
\label{merge}
\end{figure}

Next, we consider a 2D droplet-merging example to test how iteration counts change with grid size.
The initial condition over $\Omega = [-1,1]^{2}$ is given by
\begin{equation}\label{initial_merge}
\rho_0(x,y)=\tanh\Bigl(\frac{r_{1} - \sqrt{(x-c_{11})^{2}+(y-c_{12})^{2}}}{\varepsilon}\Bigr)
            + \tanh\Bigl(\frac{r_{2} - \sqrt{(x-c_{21})^{2}+(y-c_{22})^{2}}}{\varepsilon}\Bigr) + 1,
\end{equation}
where $c_{11} = -\frac{1}{3},\; c_{12} = 0,\; r_{1} = \frac{1}{2.3},\; 
c_{21} = \frac{2}{5},\; c_{22} = 0,\; r_{2} = \frac{1}{4}$. The evolution of solutions for $t\in[0, 10]$ and the iteration counts v.s. grid number are presented in \autoref{merge}. We observe that the iteration count for VPTPD is approximately one-tenth that of PrePDJKO, but both methods exhibit a similar dependence on the grid resolution. Compared to the behavior for 1D porous medium equation in \autoref{fig:performance_comparison}, the significant mesh-dependence for this 2D Cahn-Hilliard equation reflects the complexity of the convergence behavior of VPTPD method, whose complete understanding is highly nontrivial and deserves further investigation in a dedicated future study.

Furthermore, to investigate the impact of the regularization \eqref{modify} on preconditioning, we compare the iteration counts for 100 JKO-steps for different regularization parameter $r$ in \autoref{parameter_r}. 
The numerical evidence shows that iteration counts exhibit small dependence on the value of regularization parameter $r$, demonstrating the robustness of the regularization approach for constructing effective preconditioners.

\begin{table}[htpb]
\centering
\noindent Parameters: $\Omega = [-1,1]^{2} \times [0,1]$, $T=1$, $\tau=0.01$, $N_x = N_y = 64$

\vspace{1em}
\begin{tabular}{ccccccc}
\toprule
$r$ & $10^{-3}$ & $10^{-4}$ & $10^{-5}$ & $10^{-6}$ & $10^{-7}$ & $10^{-8}$\\
\midrule
Number of iterations & 6437 & 6453 & 6846 & 7241 & 7180 & 7193\\
\bottomrule
\end{tabular}
\caption{Number of iterations for 2D droplet-merging experiment under different regularization parameters.}\label{parameter_r}
\end{table}

\subsubsection{3D wetting phenomenon of droplet}
In this section, we consider the wetting phenomenon of 3D droplets with moving contact line described by \cite{carrillo2024structure}:
\begin{align}\label{3D_dewetting}
	\begin{cases}
		$$\rho_{t} = -\nabla\cdot\left(M\left(\rho\right)\left(-\nabla\dfrac{\delta }{\delta \rho}\mathcal{J\left(\rho\right)}\right)\right),$$\\
		$$\displaystyle\mathcal{J\left(\rho\right)} = \int_{\Omega}\mathcal{H}\left(\rho\right)\mathrm{d}x +  \dfrac{\epsilon^2}{2}\int_{\Omega}\vert \nabla\rho\vert^{2}\mathrm{d}x  + \int_{\mathcal{T}_{w}}f_{w}\left(\rho,\beta_{w}\right)\mathrm{d}s,$$
	\end{cases}
\end{align}
where $M\left(\rho\right) = 1 - \rho^{2}$, $\mathcal{H}\left(\rho\right) = \dfrac{1}{4}\left(1 - \rho^{2}\right)^{2}$, and $f_w(\rho,\beta_w)=\frac{\epsilon}{\sqrt{2}}\cos\beta_w(\rho^3/3-\rho)$ is the cubic-polynomial wall energy \cite{lee2011accurate} defined on the substrate $\mathcal{T}_{w}\subset\partial\Omega$, where $\beta_{w}$ is the equilibrium contact angle between the free interface and the substrate. In this example, we consider the domain $\Omega = [-0.5, 0.5]\times[-0.5, 0.5]\times[0, 0.4]$ and take the interfacial parameter $\epsilon = 0.01$, the contact angle $\beta_{w} = \pi/4$ and $\beta_{w} = 3\pi/4$, respectively. The initial value is taken as a semi-sphere:
\begin{equation}
	\rho_{0}\left(x,y,z\right) = -\tanh{\dfrac{\sqrt{x^{2} + y^{2} + z^{2}} - 0.25}{\sqrt{2}\epsilon}}.
\end{equation}

\begin{figure}[htbp]
	\centering  
	\includegraphics[width=0.3\textwidth]{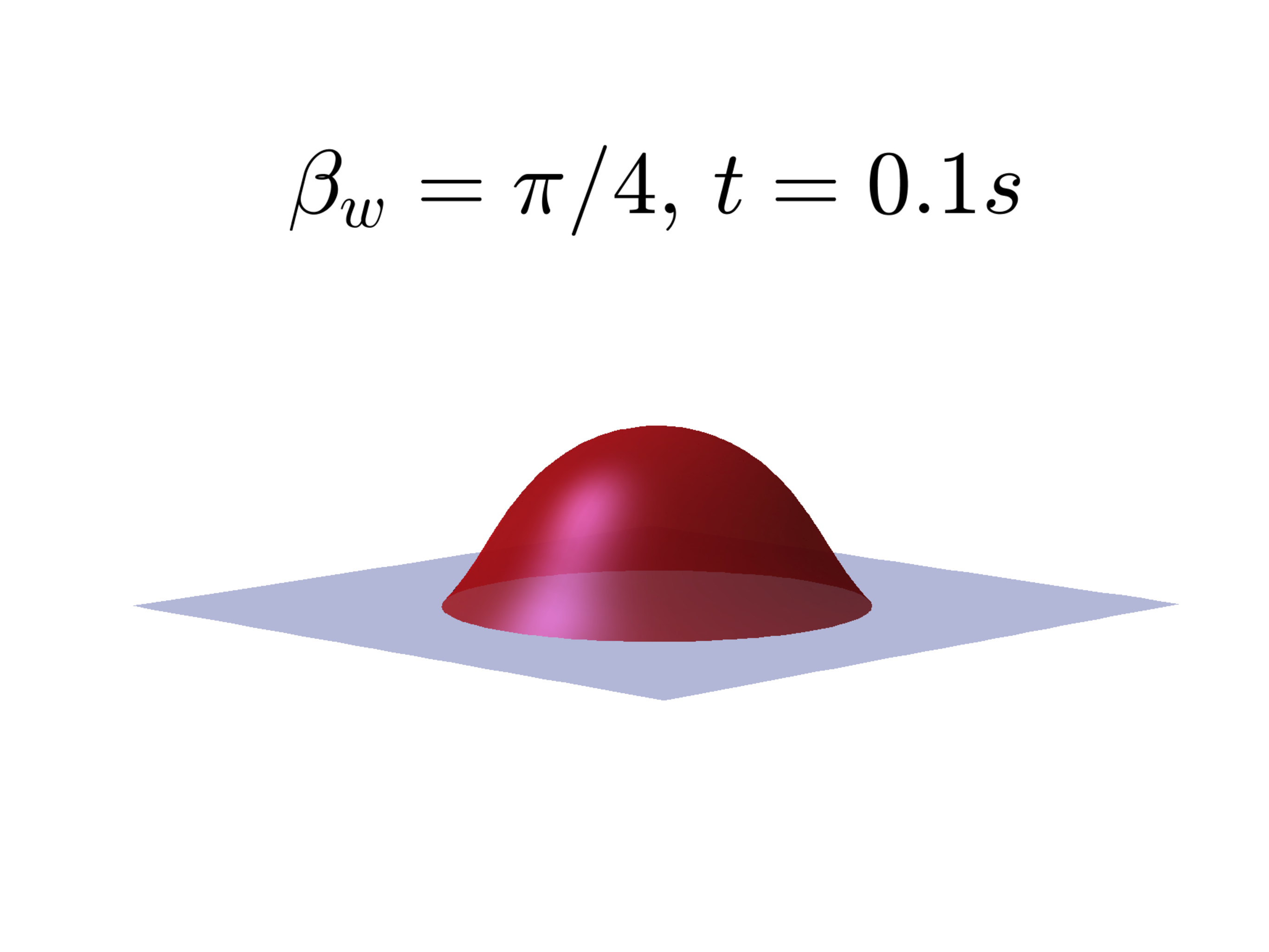}
	\includegraphics[width=0.3\textwidth]{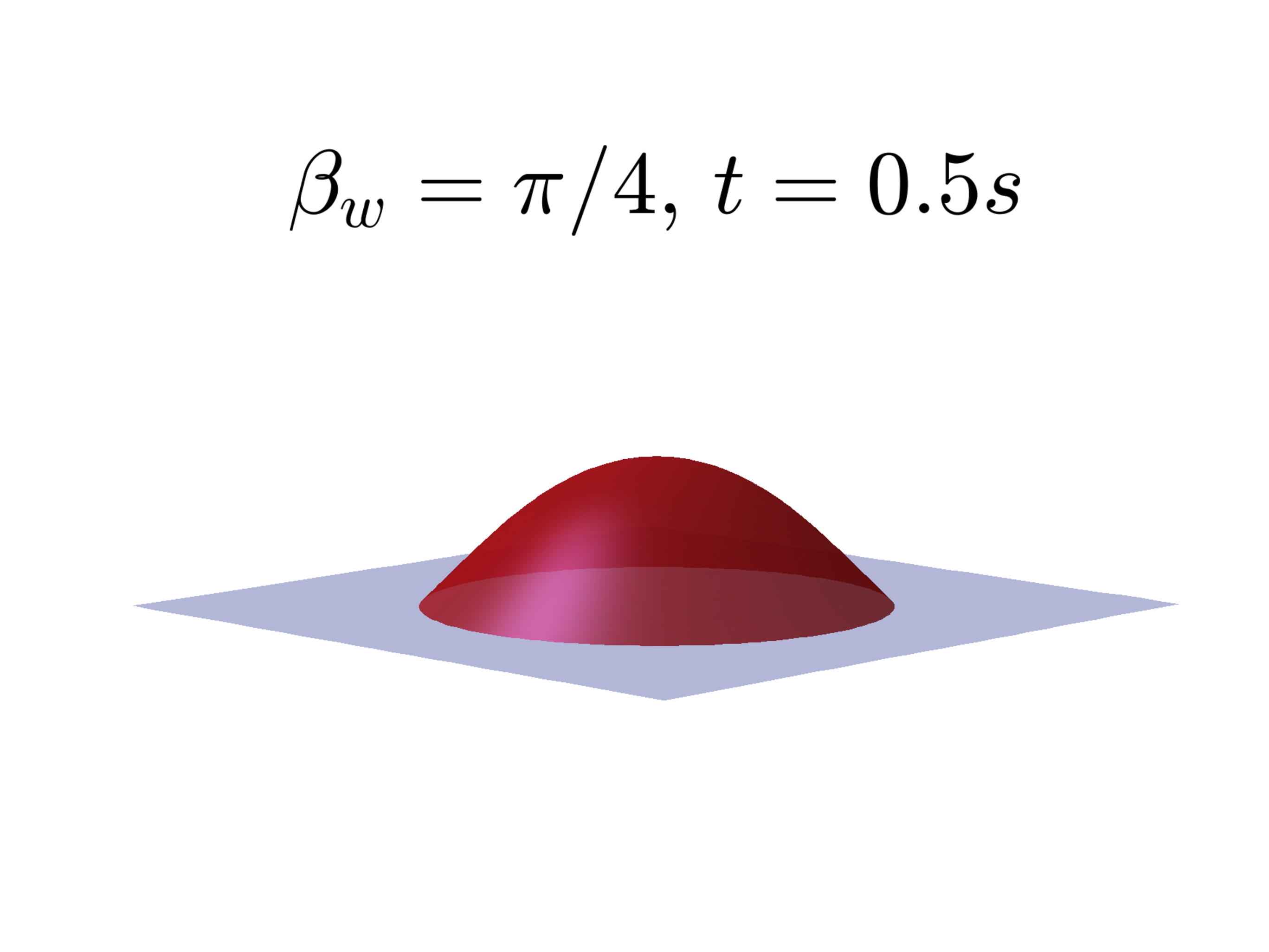}
	\includegraphics[width=0.3\linewidth]{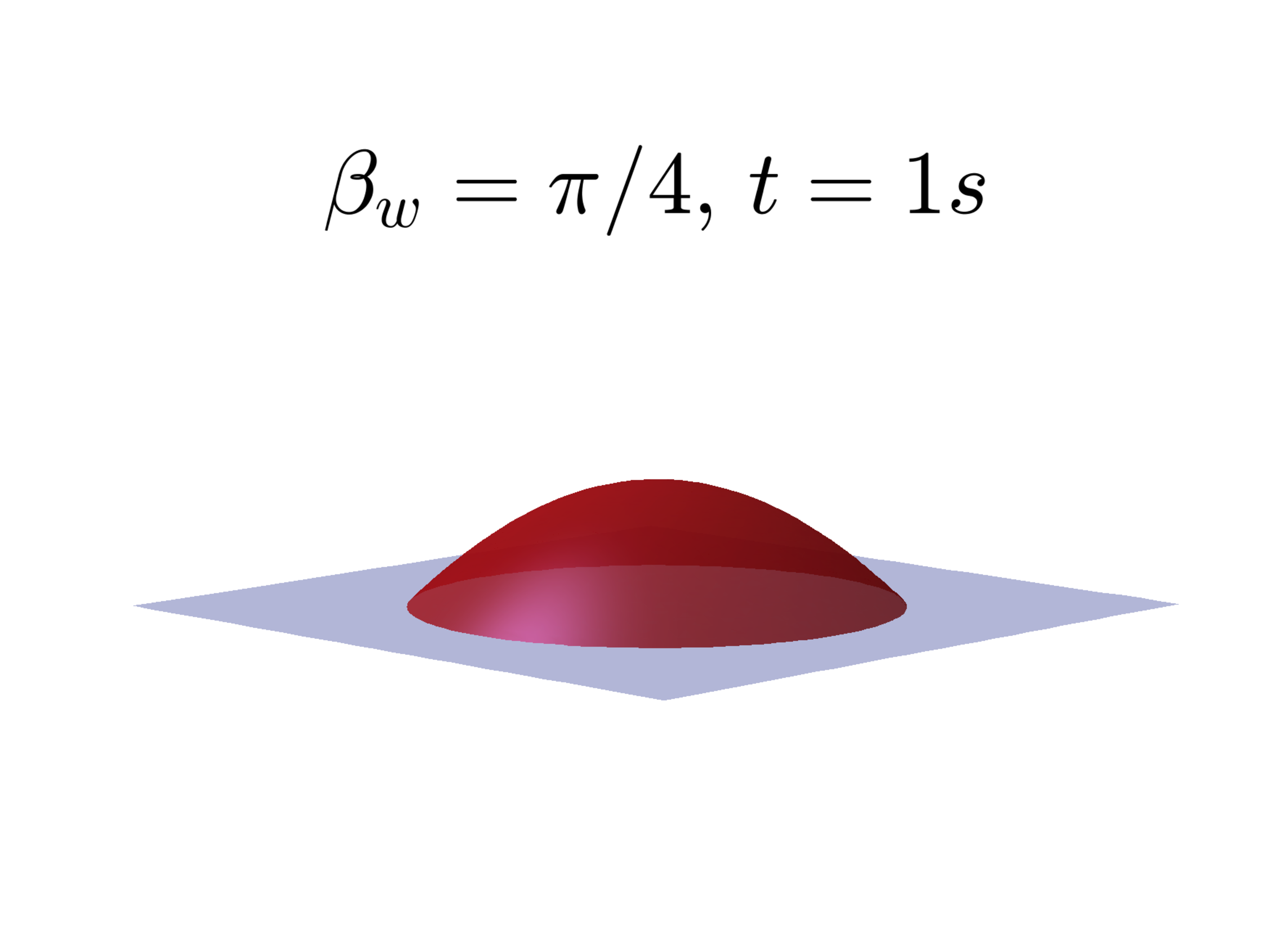}\\
	\includegraphics[width=0.3\linewidth]{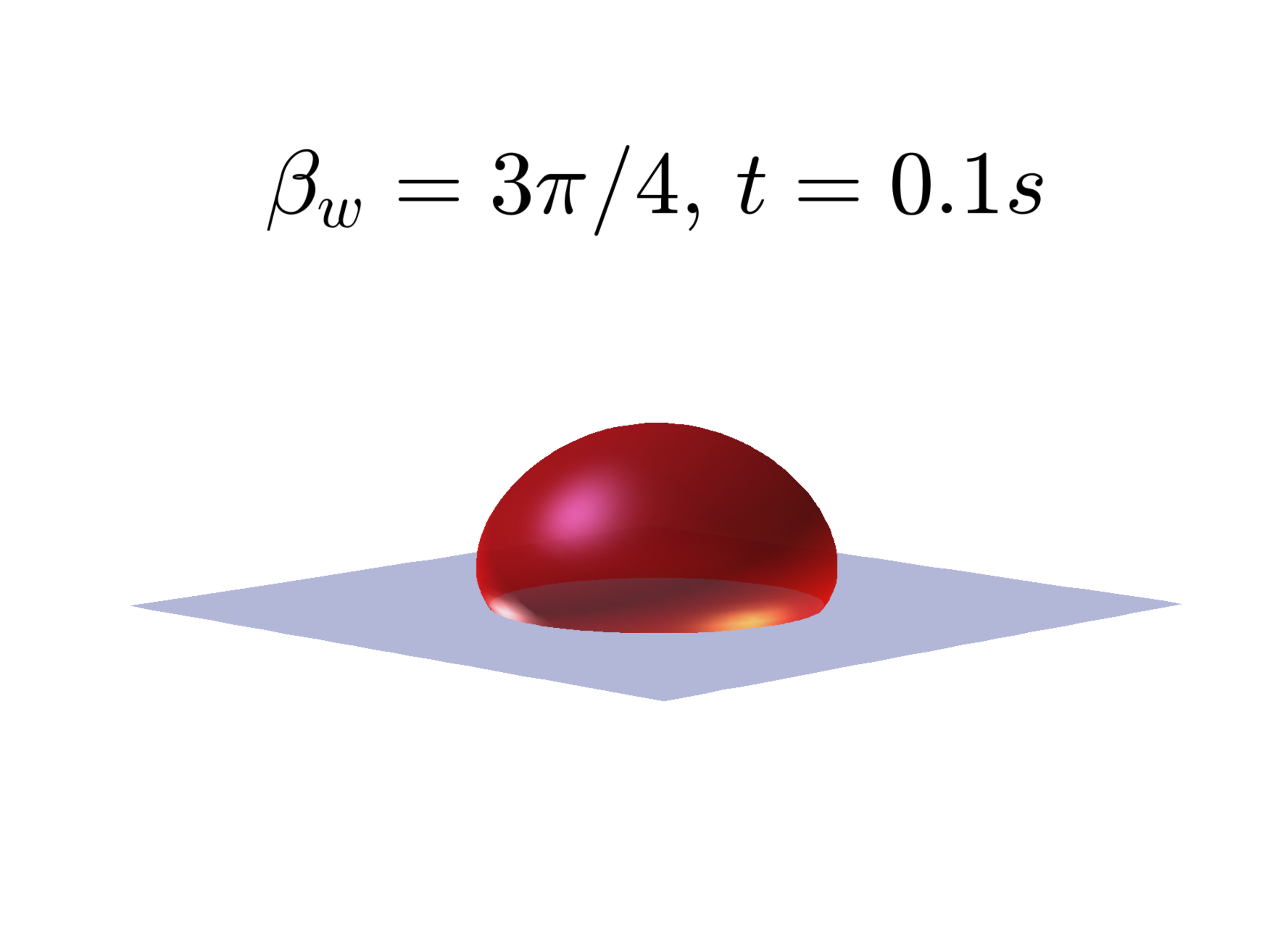}
	\includegraphics[width=0.3\linewidth]{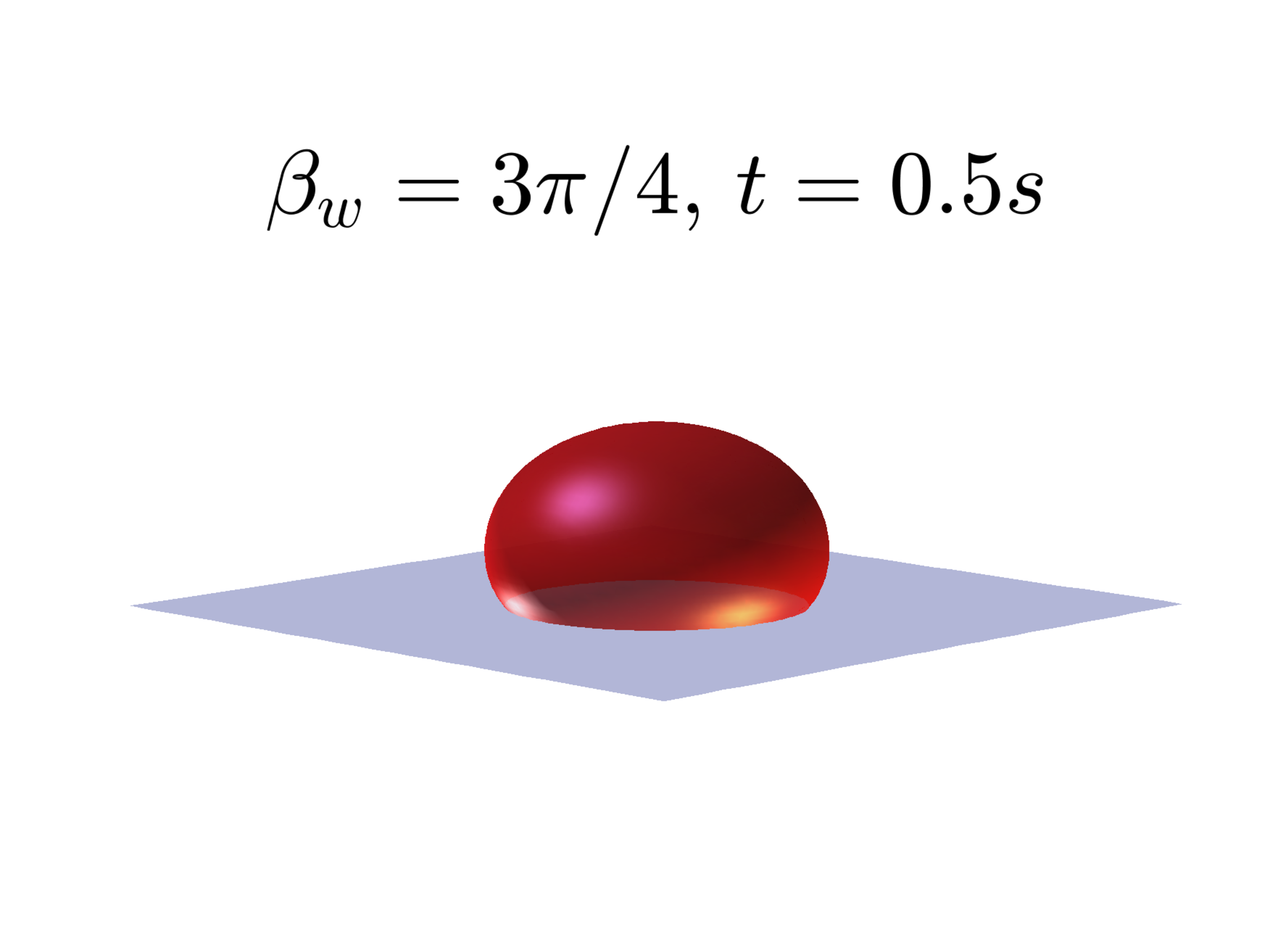}
	\includegraphics[width=0.3\linewidth]{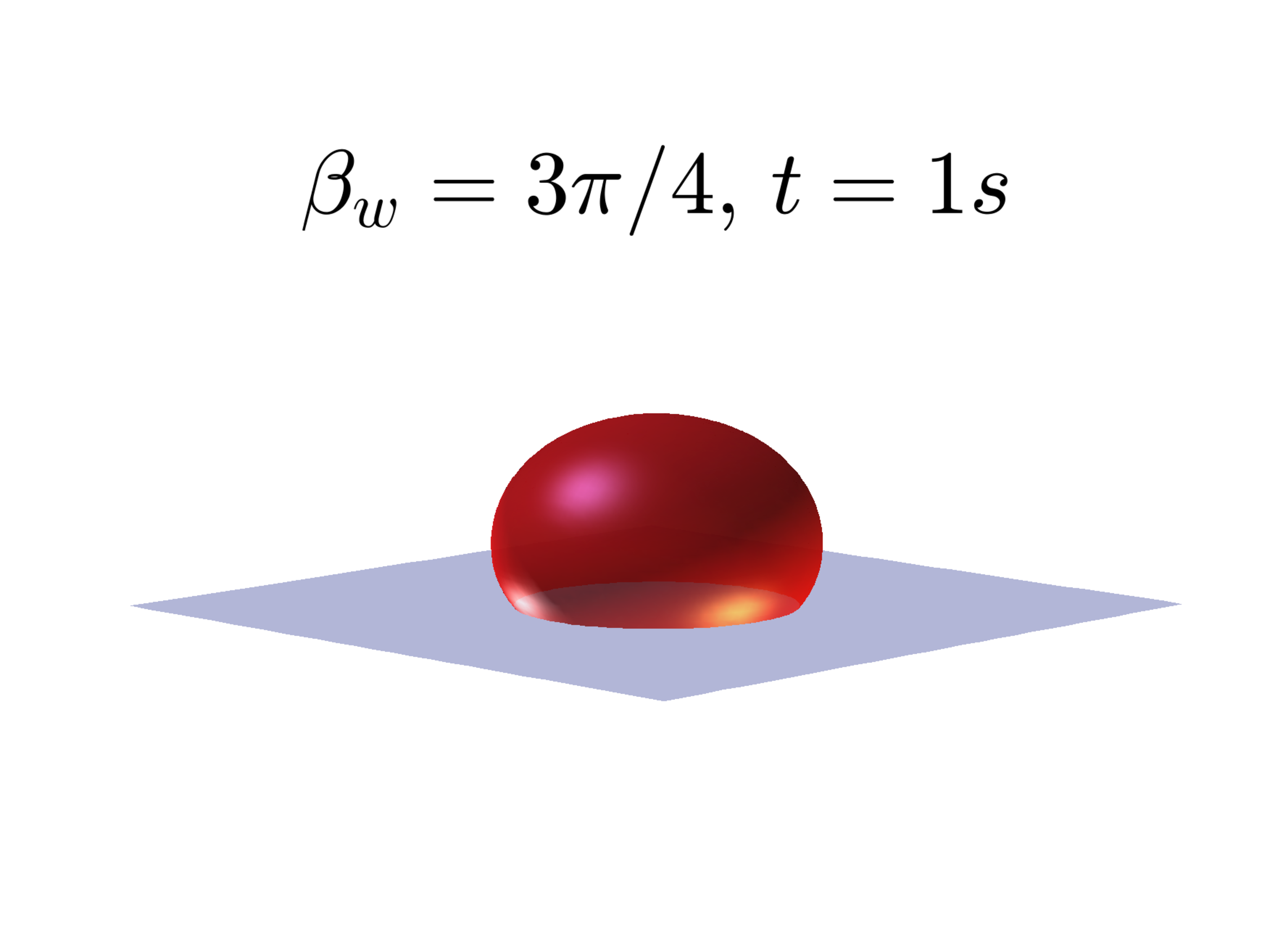}
	\caption{The shape evolution of 3D droplet wetting with $\beta_{w} = \pi/4$ (top) and $3\pi/4$ (bottom) for $t \in [0, 1]$ and $\Omega = [-0.5,0.5]^{2} \times [0,0.4]$. $N_{x} = 64,\; N_{y} = 64,\; N_{z} = 40,\;\tau = 0.1$.}
	\label{fig:3D wetting phenomenon of droplet}
\end{figure}

\begin{figure}[htbp]
	\centering  
	\includegraphics[width=0.35\textwidth]{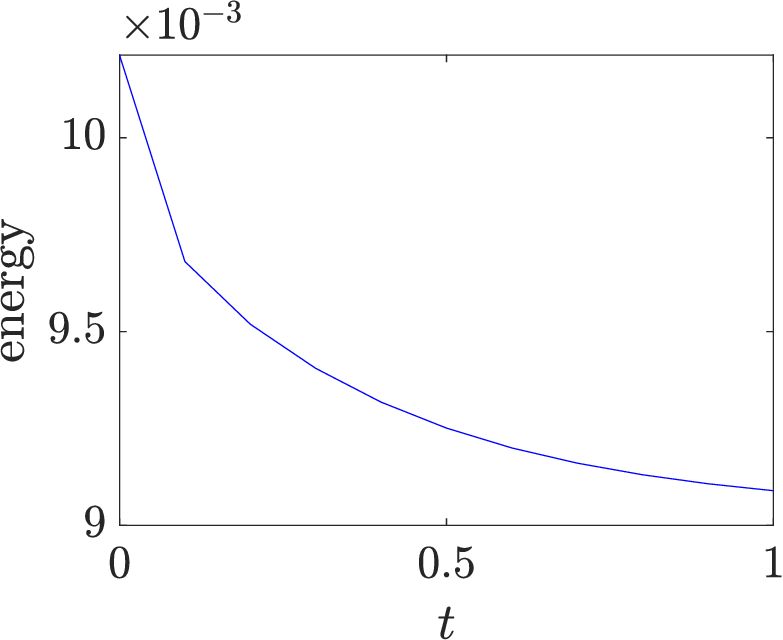}
	\includegraphics[width=0.35\textwidth]{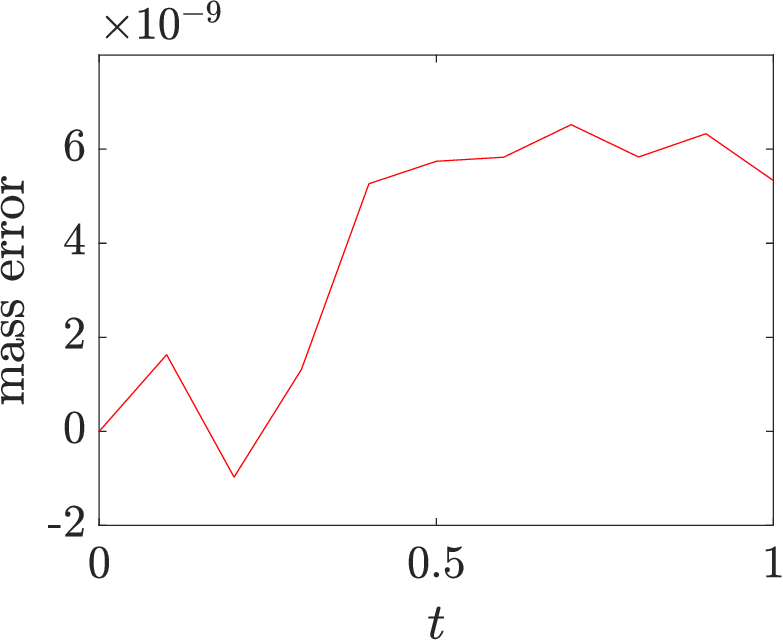}
	\caption{Evolution of energy and relative mass error for 3D droplet wetting for $t \in [0, 1]$, $\beta_{w}=\pi/4$. $N_{x} = 64,\; N_{y} = 64,\; N_{z} = 40,\;\tau = 0.1$.}
	\label{fig:energy and mass error of 3D wetting phenomenon of droplet}
\end{figure}

We simulate the wetting of droplets for different contact angles with a quite large time step $\tau=0.1$, and plot the shape evolution (outlined by the contour of $\rho = 0$) in \autoref{fig:3D wetting phenomenon of droplet}. The energy dissipation and mass conservation is demonstrated in \autoref{fig:energy and mass error of 3D wetting phenomenon of droplet}. 
The number of iterations and CPU time of different methods are shown in \autoref{3D_dewetting1}. Notably the VPTPD method exhibits significantly faster convergence (with $99.2\%$ less iterations) compared to the PrePDJKO method in this 3D example, which indicates the increasingly significant improvement by VPTPD method in the convergence rate as the dimension of problem increases. Additionally, the adaptive step-size strategy still shows a significant acceleration effect. We note that the save rate of CPU time by VPTPD method is $94.5\%$, which is already outstanding but not as high as that of iterations due to the computational bottleneck of inverting large matrices. We would expect the VPTPD could turbocharge the computational efficiency enhancement for large systems combined with high-performance computing techniques (e.g., parallel algorithms).

\begin{table}[htbp]
	\centering
	\small
	\begin{tabular}{ccccc} 
		\hline
		\multicolumn{5}{c}{\textbf{3D wetting experiment for $\beta_{w}=\pi/4$}}\\
		\hline
		\multirow{1}{*}{JKO steps}	& & PrePDJKO & VPTPD & VPTPD(s)\\ 
		\cline{2-5}
		\hline
		\multirow{2}{*}{1}	& Number of iterations  &66386   &1334(-98.0\%) & 361(-99.5\%)   \\
		&CPU time  & 20717   & 2800(-86.5\%)& 750(-96.4\%) \\
		\hline
		\multirow{2}{*}{10}	 & Number of iterations  & 568400  & 4445(-99.2\%) & 3392(-99.4\%)   \\
		&CPU time  &177380   & 9280(-94.5\%)    & 6940(-96.1\%)\\
		\hline
	\end{tabular}
	\caption{Number of iterations and CPU time of different methods for simulating 3D droplet wetting \eqref{3D_dewetting} for $\beta_{w}=\pi/4$.$\;N_{x} = 64,\;N_{y} = 64,\;N_{z}=40,\;\tau = 0.1,\;\lambda_{0}=0.5$.}
	\label{3D_dewetting1}
\end{table}

\begin{figure}[hbtp]
\centering
\includegraphics[width=0.45\textwidth]{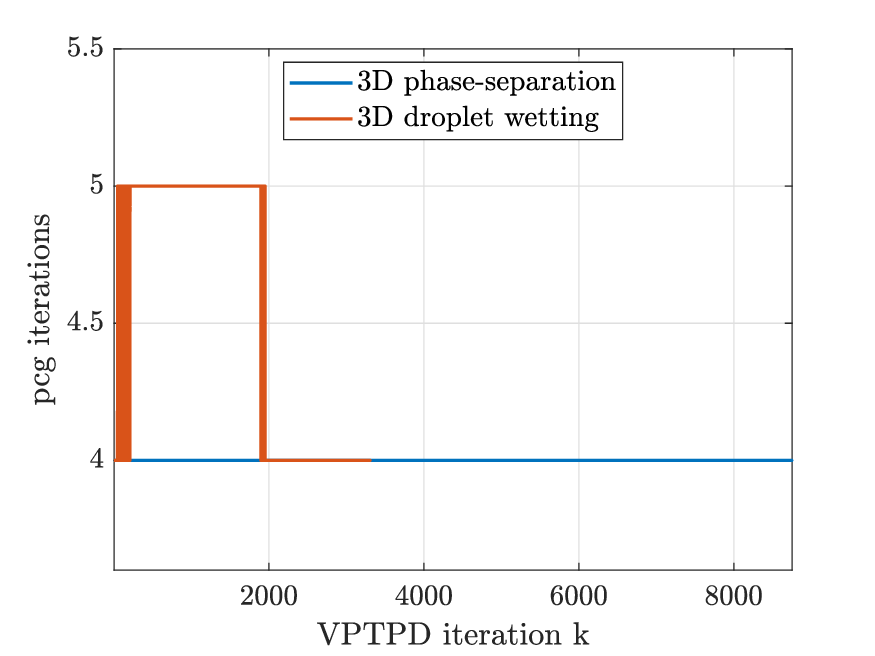}
\caption{Number of PCG iterations per VPTPD iteration in 10 JKO steps for 3D phase-separation simulation (random initial) with $N_x = 60, N_y = 60, N_z = 60$ and 3D droplet wetting ($\beta=\pi/4$) with $N_x = 100, N_y = 100, N_z = 40$.}
\label{variable mobility}
\end{figure}

Moreover, we further investigate the computational burden of applying $I^{-1}_{p}$ in 3D simulations, measured by the PCG iteration counts. The results show that the PCG iteration count is consistently bounded by a small number ($\leq 5$) across a wide range of grid resolutions (\autoref{tab:pcg_iterations}) and mobility values (\autoref{variable mobility}). This demonstrates the robustness of the proposed preconditioner and highlights the strong scalability and efficiency of the VPTPD method.

\begin{table}[ht]
\centering
\begin{tabular}{ccccc}
\toprule
$N_x:N_y:N_z$ & $32:32:20$ & $64:64:40$ & $100:100:40$ & $150:150:60$ \\
\midrule
Average PCG iterations & 2.01 & 4.63 & 4.45 & 4.75 \\
\bottomrule
\end{tabular}
\caption{Average PCG iterations for 3D droplet wetting simulations under different grid resolutions $N_x:N_y:N_z$.}
\label{tab:pcg_iterations}
\end{table}

\subsection{Extended numerical experiments}\label{extended example}
A suite of extended numerical experiments are performed to demonstrate the flexibility and applicability of our VPTPD method to several challenging problems.

\subsubsection{1D Keller-Segel equation}
Consider the 1D Keller-Segel equation:
\begin{align}\label{Keller-Segel equation}
	\begin{cases}
		$$\rho_{t} = -\nabla\cdot\left(\rho\left(-\nabla\dfrac{\delta }{\delta \rho}\mathcal{J\left(\rho\right)}\right)\right),$$\\
		$$\displaystyle\mathcal{J\left(\rho\right)} = \int_{\Omega}\mathcal{H}\left(\rho\right)\mathrm{d}x + \dfrac{1}{2}\int_{\Omega^{2}} \mathcal{W}\left(x-y\right)\rho\left(x\right)\rho\left(y\right)\mathrm{d}x\,\mathrm{d}y,$$
	\end{cases}
\end{align}
where $\mathcal{H} = \rho\log{\rho}$ and $\mathcal{W}\left(x\right) = \dfrac{1}{2\pi}\ln{\vert x\vert}$ are chosen for a balanced regime \cite{carrillo2015finite}. In this regime, the attraction due to the nonlocal interaction $\mathcal{W}$ and the repulsion due to diffusion $\mathcal{H}$ reaches a delicate balance such that there exists a critical mass separating the diffusive behavior and the blow-up behavior. 

\begin{figure}[htbp]
	\centering  
	\includegraphics[width=0.4\textwidth]{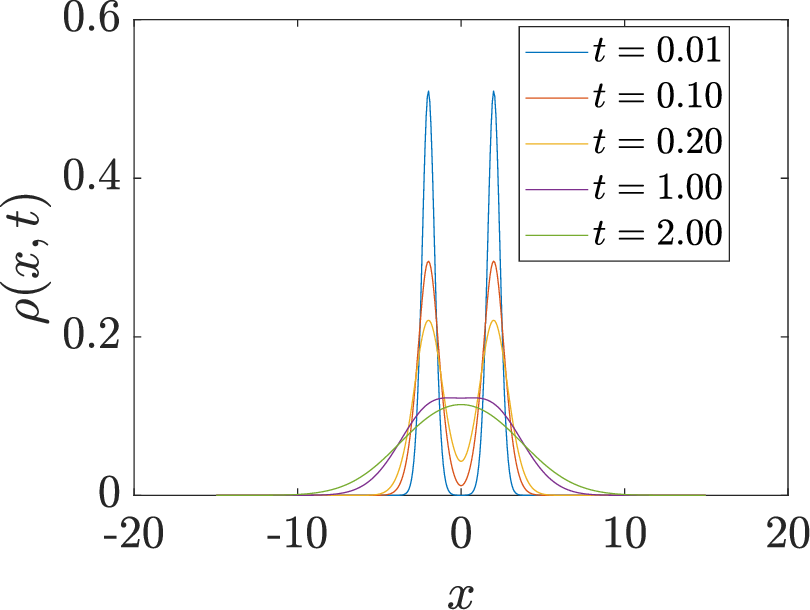}
	\includegraphics[width=0.4\textwidth]{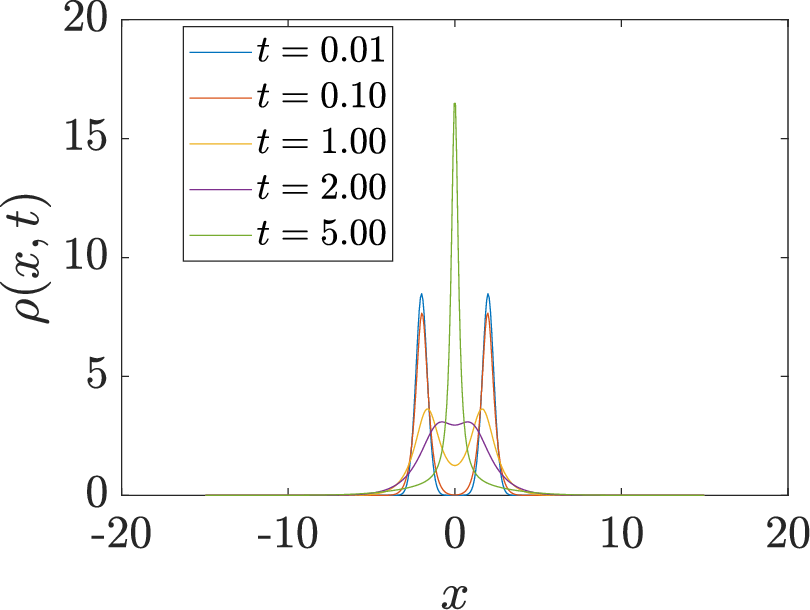}
	\caption{Evolution of solutions of the 1D Keller-Segel equation \eqref{Keller-Segel equation} with different initial values $C = 1$ (left) and $C=15$ (right) for $t \in [0, 2]$. $N_{x} = 800,\; \tau = 0.01$.}
	\label{fig:1D_KS_exp1}
\end{figure}

In this example, we consider the domain $\Omega=[-15,15]$ and the initial value \cite{CHENG2025113696}
\begin{equation*}
	\rho_{0}\left(x\right) = \dfrac{C}{\sqrt{\pi}}\left(
	\mathrm{e}^{-4\left(x-2\right)^{2}} + \mathrm{e}^{-4\left(x+2\right)^{2}}
	\right) + 10^{-8},
	\label{eq:initial_value}
\end{equation*}
where $C$ controls the dynamic behavior.
We solve \eqref{Keller-Segel equation} by the VPTPD method for $C = 1$ and $C = 15$. In particular, to address the singularity of $\mathcal{W}(0)$, we replace it by the approximation $\mathcal{W}(0)\approx \int_{\Omega_{\boldsymbol{i}_0}} \mathcal{W}(x) dx$ integrated over the cell $\Omega_{\boldsymbol{i}_0}$ containing $x=0$ as in \cite{carrillo2022primal}. The results are shown in \autoref{fig:1D_KS_exp1}. When the initial mass is smaller than threshold ($C=1$), the solution decays while aggregating towards $x=0$ and eventually stabilizes due to diffusion; when the initial mass is larger than threshold ($C=15$), the solution will continuously concentrate near $x=0$ and eventually blow up. This numerical experiment exhibits the effectiveness of our structure-preserving VPTPD method in capturing the subtle balance between competing energies. 

\subsubsection{2D aggregation drift equation}
Consider the 2D aggregation-drift equation:
\begin{align}\label{2D aggregation drift equation}
	\begin{cases}
		$$\rho_{t} = -\nabla\cdot\left(\rho\left(-\nabla\dfrac{\delta }{\delta \rho}\mathcal{J\left(\rho\right)}\right)\right),$$\\
		$$\displaystyle\mathcal{J\left(\rho\right)} = \int_{\Omega}\mathcal{V}\left(x\right)\rho\left(x\right)\mathrm{d}x + \dfrac{1}{2}\int_{\Omega^{2}} W\left(x-y\right)\rho\left(x\right)\rho\left(y\right)\,\mathrm{d}x\,\mathrm{d}y.$$
	\end{cases}
\end{align}
where $\mathcal{V}\left(x\right) = -\frac{1}{4}\ln\lvert x\rvert$ provides a drift away from $x=0$ and $W\left(x\right) = \dfrac{x^{2}}{2} - \ln{\vert x\vert}$ provides repulsion at short distances and attraction at longer distances. We consider the domain $\Omega = [-1.5, 1.5]^2$ and the initial value of five Gaussians as follows: 
\begin{align}
	\begin{cases}
		$$ \rho_{0}\left(x,y\right) =10^{-8} + \sum_{i = 1}^{5} k\left(1-\tanh{\dfrac{\sqrt{(x+p_{i})^{2} + (y+q_{i})^{2}} - c_{\mathrm{agg}}}{\sqrt{2}\epsilon}}\right)$$\\
		$$p_{1} = 0,\,p_{2} = p_{4} = 0.5, \,p_{3} = p_{5} = -0.5,$$\\
		$$q_{1} = 0,\, q_{2} = q_{3} = -0.5,\, q_{4} = q_{5} = 0.5,$$\\
		$$k = 0.5,\, c_{\mathrm{agg}} = 0.2,\,\epsilon = 0.1$$.
	\end{cases}
\end{align}

\begin{figure}[htb]
	\centering  
	\includegraphics[width=0.3\textwidth]{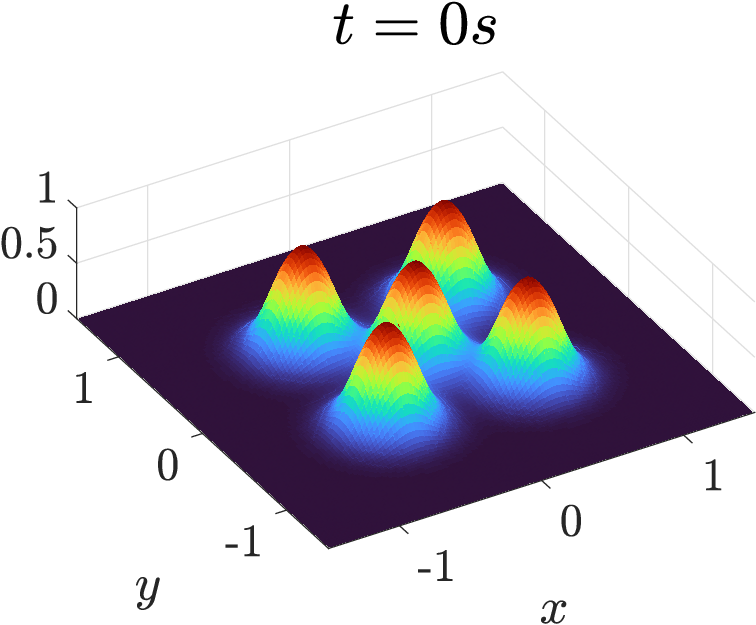}
	\includegraphics[width=0.3\textwidth]{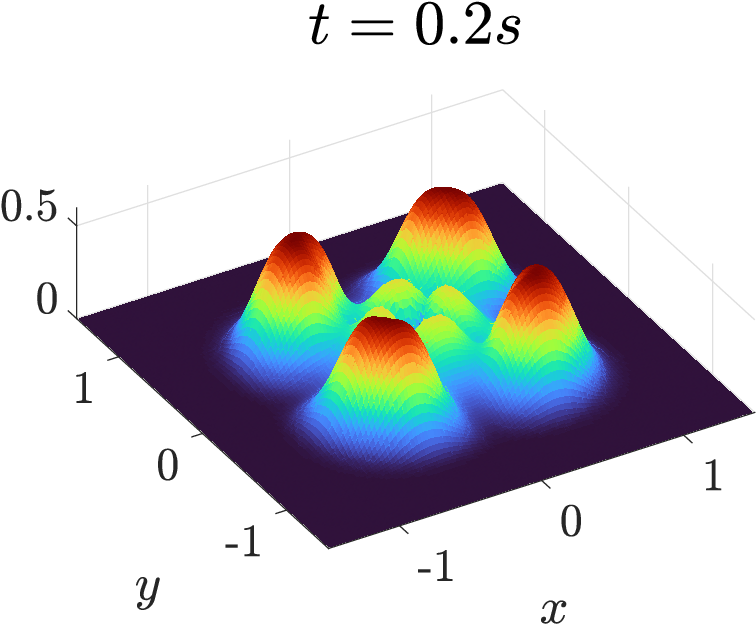}
	\includegraphics[width=0.3\textwidth]{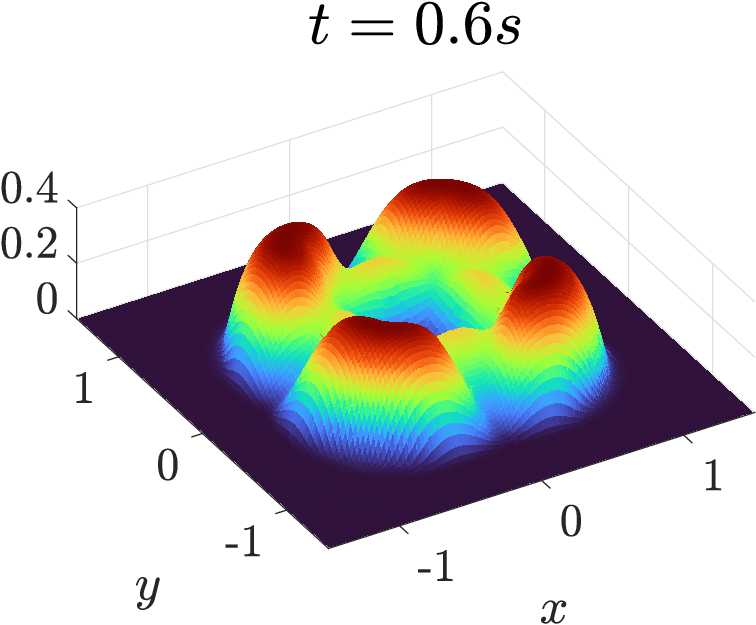}\\
	\includegraphics[width=0.3\textwidth]{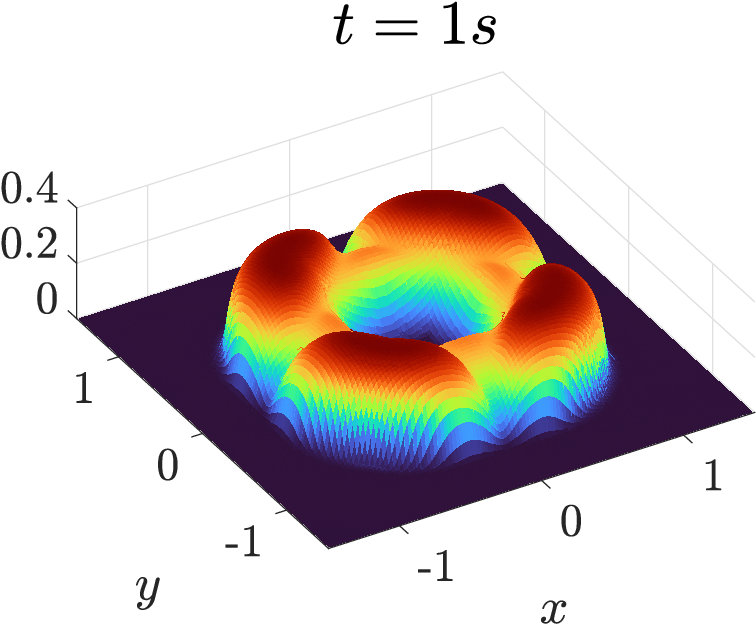}
	\includegraphics[width=0.3\textwidth]{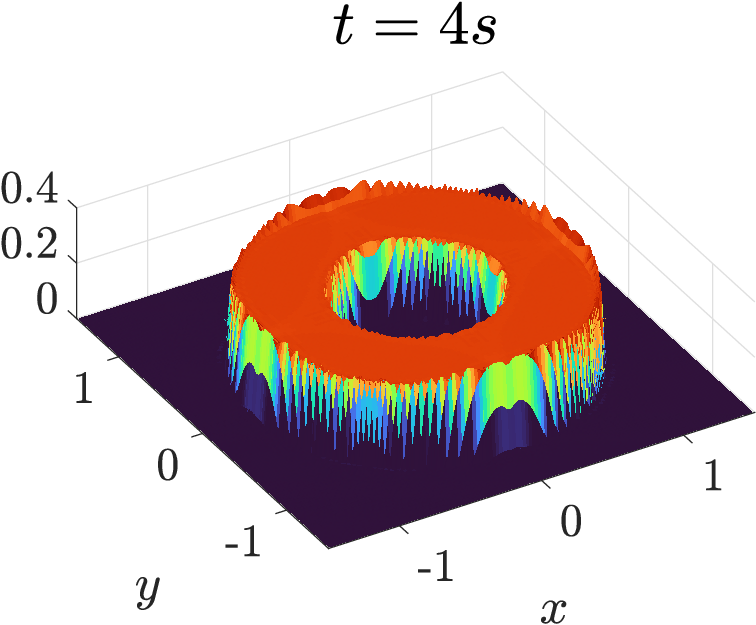}
	\includegraphics[width=0.3\textwidth]{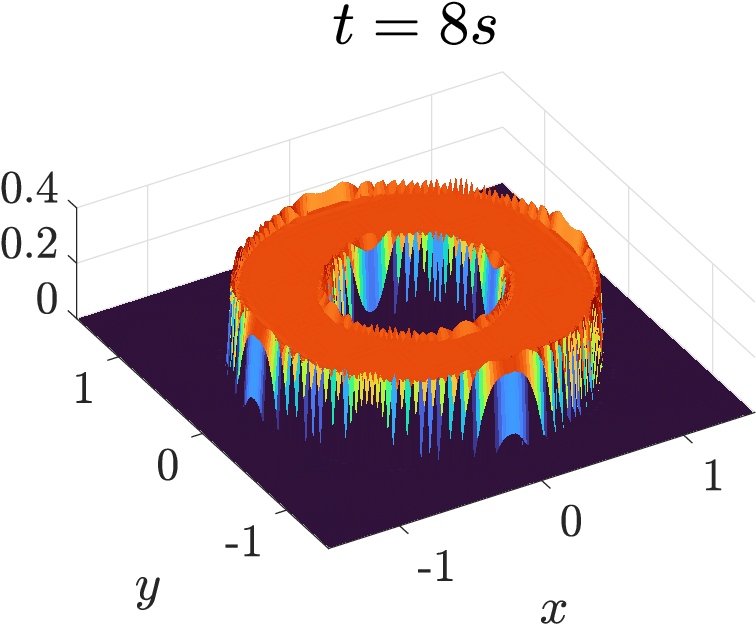}
	\caption{Evolution of solutions for the 2D aggregation drift equation \eqref{2D aggregation drift equation} for $t \in [0, 8]$. $N_{x} = 128,\;N_{y} = 128,\;\tau = 0.2$.}
	\label{fig:2D_KS}
\end{figure}

\begin{figure}[htb]
	\centering  
	\includegraphics[width=0.35\textwidth]{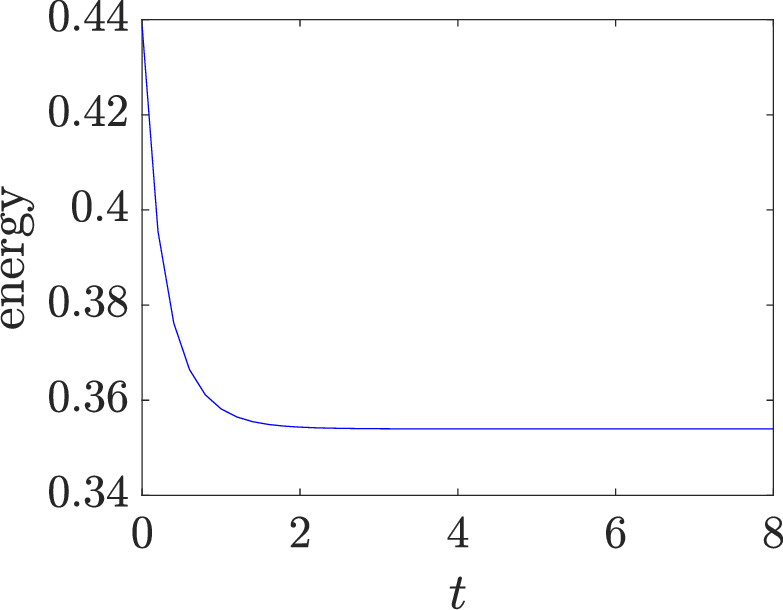}
	\includegraphics[width=0.35\textwidth]{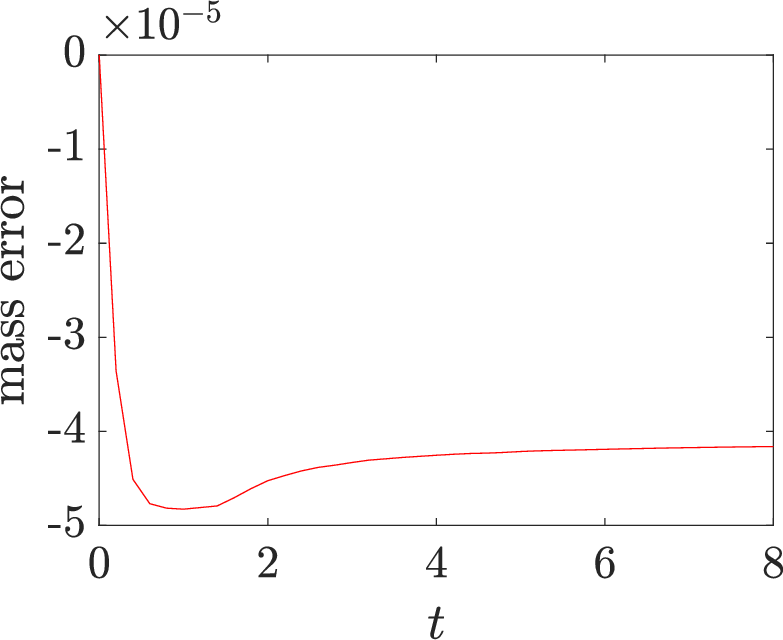}
	\caption{Evolution of energy and relative mass error for the 2D aggregation-drift equation \eqref{2D aggregation drift equation} for $t \in [0, 8]$. $N_{x} = 128,\;N_{y} = 128,\;\tau = 0.2$. }
	\label{fig:2D_KS_E_M}
\end{figure}

We simulate the evolution of solutions for $t\in[0,8]$ by our VPTPD method, and the numerical results for the solution dynamics, energy dissipation and mass conservation are shown in \autoref{fig:2D_KS} and \autoref{fig:2D_KS_E_M}. As expected from analytical and numerical results \cite{CARRILLO_HUANG_MARTIN_2014,carrillo2022primal}, the middle Gaussian splits into parts, drifted away from $x=0$ and merging with the surrounding Gaussians, and the solution eventually reaches the equilibrium profile of ``rotating mill''.

\subsubsection{3D fracture experiment}
Finally, we consider a 3D fracture experiment described by the Cahn-Hilliard equation \eqref{2DCH_random} with nonlinear mobility $M\left(\rho\right)=1-\rho^{2}$. In particular, we simulate the fracture dynamics of a 3D stick driven by diffusion over the domain $\Omega = [-2,2]\times[-0.5,0.5]\times[-0.5,0.5]$
\begin{equation}\label{3D_f_initial}
	\rho_{0}\left(x,y,z\right)=\begin{cases}
		$$1,\quad &\text{if}\;\vert x\vert<1.8\;\text{and}\;\vert y \vert,\vert z\vert < 0.1$$,\\
		$$-1,\quad&\text{otherwise}$$,
	\end{cases}
\end{equation}
with two different energy potential, the double-well potential $\mathcal{H}\left(\rho\right)$:
\begin{equation}
	\mathcal{H}_{\mathrm{dw}}\left(\rho\right)=\dfrac{1}{4}\left(1-\rho^{2}\right)^{2},
\end{equation}
and the logarithmic potential energy $\mathcal{H}_{\mathrm{log}}$:
\begin{equation}
	\mathcal{H}_{\mathrm{log}}\left(\rho\right)=\dfrac{\theta}{2}\left[\left(1+\rho\right)\ln\left(\dfrac{1+\rho}{2}\right)+\left(1-\rho\right)\ln\left(\dfrac{1-\rho}{2}\right)\right] + \dfrac{\theta_{c}}{2}\left(1-\rho^{2}\right),
\end{equation}
where we set $\theta = 0.3$, $\theta_{c}=1$, and $\epsilon=0.02$ for the Dirichelet energy.

\begin{figure}[htbp]
	\centering
	\subfigure{
		\begin{minipage}[b]{0.3\textwidth}
			\centering
			\includegraphics[width=1\linewidth]{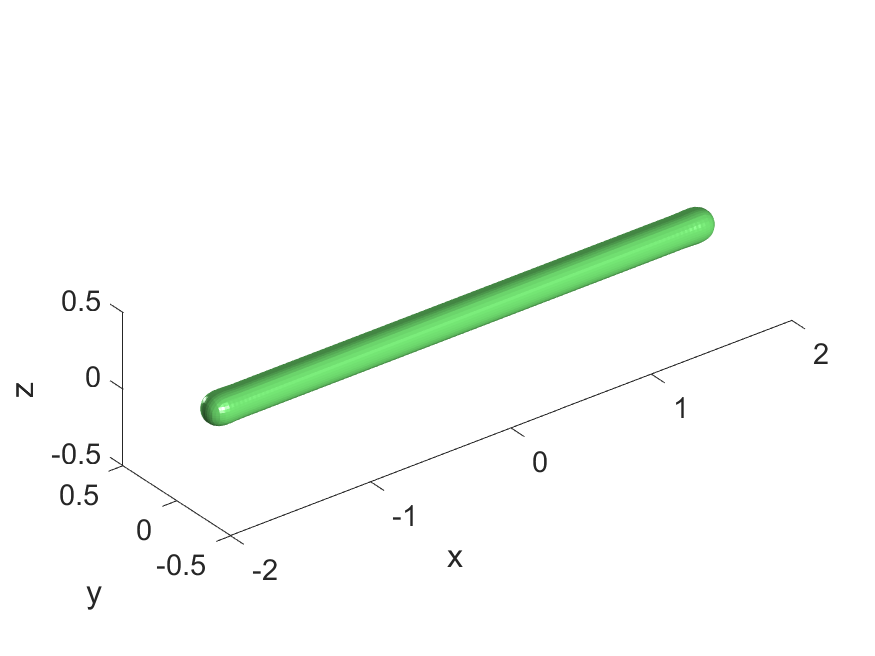} 
		\end{minipage}
	}
	\subfigure{
		\begin{minipage}[b]{0.3\textwidth}
			\centering
			\includegraphics[width=1\linewidth]{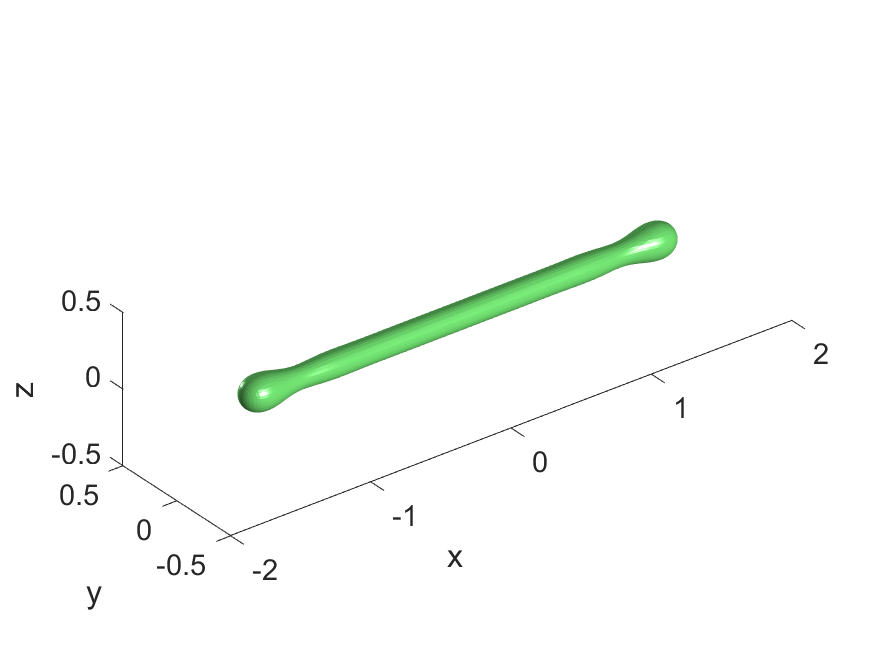} 
		\end{minipage}
	}
	\subfigure{
		\begin{minipage}[b]{0.3\textwidth}
			\centering
			\includegraphics[width=1\linewidth]{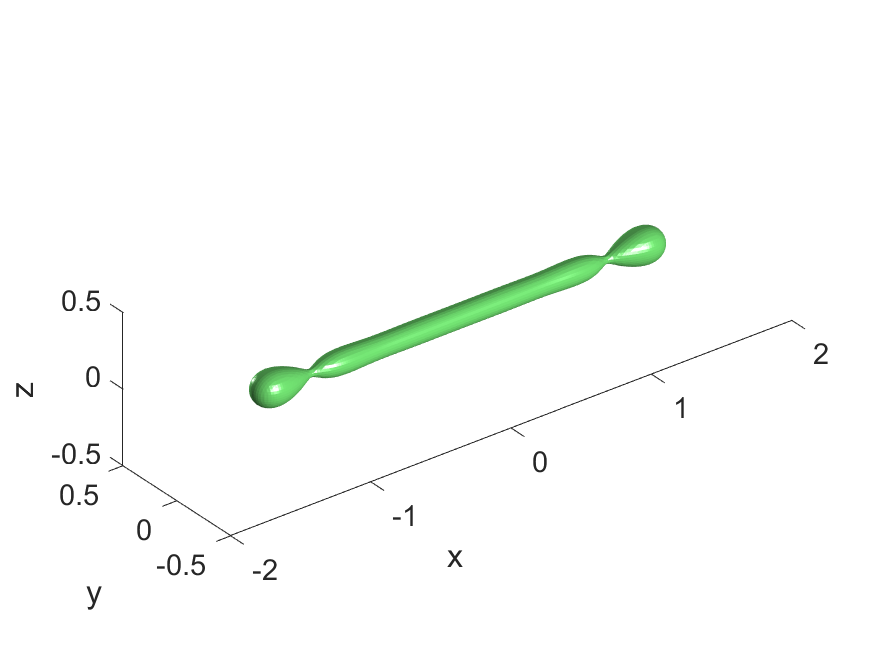} 
		\end{minipage}
	}
	\subfigure{
		\begin{minipage}[b]{0.3\textwidth}
			\centering
			\includegraphics[width=1\linewidth]{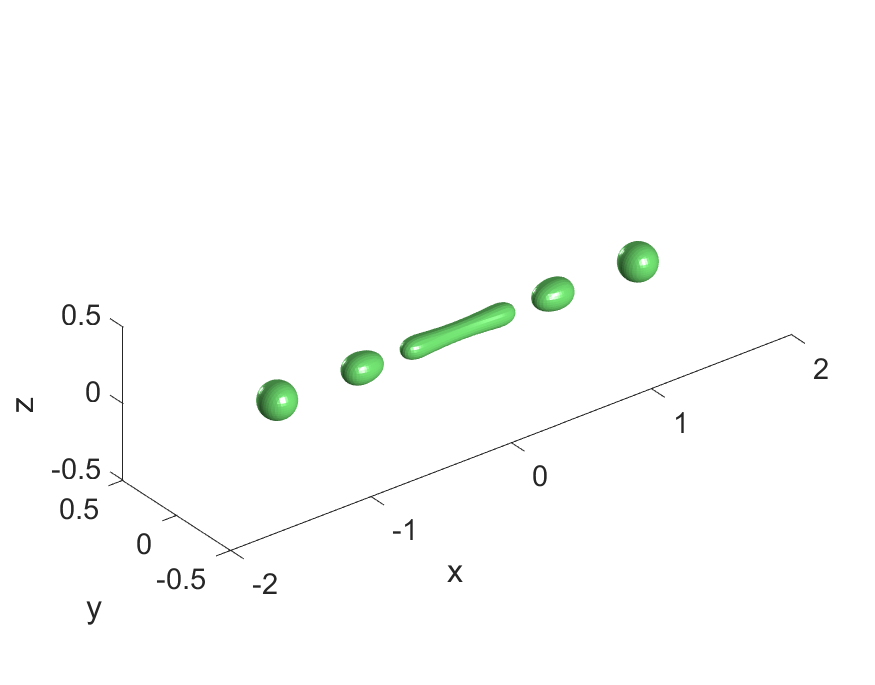} 
		\end{minipage}
	}
	\subfigure{
		\begin{minipage}[b]{0.3\textwidth}
			\centering
			\includegraphics[width=1\linewidth]{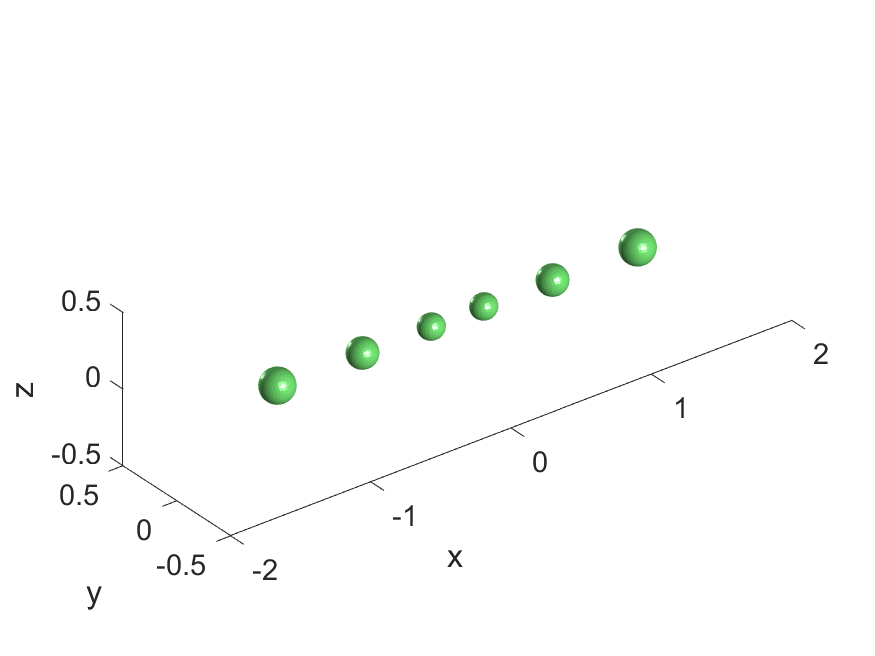} 
		\end{minipage}
	}
	\subfigure{
		\begin{minipage}[b]{0.3\textwidth}
			\centering
			\includegraphics[width=1\linewidth]{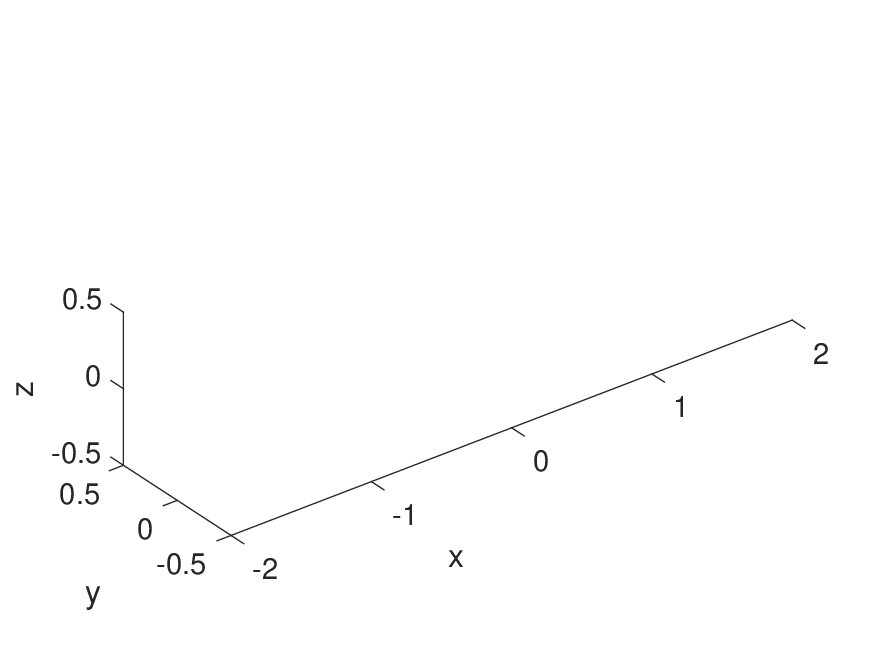} 
		\end{minipage}
	}
	\caption{Evolution of interface in 3D fracture experiment with double-well potential $\mathcal{H}_{\mathrm{dw}}$ for $t \in [0, 2]$. $N_{x} = 256,\; N_{y} = 64,\; N_{z} = 64,\;\tau = 0.001,\;\text{tol}=10^{-6},\;\text{TOL}=10^{-6}$.}
	\label{fig:3D fracture experiment1}
\end{figure}

\begin{figure}[htbp]
	\centering
	\subfigure{
		\begin{minipage}[b]{0.3\textwidth}
			\centering
			\includegraphics[width=1\linewidth]{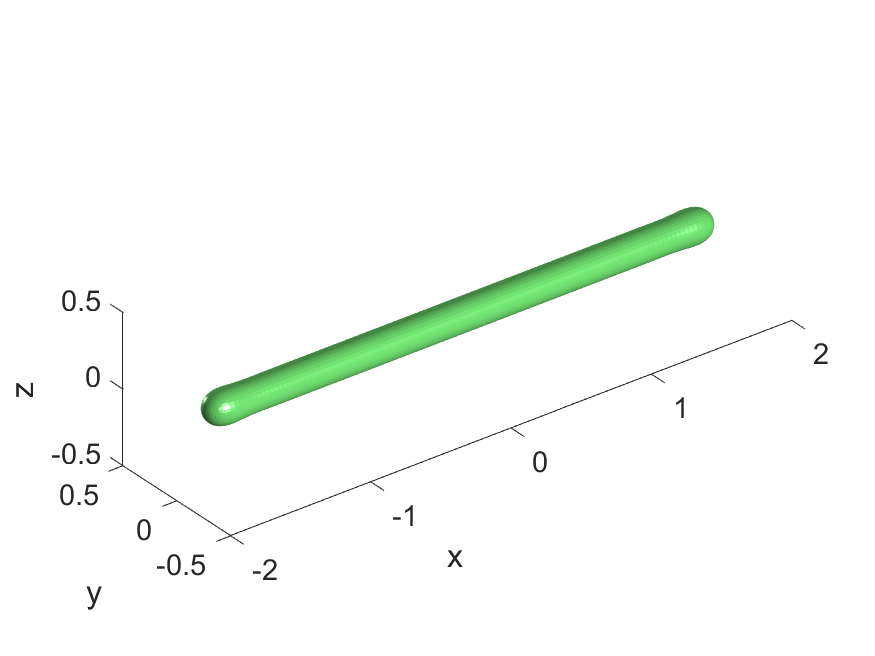} 
		\end{minipage}
	}
	\subfigure{
		\begin{minipage}[b]{0.3\textwidth}
			\centering
			\includegraphics[width=1\linewidth]{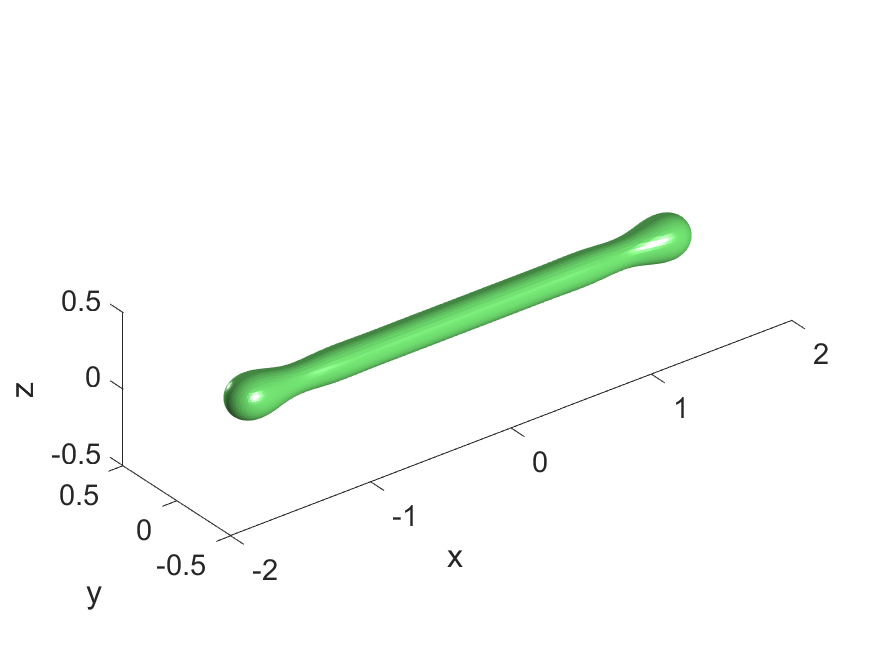} 
		\end{minipage}
	}
	\subfigure{
		\begin{minipage}[b]{0.3\textwidth}
			\centering
			\includegraphics[width=1\linewidth]{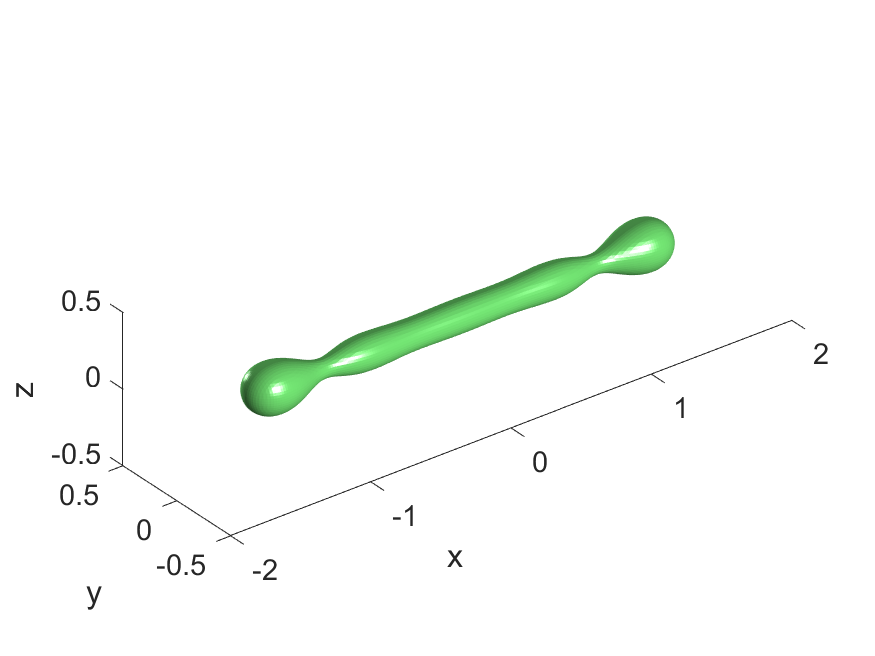} 
		\end{minipage}
	}
	\subfigure{
		\begin{minipage}[b]{0.3\textwidth}
			\centering
			\includegraphics[width=1\linewidth]{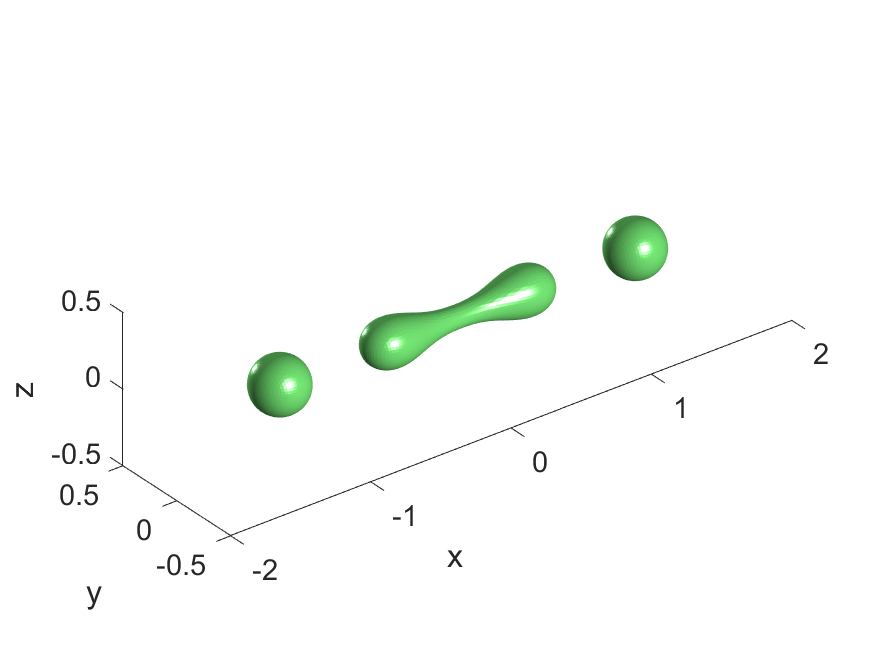} 
		\end{minipage}
	}
	\subfigure{
		\begin{minipage}[b]{0.3\textwidth}
			\centering
			\includegraphics[width=1\linewidth]{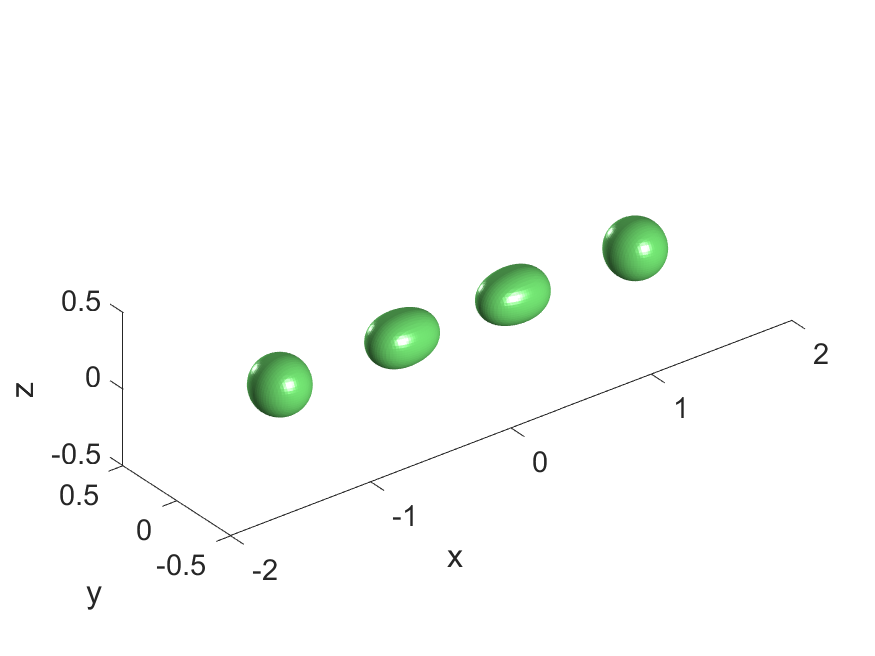} 
		\end{minipage}
	}
	\subfigure{
		\begin{minipage}[b]{0.3\textwidth}
			\centering
			\includegraphics[width=1\linewidth]{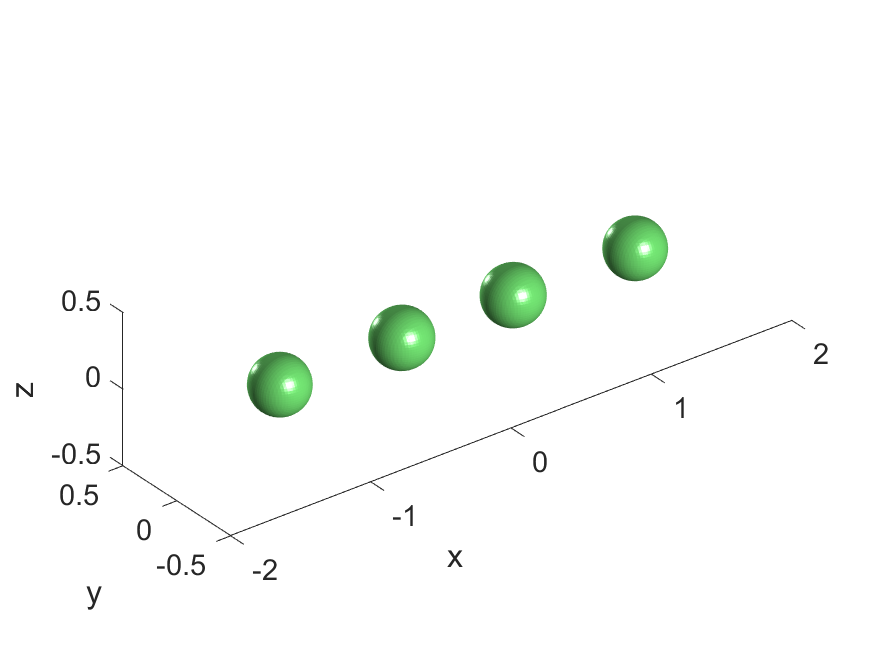} 
		\end{minipage}
	}
	\caption{Evolution of interface in 3D fracture experiment with logarithmic potential $\mathcal{H}_{\mathrm{log}}$ for $t \in [0, 6]$. $N_{x} = 256,\; N_{y} = 64,\; N_{z} = 64,\;\tau = 0.001,\;\text{tol}=10^{-6},\;\text{TOL}=10^{-6}$.}
	\label{fig:3D fracture experiment2}
\end{figure}

\begin{figure}[htbp]
	\centering  
	\includegraphics[width=0.35\textwidth]{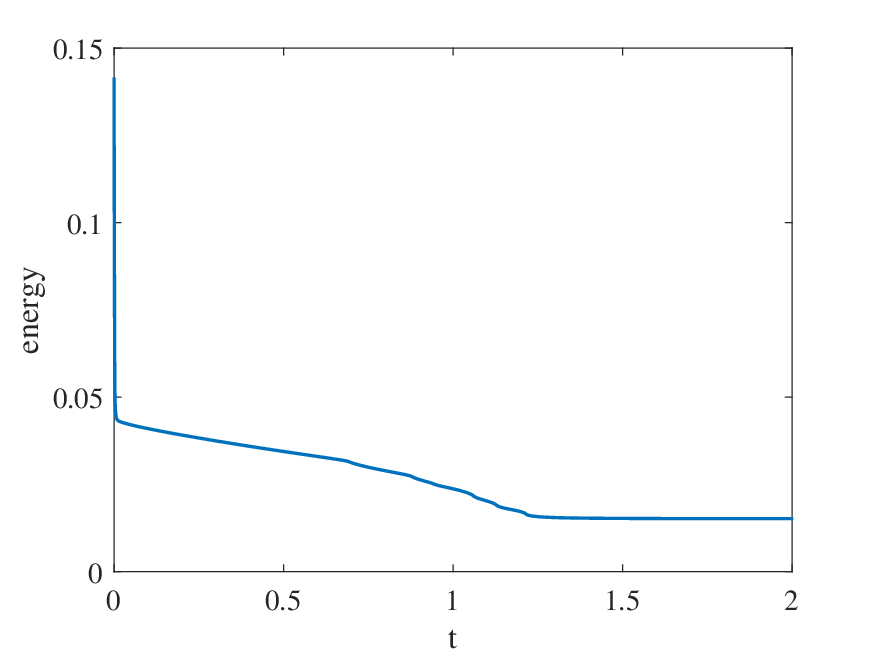}
	\includegraphics[width=0.35\textwidth]{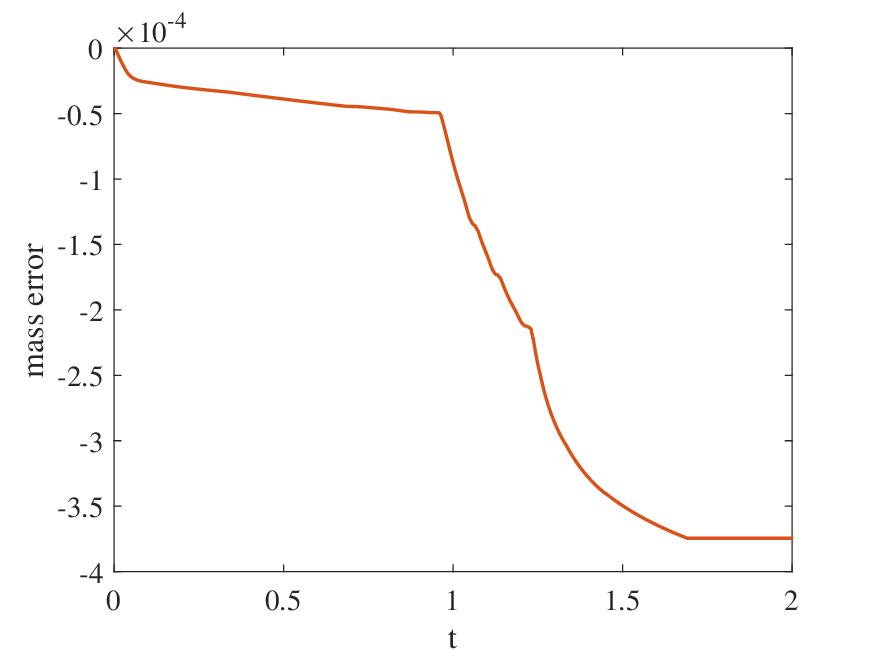}
	\caption{Evolution of energy and relative mass error of the 3D fracture experiment with double-well potential.}
	\label{fig:3D fracture EM}
\end{figure}

The high efficiency of the VPTPD method enables us to perform the challenging 3D simulations and investigate the applicability of the phase-field model with different energy potentials for interface dynamics driven by surface diffusion. As shown in \autoref{fig:3D fracture experiment1}, the numerical results for the model with double-well potential $\mathcal{H}_{\mathrm{dw}}$ and quadratic mobility exhibit a pronounced spontaneous shrinkage effect due to undesired diffusion \cite{Lee2015Degenerate}: the interface gradually disappears and the enclosed domain eventually vanishes when its radius decreases below a critical size \cite{zhou2025newphasefieldmodelanisotropic}. This entire process of energy dissipation and mass conservation is shown in \autoref{fig:3D fracture EM}, which indicates that this "unphysical" behavior is not due to numerical mass loss or instability but due to the model itself. In contrast, as shown in \autoref{fig:3D fracture experiment2}, the model with the logarithmic potential  $\mathcal{H}_{\mathrm{log}}$ and quadratic mobility accurately describes the interface dynamics driven by surface diffusion \cite{Elliott1996CahnHilliard,Lee2015Degenerate,huang2023structure}, capable of capturing the evolution of small or thin interfaces and preserving the enclosed volume more accurately.  
\FloatBarrier


\section*{Declaration of competing interest}
The authors declare that they have no known competing financial interests or personal relationships that could have appeared to influence the work reported in this paper.

\section*{Data availability}
Data will be made available on request.

\section*{Acknowledgments}
Chaozhen Wei is supported by the National Natural Science Foundation of China under grants 12371392 and 12431015.
Ruchi Guo is supported by National Key R\&D Program of China 2025YFA1018700, NSFC grant 12571436, 
the Fundamental Research Funds for the Central Universities, NSF grant DMS-2012465.

\bibliographystyle{elsarticle-num.bst}
\bibliography{references.bib}

\appendix
\section{Temporal and spatial accuracy of the fully discrete JKO scheme}\label{JKO_accuracy}
To demonstrate the order of spatial and temporal accuracy of the JKO scheme with VPTPD method, we have selected two test examples from \cite{Fu2023JCP} and carried out similar numerical tests. 

\noindent Example 1: 3D Heat Equation.
$M(\rho)=\rho$, \quad $\Omega = [0,1]^{3}$, \quad $T = 0.1$

\[
\left\{
\begin{array}{l}
\displaystyle
\frac{\partial \rho}{\partial t} = \Delta \rho \\
\rho_{0} = 1 + \frac{1}{2} \cos (\pi x) \cos (\pi y) \cos (\pi z), \quad t = 0
\end{array}
\right.
\]

\noindent
whose exact solution is
\[
\rho(\boldsymbol{x},t) = 1 + \frac{1}{2} \cos (\pi x) \cos (\pi y) \cos (\pi z) e^{-3\pi^{2} t}.
\]

\noindent
Setting the spatial mesh size to $\Delta x = \Delta y = \Delta z = 1/32$, we compute the $L^{2}$-error between the numerical solution and the exact solution at the final time $T=0.1$ for various temporal step sizes $\tau$. The resulting convergence order is summarized in the following table:

\begin{table}[htbp]
\centering
\caption{Temporal convergence rate of VPTPD-JKO for the $L^{2}$-error at $T=0.1$ for 3D heat equation.}
\begin{tabular}{ccc}
\toprule
$\tau$ & Error & Order \\
\midrule
1/4   & $1.04 \times 10^{-2}$ & -- \\
1/8  & $5.42 \times 10^{-3}$ & 0.98 \\
1/16 & $2.87 \times 10^{-3}$ & 0.98 \\
1/32 & $1.61 \times 10^{-3}$ & 0.95 \\
\bottomrule
\end{tabular}
\end{table}

\noindent Example 2: 2D Nonlinear Fokker-Planck Equation. $M(\rho) = \rho$, \quad $\Omega = [-1, 1]^{2}$, \quad $T = 10$

\[
\left\{
\begin{array}{l}
\displaystyle
\frac{\partial \rho}{\partial t} - \Delta \rho^{3} = \nabla \cdot (\rho \boldsymbol{x}), \quad t \in (0, T] \\[6pt]
\rho_{0} = C, \quad t = 0
\end{array}
\right.
\]

\noindent
This problem reaches a steady state solution:
\[
\rho_{\text{steady}}(\boldsymbol{x}) = \sqrt{\frac{\left(2C - (x^{2} + y^{2})\right)_{+}}{3}}.
\]

\noindent
We fix the time step size and compute the $L^{2}$-error between the numerical solution at $T = 10$ and the steady state solution for different spatial mesh sizes. The convergence order with respect to the grid size is summarized in the following table:

\begin{table}[htbp]
\centering
\caption{Temporal convergence rate of VPTPD-JKO for the $L^{2}$-error in steady state ($T=10$) for 2D Nonlinear Fokker-Planck Equation.}
\begin{tabular}{ccc}
\toprule
$\Delta x = \Delta y$ & Error & Order \\
\midrule
$1/8$   & $3.70 \times 10^{-3}$ & -- \\
$1/16$  & $9.27 \times 10^{-4}$ & 1.9972 \\
$1/32$  & $2.32 \times 10^{-4}$ & 1.9981 \\
$1/64$  & $5.81 \times 10^{-5}$ & 1.9975 \\
\bottomrule
\end{tabular}
\end{table}

\section{Equivalence between G-prox PDHG method and PrePD method}\label{Gprox-PrePD}
The one-step JKO optimization problem \eqref{generalized discrete JKO} can be recast as a different saddle-point system:
\begin{equation}\label{saddle}
\min_{u=(\rho,\mathbf{m})} \max_{p} \left[ \mathcal{F}_h(u) + \langle (\rho - \rho^k) + \nabla_h \cdot \mathbf{m}, p \rangle \right] \overset{\text{integration by parts}}{\implies} -\min_{p} \max_{u} \, \langle Kp, u\rangle + g^*(p) - \mathcal{F}_h(u) ,
\end{equation}
where 
\begin{equation*}
\mathcal{F}_h(u) = \Phi_h(u)/2\tau + \mathcal{J}_h(u), \quad
g^*(p) = \langle p, \rho^k\rangle, \quad
K = (-I; \nabla_h).
\end{equation*}
The iterative scheme of G-prox PDHG \cite[Algorithm 2]{Jacobs2019Gprox} for the saddle-point system above is:
\begin{align}\label{gprox}
	\begin{cases}
    & p_{n+1} = \arg\min \limits_{p} \, g^*(p) + \langle Kp, \bar{u}_n \rangle + \frac{1}{2\sigma} \|K(p - p_n)\|^2,\\
    & u_{n+1} = \arg\min \limits_{u} \, \mathcal{F}_h(u) - \langle Kp_{n+1}, u \rangle + \frac{1}{2\lambda} \|u - u_n\|^2 ,\\
    & \bar{u}_{n+1}=u_{n+1}+(u_{n+1}-u_n),
    \end{cases}
\end{align}

The PrePDJKO \cite[Algorithm 2]{carrillo2024structure} treats the one-step JKO optimization as another saddle-point problem:
\begin{equation}\label{saddle2}
\min_{u} \max_{p} \, \langle Bu, p\rangle +\mathcal{F}_h(u) - g^*_\delta(p) ,
\end{equation}
where $g^*_\delta(p)$ corresponds to a relaxed penalty term $g(Bu)=i_\delta(Bu)$ for the constraint of continuity equation:
\begin{equation*}
B = (I, \nabla_h\cdot), \quad 
g^*_\delta(p)=\delta \|p\| + \langle \rho^k, p \rangle.
\end{equation*}
The iterative scheme of PrePDJKO for the saddle point \eqref{saddle2} gives:
\begin{align}\label{PrePD}
	\begin{cases}
		& p_{n+1} = \arg \min \limits_{p} \ g^{*}_\delta\left(p\right) - \langle B\bar{u}_{n},p \rangle + \dfrac{1}{2\sigma}\left\Vert p-p_{n}\right\Vert_{BB^\top}^{2}, \\
        & u_{n+1} = \arg \min \limits_{u} \  \mathcal{F}_h(u) + \langle u,B^{\top} p_{n+1} \rangle + \dfrac{1}{2\lambda}\left\Vert u-u_{n}\right\Vert^{2},\\& \bar{u}_{n+1}=u_{n+1}+(u_{n+1}-u_n).
	\end{cases}
\end{align}
By recognizing $B=(I, \nabla_h\cdot)=-K^\top$,  $B^T=(I; -\nabla_h)=-K$ and hence $BB^\top=K^\top K$, the PrePDJKO iterative scheme \eqref{PrePD} and the saddle point \eqref{saddle2} (with relaxation $\delta=0$) are essentially {\it identical} to the G-prox PDHG iterative scheme \eqref{gprox} and corresponding saddle-point problem \eqref{saddle}. 

However, the computation of proximal operator for general $\mathcal{F}_h$ is not a trivial task, especially when the energy functional $\mathcal{J}_h$ involves Dirichelet energy and nonlocal interaction energy.
Since $\mathcal{F}_h$ consists of smooth part $\mathcal{J}_h$ and nonsmooth part $\Phi_h$, we instead apply proximal splitting descent in the primal step for general Wasserstein gradient flows:
\begin{align}\label{PrePDJKO}
	\begin{cases}
        & p_{n+1} = \arg \min \limits_{p} \ g^{*}_\delta\left(p\right) - \langle B\bar{u}_{n},p \rangle + \dfrac{1}{2\sigma}\left\Vert p-p_{n}\right\Vert_{BB^\top}^{2}, \\
		& u_{n+1} = \arg \min \limits_{u} \Phi_h(u) + \mathcal{J}_h(u^n) + \langle u , \nabla \mathcal{J}_h(u^n)\rangle + \langle u,B^{\top} \bar{p}_{n} \rangle + \dfrac{1}{2\lambda}\left\Vert u-u_{n}\right\Vert^{2}, \\
        & \bar{u}_{n+1} = u_{n+1} + \left[ (u_{n+1}-\lambda\nabla\mathcal{J}_h(u_{n+1}))- (u_{n} - \lambda\nabla\mathcal{J}_h(u_n)) \right].
	\end{cases}
\end{align}
Therefore, PrePDJKO method \eqref{PrePDJKO} can be regarded as a proximal-splitting variant of G-prox PDHG method \eqref{gprox} for three-operator problems, such as one-step JKO minimization. 

\section{Implementation of G-prox PDHG}\label{Gprox_prox}
To implement G-prox PDHG \eqref{gprox} for JKO scheme, one has to compute the proximal operator of $\mathcal{F}_h=\Phi_h/2\tau+\mathcal{J}_h$ in the primal step, rather than that of $\Phi_h$ for PrePDJKO \eqref{PrePDJKO}. For 1D porous medium equation, the corresponding proximal operator of $\mathcal{F}$ is given by:
\begin{equation*}
\operatorname{Prox}_{\lambda \mathcal{F}_h}(\hat{\rho},\hat{\mathbf{m}}) = \arg\min \limits_{\rho,\mathbf{m}} \  
 \dfrac{1}{2}\left\Vert \rho-\hat{\rho} \right\Vert^{2} + \dfrac{1}{2}\left\Vert \mathbf{m}-\hat{\mathbf{m}} \right\Vert^{2} + \sum \limits_{i} \Big(\dfrac{\zeta}{2}\dfrac{\|\mathbf{m}_i\|^2}{\rho_i}+\xi \rho_i^2 \Big) +\lambda \big(\langle \rho, p_{n+1} \rangle - \langle \mathbf{m}, \nabla_h p_{n+1} \rangle\big),
\end{equation*}
where $\zeta = \lambda \Delta V/\tau$ and $\xi = \lambda \Delta V$. Fortunately, this proximal operators is separable and can be efficiently computed. In fact, one can verify that the proximal operator reduces to finding the largest positive root of a cubic polynomial:
\begin{align*}
	\big(\operatorname{Prox}_{\lambda \mathcal{F}_h}(\hat{\rho},\hat{\mathbf{m}}) \big)_i = \begin{cases}
		\displaystyle 
        \Big( \rho^*_i, \frac{\rho^*_i}{\rho^*_i+\zeta} \big(\hat{\mathbf{m}}_i + \lambda (\nabla_h p_{n+1})_i \big) \Big) ,& \mathrm{if}\;  \hat{\rho}_i \geq \lambda p_{n+1}-\dfrac{1}{2\zeta} \|\hat{\mathbf{m}}_i + \lambda (\nabla_h p_{n+1})_i \|^{2},\\
		(0,\mathbf{0}),& \mathrm{otherwise},
	\end{cases}
\end{align*}
where $\rho_i^*$ is the only positive root of the cubic polynomial $L(\rho)$ given by:
\[
L(\rho) = [(2\xi+1)\rho + \lambda (p_{n+1})_i - \hat{\rho}_i]\cdot(\rho+\zeta)^2 - \frac{\zeta}{2}(\hat{\mathbf{m}}_i + \lambda (\nabla_h p_{n+1})_i)^2.
\]
However, in general, the tractability of computing $\mathrm{prox}_{\lambda\mathcal{F}_h}$ depends on the specific form of the energy functional $\mathcal{J}_h$. Designing an efficient solver for this proximal operator is a nontrivial task, particularly when $\mathcal{J}_h$ involves differential operators or convolutions as in the case of Ginzburg-Landau free energy or nonlocal interaction energy. In such cases, one typically needs to solve a highly nonlinear coupled system, with particular care taken to ensure the convergence of the iterative solver and the boundedness of the computed solution.

Although numerical experiments in \autoref{fig:performance_comparison} indicate mesh-independent convergence of G-prox PDHG (as well as PrePDJKO) for the porous medium equation, we noticed that this example does not seem to satisfy a key assumption in the mesh-independent convergence analysis in \cite{Jacobs2019Gprox}:
\begin{equation}\label{supnorm}
\sup_{p} \lVert u(p) \rVert < C < \infty, \quad \text{where } u(p) \in \arg\max_{u} \bigl( \langle Kp, u \rangle - \mathcal{F}(u) \bigr).
\end{equation}
This condition demands that the maximizer $u(p)$ to remain uniformly bounded in norm by a constant $C$, regardless of the operator $K$ and its discretization on grids. 

However, one can obtain the explicit expression of the maximazer $u(p)$ for porous medium equation:
\begin{equation}
    \label{u_p}
    u(p) = \arg\max\limits_{\rho\geq0,m} \left\{ -p\cdot\rho + \nabla p \cdot \mathbf{m} -\frac{\lvert \mathbf{m} \rvert^2}{2\tau\rho} - \rho^2 \right\}=\left(\frac{1}{2}\left[\frac{\tau}{2} \lvert\nabla p\rvert^2 - p \right]_{+}, \tau \rho \nabla p \right)
\end{equation}
whose supremum norm $\sup_p \lVert u(p) \rVert$ depends on $|\nabla p|$, which is unbounded at the continuous level and whose discrete norm generally grows under mesh refinement, thereby violating assumption \eqref{supnorm}. 
This indicates a gap between the current theory and the observed numerical behavior. The mesh-dependence behavior of G-prox PDHG and VPTPD is an important and subtle issue that deserves further investigation. A complete understanding of this question is highly nontrivial and lies beyond the scope of the present paper. We therefore leave it to a dedicated future study.

\end{document}